\documentclass[aos,seceqn,MSNbibl,citesort,dvips]{arximspdf}
\usepackage{dcolumn}
\usepackage{graphicx}

\doi{10.1214/10-AOS863}
\volume{39}
\issue{2}
\pubyear{2011}
\firstpage{1125}
\lastpage{1179}

\makeatletter

\def\bsuffix #1{#1}

\def\bptnote#1{}

\def\@bmisc[#1]{%
  \get@battribute{unstr}%
  \common@pub@types%
  \let\bauthor\bbl@bauthor%
  \let\bhowpublished\@firstofone%
  \def\borganization##1{{\bauthor@style ##1}}%
}

\newcolumntype{d}[1]{D{.}{.}{#1}}

\newproclaim{conditionE1}{E.1. \textit{Random sampling}}
\newproclaim{conditionE2}{E.2. \textit{Regular partial order}}
\newproclaim{conditionE3}{E.3. \textit{Identification condition}}
\newproclaim{conditionE4}{E.4. \textit{Continuity of partial quantile points}}
\newproclaim{conditionE5}{E.5. \textit{Empirical error of probability of comparisons}}
\newproclaim{conditionE6}{E.6. \textit{Gaussian process in $\mathcal{T}$}}

\newtheorem{theorem}{Theorem}
\newtheorem{corollary}{Corollary}
\newtheorem{lemma}{Lemma}
\newtheorem{proposition}{Proposition}

\newproclaim{definition}{Definition}
\newproclaim{remark}{Comment}[section]
\newproclaim{example}{Example}

\renewcommand{\hat}{\widehat}

\newcommand{\Pn}{\mathbb{P}_n}
\newcommand{\Pp}{P}
\newcommand{\mS}{\mathcal{S}}
\newcommand{\RR}{\mathbb{R}}

\makeatother

\begin{document}
\begin{frontmatter}

\title{On multivariate quantiles under partial orders}
\runtitle{Partial quantiles}

\begin{aug}
\author[A]{\fnms{Alexandre} \snm{Belloni}\corref{}\ead[label=e1]{abn5@duke.edu}} and
\author[A]{\fnms{Robert L.} \snm{Winkler}\ead[label=e2]{rwinkler@duke.edu}}
\runauthor{A. Belloni and R. L. Winkler}
\affiliation{Duke University}
\address[A]{Fuqua School of Business\\
Duke University\\
1 Towerview Drive\\
Durham, North Carolina 27708-0120\\
PO Box 90120\\
USA\\
\printead{e1}\\
\phantom{E-mail: }\printead*{e2}}
\end{aug}

\received{\smonth{7} \syear{2010}}

%
\begin{abstract}
This paper focuses on generalizing quantiles from the ordering point of
view. We propose the concept of \textit{partial quantiles}, which are
based on a given partial order. We establish that partial quantiles
are equivariant under order-preserving transformations of the data,
robust to outliers, characterize the probability distribution if the
partial order is sufficiently rich, generalize the concept of efficient
frontier, and can measure dispersion from the partial order
perspective.

We also study several statistical aspects of partial quantiles. We
provide estimators, associated rates of convergence, and asymptotic
distributions that hold uniformly over a continuum of quantile indices.
Furthermore, we provide procedures that can restore monotonicity
properties that might have been disturbed by estimation error,
establish computational complexity bounds, and point out a
concentration of measure phenomenon (the latter under independence and
the componentwise natural order).

Finally, we illustrate the concepts by discussing several theoretical
examples and simulations. Empirical applications to compare intake
nutrients within diets, to evaluate the performance of investment
funds, and to study the impact of policies on tobacco awareness are
also presented to illustrate the concepts and their use.
\end{abstract}

%

\end{frontmatter}

\section{Introduction}\label{Sec:Intro}

The quantiles of a univariate random variable have proved to be a
valuable tool in statistics. They provide important notions of location
and scale, exhibit robustness to outliers, and completely characterize
the random variable. Moreover, quantiles also play a significant role
in applications. Naturally, the quantiles of a
multivariate random variable are also of interest, and the search for a
multidimensional counterpart of the quantiles of a random variable has
attracted considerable attention in the statistical literature. Various
definitions have been proposed and studied.

Barnett \cite{Barnett1976}, Serfling \cite{Serfling2002} and Koenker
\cite{K2005} provide valuable
comparisons and surveys of different methods. Some interesting recent
work is presented in Hallin, Paindaveine and Siman \cite{HPS2009}
(with discussions \cite
{HPS2009Rejoinder,SerflingZuoDisscussion2009,WeiDiscussion2009}), Kong
and Mizera \cite{KongMizera2008} and Serfling \cite{Serfling2009}. A
substantial part of the literature focuses on developing relevant
measures to characterize location and scale information of the
multivariate random variable of interest. This is usually accomplished
by defining a suitable nested family of sets. As discussed below, our
focus will be on a given partial order between points instead. The
incorporation of this additional information is the distinctive feature
of this work. Therefore, our approach is different and hence
complementary to previous work that focuses on location and scale measures.

The fundamental difficulty in reaching agreement on a suitable
generalization of univariate quantiles is arguably the lack of a
natural ordering in a multidimensional setting. Serfling
\cite{Serfling2002} points out that, as a result, ``various
\textit{ad hoc} quantile-type multivariate methods have been
formulated, some vector-valued in character, some univariate, and the
term ``quantile'' has acquired rather loose usage'' (page~214). The simplest
notion of a multivariate quantile is that of a vector of the
corresponding univariate quantiles, but this fails to reflect any
multivariate features of the random vector. More often than not,
attempts to take into account such multivariate features have been
influenced by the justifiable temptation to exploit some geometric
structure of the underlying space. For example, many approaches are
based on the use of specific metrics to collapse the multivariate
setting into a univariate measure.
Many definitions of
multivariate quantiles that use notions such as the distance from a
central measure, norm minimization, or gradients immediately make the
values relevant. In contrast, for univariate quantiles only the
ordering matters, and the actual values of the variable away from the
quantile of interest are irrelevant.

In our work, within the definition of multivariate quantiles, the crux
is the concept of ordering, which might or not be related to geometric
notions of the underlying space. Our starting point will be to
detach our concept from the geometry of the random variable, and assume
that a partial order is
provided which will be used to define the \textit{partial quantiles}. This
allows our work to focus on the minimum structure for which
the problem makes sense. With a~general partial order, as
opposed to a complete order, we recog\-nize that some points simply
cannot be compared. Our key insight is to rely on a~family of
conditional probabilities induced by the partial order to circumvent
the lack of comparability. Such approach yields a distinguishing
feature of the proposed \textit{partial quantiles}: the reliance on
the partial order. Our analysis is close in spirit to, but still quite
different from, the important work of Einmahl and Mason
\cite{EinmahlMason1992}, who proposed a broad class of generalized quantile
processes. We
defer a detailed discussion to Section \ref{sec:add} but we anticipate
that our definitions do not fit within the framework of
\cite{EinmahlMason1992} and most of our results have no parallels in
\cite{EinmahlMason1992}.

Our main contributions are as follows. First, we propose a
generalization of quantiles based on a given partial order on the space
of values of the random variable of interest. Index, point, surface,
and comparability notions of the partial quantiles are studied. We
establish that these partial quantiles have several desirable features:
equivariance under monotone mappings with respect to the chosen partial
ordering (an instrumental feature of the univariate case);
generalization of the efficient frontier concept; meaningfulness not
only in high-dimensional Euclidean spaces but also in arbitrary sets
(relevant for decision making, where metrics are not available); and
applicability even to general binary relations.

Second, we investigate statistical estimation and inference based on
finite samples. We derive results on rates of convergence that hold
uniformly over infinitely many quantile indices. In the analysis of the
estimation problems, we have to accommodate discontinuous criterion
functions, potential nonuniqueness of the true parameter, and a
restricted identification condition. These difficulties lead to
nonstandard rates of convergence. Also, we derive the asymptotic
distribution for the partial
quantile indices process (indexed by a subset of the underlying space)
and for the partial quantile comparability where non-Gaussian limits
are possible due to nonuniqueness.

Several other results are established. Partial quantile indices and
probabilities of comparisons are robust to outliers and we study when
they characterize the underlying probability distribution, both
important properties of univariate quantiles. Due to sampling error,
the estimated partial quantile points could violate the partial order,
as can happen with (univariate) quantile regression \cite{K2005}. In
quantile regression, Chernozhukov, Fern\'{a}ndez-Val and Galichon
\cite{CFG2007,CFG2009} based on rearrangement,
Dette amd Volgushev based on smoothing and monotonization~\cite{DetVol08},
and Neocleous and
Portnoy \cite{NeocleousPortnoy2008} based on interpolation, show how
to obtain monotone estimates of quantile curves. In the context of
partial quantiles within lattice spaces, we propose a new procedure to
correct for this estimation error that leads to partial quantile point
estimates that are monotone with respect to the partial order. (Under
the componentwise natural ordering, we build upon the use of
rearrangement in Chernozhukov, Fern\'{a}ndez-Val and Galichon \cite
{CFG2007,CFG2009} to achieve an improvement on the estimation under
suitable mild conditions.) Under independence and the componentwise
natural ordering, we also point out a concentration of measure and a
possible ``curse of dimensionality'' for comparisons. We also define
dispersion measures based on partial quantile regions. Moreover, we
study the computational requirements associated with approximating
partial quantiles. We provide interesting primitive conditions
under which computation can be carried out efficiently. Finally, we
illustrate these concepts through applications to evaluate the intake
of nutrients within diets, the performance of investment funds, and the
impact of different policies on tobacco awareness.

\section{Partial quantiles}\label{sec:pq}

In this section, we propose a generalization of quantiles and derive
basic probabilistic properties implied by the definition of partial quantiles.

\subsection{Definitions}

Let $X$ denote an $\mathcal{S}$-valued random variable defined on a
probability space $(\Omega,\mathcal{A},P)$, where $\mathcal{S}$ is an
arbitrary set. Moreover, let $\succcurlyeq$ denote a partial order
defined on $\mathcal{S}$ ($x \preccurlyeq y$ if $x$ precedes $y$).
Throughout the paper, we assume that for all $x \in\mathcal{S} $, the
events $\{X\succcurlyeq x \}$ and $\{X\preccurlyeq x \}$ are $\mathcal
{A} $-measurable. We begin by defining the set of points that can be compared
with a fixed element $x\in\mathcal{S}$ given the partial order.
\begin{definition}
\label{Def:CompSets} For any $x\in\mathcal{S}$, the set of points
comparable with $x$ under the partial order $\succcurlyeq$ is defined as
$
\mathcal{C}(x)=\{y\in\mathcal{S}\dvtx  y\succcurlyeq x$ or
$y\preccurlyeq x\}$. Let $p_{x}=P(X \in\mathcal{C}(x))$
denote the probability of comparison of $x$.
\end{definition}
\begin{remark}
It follows that all definitions and results can be derived for general
binary relations $\preccurlyeq$. We focus on partial orders since
these binary relations encompass our applications and to simplify the
exposition. A binary relation $\preccurlyeq$ is a partial order if it
is (i) reflexive ($x\preccurlyeq x$), (ii) transitive ($x\preccurlyeq
y$ and $y\preccurlyeq z$ implies $x\preccurlyeq z$) and (iii)
antisymmetric ($x\preccurlyeq y$ and $y \preccurlyeq x$ implies $x=y$).
Unless otherwise noted, we will assume that the binary relation
$\preccurlyeq$ is a partial order.
\end{remark}

The probability of comparison $p_x$ is simply the probability of
drawing a point comparable with $x$. The usefulness of $\mathcal
{C}(x)$ relies on the fact that conditional on
the event $\{\omega\in\Omega\dvtx X(\omega)\in\mathcal{C}(x)\}$, which
hereafter we denote simply by $\mathcal{C}(x)$, we have
\[
P\bigl(X\succ x | \mathcal{C}(x)\bigr)+P\bigl(X\thicksim x | \mathcal
{C}(x)\bigr)+P\bigl(X\prec
x | \mathcal{C}(x)\bigr)=1.
\]
That is, conditioning on $\mathcal{C}(x)$ avoids points that are
incomparable with $x$ making the partial order $\preccurlyeq$
``complete'' with respect to $x$ [for every $y \in\mathcal{C}(x)$
either $x\preccurlyeq y$ or $y \preccurlyeq x$]. Under this
conditioning, a sensible definition for $x$ being a quantile of $X$
should involve $%
P(X\preccurlyeq x | \mathcal{C}(x))$ and $P(X\succcurlyeq x |
\mathcal{C}%
(x))$, the probabilities of drawing a point preceding $x$ and
succeeding $x$, respectively, under the partial order. Next, we
formally define the concept of partial quantile surfaces and indices.
\begin{definition}
\label{Def:PartialSurface}For each $x\in\mathcal{S}$, we define its
partial quantile index as
%
\begin{equation}\label{DefAlt:tau_x}
\tau_{x}=P\bigl(X\preccurlyeq x|\mathcal{C}(x)\bigr).
\end{equation}
Moreover, for $\tau\in(0,1)$, the $\tau$-partial quantile surface
is defined as
%
\begin{equation}\label{Def:Eq:PQsurfaces}
\mathcal{Q}(\tau)=\bigl\{ x\in\mathcal{S}\dvtx
P\bigl(X\succcurlyeq x|\mathcal{C}(x)\bigr)\geq(1-\tau),
P\bigl(X\preccurlyeq x|\mathcal{C}(x)\bigr)\geq\tau%
\bigr\} .
\end{equation}
\end{definition}

Partial quantile indices provide an ordering notion
for each element of $\mathcal{S}$ relative to its comparable points.
Definition \ref{Def:PartialSurface} also defines a subset of $\mathcal
{S}$ associated with each quantile index $\tau\in(0,1)$. In the case
of a univariate random variable under the natural ordering, $\mathcal
{Q}(\tau)$ is simply the set of $\tau$-quantiles of $X$. Note that we
can have
$x\in\mathcal{Q}(\tau)$ for more than one value of $\tau$ only if
$P(X\sim x|\mathcal{C}(x))>0$. (The same would happen in the
univariate quantile case.)

Next, we select a meaningful representative point, called a $\tau
$-partial quantile point, from each $\tau$-partial quantile surface.
To do that, we use the criterion of maximizing the probability of
drawing a comparable point.
\begin{definition}
\label{Def:PartialPoint} For $\tau\in(0,1)$, a $\tau$-partial
quantile point, or simply a~$\tau$-partial quantile, is defined
as any maximizer of $p_x$ over
$\mathcal{Q}(\tau)$, namely,
%
\begin{equation}\label{Def:Beta-Prob}
x_{\tau}\in\mathop{\arg\max}_{x\in\mathcal{S}} p_{x} \qquad
\mbox{s.t. } x\in\mathcal{Q}(\tau).%
\end{equation}
Also, for each $\tau\in(0,1)$, let $p_{\tau}=p_{x_{\tau}}=P(X \in
\mathcal{C}%
(x_{\tau}))$ be the probability measure of the points comparable with
any $%
\tau$-partial quantile $x_{\tau}$. The set of all $\tau$-partial quantile
points is denoted by $\mathcal{Q}^{\ast}(\tau)=\{x\in\mathcal
{Q}(\tau
)\dvtx p_{x}=p_{\tau}\}$.
\end{definition}

The lack of a complete order in $\mathcal{S}$ is exploited to select a
representative point within the partial quantile surface. This approach
is detached from any geometric aspect of $\mathcal{S}$, yet it
reflects the multivariate nature of the situation as well as the
partial order. Also, note that if we have a complete
order, in which $p_{x}=1$ for all $x\in\mathcal{S}$, then any $x\in
\mathcal{Q}(\tau)$ is a $\tau$-partial quantile. This is exactly what
happens in the univariate case, where multiplicity can also occur.

Partial quantile points $x_\tau$ can also be interpreted as
``approximate quantiles'' in the sense that
\[
P(X\preccurlyeq x_\tau) \geq p_{x_\tau} \cdot\tau \quad\mbox{and}\quad
P(X \succcurlyeq x_\tau) \geq p_{x_\tau} \cdot(1-\tau)
\]
and that the balance is ``correct'' within comparable points
\[
P\bigl(X\preccurlyeq x_\tau| \mathcal{C}(x_\tau) \bigr) \geq \tau \quad\mbox
{and}\quad P\bigl(X \succcurlyeq x_\tau| \mathcal{C}(x_\tau)\bigr) \geq
(1-\tau).
\]
In fact, $x_\tau$ is the ``best approximate quantile'' since it is the
maximizer of the probability of comparisons given the restrictions.

The probability of comparison plays an important role in our
definitions and, consequently, in the interpretation of partial
quantiles. It will allow us to quantify the gap between the
interpretation of partial quantiles and the interpretation of
traditional quantiles where all points are comparable to each other. We
will focus on the following
quantity that characterizes the overall comparability of
partial quantile points uniformly over different quantiles.
\begin{definition} The \textit{partial quantile
comparability} is the minimum proba\-bility of comparison associated with
partial quantile points, namely
%
\begin{equation}\label{Def:p*}
\wp=\min_{\tau\in(0,1)}p_{\tau}.
\end{equation}
\end{definition}

When the comparability $\wp$ is large, the interpretation of partial
quantile points is very similar to traditional quantiles. On the other
hand, if $\wp$ is small, there are partial quantile indices for which
the interpretation of partial quantile points deviates considerably
from that for the traditional quantile since drawing a point that is
incomparable to at least some $\tau$-partial quantile point is likely.
Clearly, if the binary relation $\preccurlyeq$ is a complete order,
like univariate quantiles, we have $\wp=1$. As a side note, (\ref
{Def:p*}) can be written as $\wp=\min_{\tau\in
(0,1)}\max_{x\in\mathcal{Q}(\tau)}p_{x}$, so that $\wp$ is a
saddle point of the probability of comparison.


\subsection{Structural properties}

Next, we move to structural properties implied by the definition. It is
notable that interesting and useful properties can be derived within
the general case.

We say that a mapping $h\dvtx\mathcal{S}\to\mathcal{S}$ is
order-preserving if $x\succcurlyeq y$ implies $h(x)\succcurlyeq h(y)$
and $x\succ y$ implies $h(x) \succ h(y)$.
\begin{proposition}[(Equivariance and invariance)]
\label{Prop:Invariance} 
Let $h\dvtx\mS\to\mS$ be a order-preserving mapping. For an $\mS
$-valued random variable $X$, let $x^X_\tau$, $\mathcal{Q}^X(\tau)$,
$\tau^X_x$, $p_x^X$, $p_\tau^X$ and $\wp^X$ denote the partial
quantile quantities.

Then partial quantile points and surfaces are equivariant under $h$,
namely
\[
x_\tau^{h(X)} = h(x_\tau^X) \quad\mbox{and}\quad  \mathcal
{Q}^{h(X)}(\tau) = h(\mathcal{Q}^X(\tau)),
\]
and partial quantile indices and probability of comparisons are
invariant under $h$, namely
\[
\tau_{h(x)}^{h(X)} = \tau_x^X, \qquad p_{h(x)}^{h(X)} = p_x^X, \qquad  p_\tau
^{h(X)} = p_\tau^X   \quad\mbox{and}\quad  \wp^{h(X)} = \wp^X.
\]
\end{proposition}

Proposition \ref{Prop:Invariance} is simple but very useful. As with
univariate quantiles under the natural ordering, any
order-preserving transformation of the data can be dealt with by
transforming the partial quantiles of $X$. For concreteness, consider
$\mathcal{S}={\mathbb{R}}^{d}$ with $a\succcurlyeq b$ only if $a\geq b$
componentwise. In this case, common examples of invariant
transformations are: translation ($x\mapsto x+z$), positive scaling
($x\mapsto tx$, where $t>0$), and componentwise monotonic transformation
[e.g., $x_{j}\mapsto\ln(x_{j})$, where $x_{j}>0$]. Note that no
assumption on the probability distribution was made in
Proposition~\ref{Prop:Invariance}.

In order to show symmetry, we also require assumptions on the
probability distribution.
\begin{proposition}[(Symmetry)]
\label{Prop:Symmetry} Assume that the probability distribution of $X$
is invariant over a order-preserving mapping $m\dvtx\mS\mapsto\mS$, that
is, $%
P(A)=P(m(A))$ for every measurable $A\subset\mathcal{S}$.
Then if $x_\tau$ is a partial quantile point, $m(x_\tau)$ is also a
partial quantile point; if $z \in\mathcal{Q}(\tau)$, then $m(z) \in
\mathcal{Q}(\tau)$; and $\tau_x = \tau_{m(x)}$.
\end{proposition}

The next lemma shows that transitivity in the partial order is
automatically transferred to the partial quantile indices.
\begin{proposition}[(Transitivity)]
\label{Prop:Transitive} Assume that the binary relation $\preccurlyeq
$ is transitive. Then we have that $%
x\succcurlyeq x^{\prime}$ implies that $\tau_x \geq\tau_{x^{\prime}}$.
\end{proposition}

\section{Estimation of partial quantiles}\label{sec:estimation}

Up to now, we have studied properties of the partial quantiles when the
probability distribution of the random variable of interest is known. Next,
we focus on exploring sample-based partial quantiles viewed as
estimates of
their population counterparts. Following standard notation in the empirical
process literature, we let $\mathbb{P}_{n}(A)=\frac{1}{n}\sum
_{i=1}^{n}1%
\{x_{i}\in A\}$. Also, we let $\mathbb{P}_{n}(A|B)=\mathbb
{P}_{n}(A\cap B)/%
\mathbb{P}_{n}(B)$ if $\mathbb{P}_{n}(B)>0$ and zero otherwise. We
carry out
all of the asymptotic analysis as $n\rightarrow\infty$. We use the
notation $a\lesssim b$ to denote that $a=O(b)$, that is, $a\leq cb$ for all
sufficiently large $n$, for some constant $c>0$ that does not depend on~$n$,
and we use $a\lesssim_P b$ to denote that $a=O_P(b)$. We also use the
notation $a\vee b=\max\{a,b\}$ and $%
a\wedge b=\min\{a,b\}$.

\subsection{Assumptions}

We base our analysis in this and the next section on high-level conditions
E.1--E.6. These high-level conditions are implied by a~variety of more
primitive conditions as discussed below.
\begin{conditionE1*}
The data $X_{i}$, $i=1,\ldots,n$,
are an i.i.d. sequence of $\mathcal{S}$-valued random
variables.
\end{conditionE1*}

The next condition imposes regularity on the family of sets induced by
the partial relation
%
\begin{equation}\label{Def:T}
\mathcal{T}=\bigl\{\mathcal{C}(x),\{y\in\mathcal{S}
\dvtx y\preccurlyeq x\}, \{y\in\mathcal{S}
\dvtx y\succcurlyeq x\}\dvtx  x\in\mathcal{S}\bigr\}.
\end{equation}
\begin{conditionE2*}
For $\bar{p}\in(0,1)$, there is a positive number $v(\bar{p})$ such that
\begin{eqnarray*}
&&\sup_{x\in\mathcal{S},p_{x}\geq\bar{p}} \biggl\vert\frac{\mathbb
{P}%
_{n}(X_i\preccurlyeq x)-P(X\preccurlyeq x)}{p_{x}}\biggr\vert\vee
\biggl\vert\frac{\mathbb{P}%
_{n}(X_i\succcurlyeq x)-P(X\succcurlyeq x)}{p_{x}}\biggr\vert \vee
\biggl\vert\frac{\hat{p}_{x}-p_{x}}{p_{x}}\biggr\vert \\
&&\qquad\lesssim_P%
\sqrt{v(\bar{p})/n}.
\end{eqnarray*}
\end{conditionE2*}

Condition E.2 ensures that the partial order is well-behaved for a
uniform law of large numbers to hold over the sets $\{ X\preccurlyeq x
\}$, $\{ X\succcurlyeq x \}$, and $\mathcal{C}(x)$ for all $x$ in
%
\begin{equation}\label{Def:Cp}
\mathcal{C}_{\bar{p}}=\{x\in\mathcal{S}\dvtx  p_{x}\geq\bar{p}\},
\end{equation}
that is, over points with a minimum requirement on the probability of
comparison. Condition E.2 is implied by several more primitive conditions
on $\mathcal{T}$ [e.g., if $\mathcal{T}$ is a
Vapnik--\v{C}ernonenkis class with VC index $v(\mathcal{T})<\infty$
and mild measurability conditions]. We refer to Alexander \cite
{Alexander}, Pollard
\cite{Pollard1995} and Gin\'{e} and Koltchinskii \cite{GineKoltchinskii}
for several results on deriving bounds for $v(\bar{p})$ under primitive
assumptions. A technical remark is that we require the normalization factor
to be $p_{x}$ for all three terms, which is considerably weaker than
using $%
P(X\preccurlyeq x)$ and $P(X\succcurlyeq x)$.

Alternatively, we could derive all of our results
under the condition
%
\begin{equation}\label{E.2alt}
\sup_{A\in\mathcal{T}}|\mathbb{P}_{n}(A)-P(A)|\lesssim_P \sqrt
{v(\mathcal{T})/n}.
\end{equation}
However, (\ref{E.2alt}) might not lead to results as sharp as E.2 achieves
when $p_{x}$ is small. We refer to Dudley \cite{Dudley2000} and van der
Vaart and Wellner \cite{vdV-W} for a~complete treatment to derive
bounds on $%
v(\mathcal{T})$ leading to (\ref{E.2alt}). Note that if condition (%
\ref{E.2alt}) holds, then condition E.2 is satisfied with $v(\bar
{p})=v(%
\mathcal{T})/\bar{p}^{2}$. It is convenient to keep in mind the case
$0<\bar{p}\leq\wp/2$, for which all partial quantile points
$x_{\tau}$ are contained in $\mathcal{C}_{\bar{p}}$ and therefore
covered by
condition~E.2.\looseness=1

Next, we consider conditions the following identification and
regularity conditions relating probability of comparisons and a metric
$d(\cdot,\cdot)$ for~$\mS$.
\begin{conditionE3*}
There are positive
constants $c$ and
$\alpha\geq1$ such that for every $x\in\mathcal{Q}(\tau)$, we have
\[
p_\tau- p_x \gtrsim c \wedge\inf_{x_{\tau}\in\mathcal
{Q}^*(\tau
)}d(x_\tau,x)^\alpha.
\]
\end{conditionE3*}
\begin{conditionE4*}
For a
compact set of
quantile indices $\mathcal{U}\subset(0,1)$, let $\tau\in\mathcal
{U}$ and
let $\tau^{\prime}$ be in a neighborhood of $\tau$. For every
$x_{\tau
}\in\mathcal{Q}^{\ast}(\tau)$, there exists $x_{\tau^{\prime}}\in
\mathcal{Q}^{\ast}(\tau^{\prime})$ such that:
\[
\textup{(i)}\quad |p_{\tau}-p_{\tau^{\prime}}|\lesssim|\tau-\tau^{\prime
}|^{\gamma
}  \quad\mbox{and}\quad \textup{(ii)}\quad d(x_{\tau},x_{\tau^{\prime}})\lesssim|\tau
-\tau
^{\prime}|.
\]
\end{conditionE4*}
\begin{conditionE5*}
We have that
\[
\sup_{\tau\in\mathcal{U}}  \sup_{x_\tau\in\mathcal{Q}^*(\tau
)} \sup_{y\in\mathcal{S},d(x_\tau,y)\leq r}|\hat{p}_{x_\tau
}-p_{x_\tau}-(\hat{p}%
_{y}-p_{y})|\lesssim_P\phi_{n}(r)/\sqrt{n},
\]
where $\phi_{n}\dvtx{\mathbb{R}}_{+}\rightarrow{\mathbb{R}}_{+}$ is
such that $%
r\mapsto\phi_{n}(r)$ is nondecreasing and concave, and $r\mapsto
\phi
_{n}(r)/r^{\kappa}$ is decreasing for some $\kappa<\alpha$.
\end{conditionE5*}

Condition E.3 is a \textit{restricted} identification condition, that is,
$x_{\tau}$ is a maximizer of the probability of comparison only over
$\mathcal{Q}(\tau)$. Moreover, it allows for partially identified
models in the spirit of Chernozhukov, Hong and Tamer~\cite{CHT2007}.
Condition~E.4 requires that the set-valued mapping $\tau\mapsto
\mathcal{Q}^{\ast}(\tau)$ of partial quantile points is a continuous
correspondence over $\mathcal{U}$. However, it does not restrict
$\mathcal{Q}^{\ast}(\tau)$ to be a singleton, convex, or even
bounded. Condition E.5 is a
standard condition on the criterion function for deriving rates of
convergence of $M$-estimators (see, e.g., van der Vaart and Wellner
\cite{vdV-W}, Theorem 3.2.5). Bounds for $\phi_n$ are available in
the literature for a variety of classes of functions (see van der Vaart
and Wellner \cite{vdV-W}).

Finally, in order to establish functional
central limit theorems, the following mild assumption is is imposed on
the class of sets $\mathcal{T}$ as defined
in~(\ref{Def:T}).\looseness=1
\begin{conditionE6*}
For each $n\geq1$, the process indexed by $\mathcal{T}$
\[
\alpha_{n}(A)= \sqrt{n}\bigl( \mathbb{P}%
_{n}(A)-P(A)\bigr),\qquad A\in\mathcal{T},
\]
converges weakly in $\ell^{\infty}(\mathcal{T})$ to a bounded, mean zero
Gaussian process $Z_P$, indexed by $\mathcal{T}$ with covariance
function $%
P(A\cap B)-P(A)P(B)$ for $A,\break B \in\mathcal{T}$.\looseness=1
\end{conditionE6*}

Condition E.6 is directly satisfied if the class of sets $\mathcal{T}$
satisfies uniform entropy or bracketing conditions and mild measurability
conditions (see \cite{vdV-W}).\looseness=1

Next, we verify these conditions for our main motivational examples.
\begin{example}[(Convex cone partial order)]
\label{Ex:Conic} Let $X$ be an $\RR^d$-valued random variable with a
bounded and differentiable probability density function. Consider the
partial order given by $a\succcurlyeq b$ only if $a-b\in K$, where $K$
is a~proper convex cone (nonempty interior, and does not contain a
line). In this case, we have $P(X\succcurlyeq x)=P(x+K)$ and
$P(X\preccurlyeq x)=P(x-K)$.
\end{example}
\begin{lemma}\label{Lemma:LCPO-assump}
Consider the convex cone partial order setup with a compact set
$\mathcal{U} \subset(0,1)$, and $X$ be an $\RR^d$-valued random
variable bounded and differentiable probability density function. Then,
under i.i.d. sampling\vspace*{1pt} of $X$ (condition \textup{E.1}), we
have that \textup{E.2} with $v(\bar{p})\lesssim d/\bar{p}^2$, \textup{E.5} with
$\phi_n(r) \lesssim (r^{1/2} + n^{-1/4})\sqrt{\log n}$ and
$d(x,y)=\|x-y\|$ and \textup{E.6} hold. Assume further that $X$ has convex
support and the probability density function is strictly positive in
the interior of the support. Then \textup{E.3} holds with $d(x,y)=\|x-y\|$ and
$\alpha= 2$, \textup{E.4(i)} holds with $\gamma =1$, and the mapping
$\tau\mapsto\mathcal{Q}^*(\cdot)$ is upper semi-continuous.
\end{lemma}
\begin{example}[(Acyclic directed graph partial order)]\label
{Ex:AcyclicGRAPHpreference}
Let $X$ be an $\mS$-va\-lued random variable where $|\mS|<\infty$. The
partial order is described by an acyclic directed graph, that is,
$x\preccurlyeq y$ if there is a directed path from $x$ to $y$ in the
graph.
\end{example}
\begin{lemma}\label{Lemma:GRAPH-assump}
Consider a space $\mS$, with $|\mathcal{S}|<\infty$, a partial order
defined over $\mathcal{S}$ by an acyclic directed graph, and let $X$
be an $\mathcal{S}$-valued random variable. Then, under i.i.d.
sampling of $X$ (condition \textup{E.1}), we have that \textup{E.2} with $v(\bar
{p})\lesssim(\log|\mS|)/\bar{p}^2$. Moreover, for $d(x,y)=1\{x\neq
y\}$, we have that \textup{E.3} with any $\alpha\geq0$, \textup{E.5} with $\phi_n(r)
\lesssim1\{ r > 0 \}\sqrt{\log|\mS|}$ and \textup{E.6} hold. Moreover,
\textup{E.4} holds with any $\gamma>0$ if the compact set $\mathcal{U}$ does not
contains a particular finite set of indices.
\end{lemma}

In Section \ref{Sec:Ex}, we discuss other examples where conditions
E.1--E.6 hold.

\subsection{Rate for partial quantile indices}
We start by considering the estimation of the partial quantile indices
$\tau_{x}$ associated with each $x\in\mathcal{S}$, as
defined in (\ref{DefAlt:tau_x}). In order to estimate this parameter,
we define the estimator 
%
\begin{equation}\label{Def:taux}
\hat{\tau}_{x}=
\mathbb{P}_{n}\bigl(X_i\preccurlyeq x|\mathcal{C}(x)\bigr)\qquad \mbox{for each }
x \in\mS. 
\end{equation}
%

A fundamental departure from the univariate case arises from the lack of
comparability between some points. This will oblige us to restrict the set
on which uniform convergence is achieved. The next result establishes that
the convergence of partial quantile indices is uniform over $\mathcal
{C}_{%
\bar{p}}$, which from (\ref{Def:Cp}) is the set of points for which the
probability of drawing a comparable point is at least $\bar{p}$.
\begin{theorem}[(Uniform rate for partial quantile indices)]
\label{Thm:RateFronteir} Assume that conditions \textup{E.1} and \textup{E.2} hold. Then for
any $\bar{p}\in(0,1)$, we have
\[
\sup_{x \in \mathcal{S}, p_{x}\geq\bar{p}}|\hat{\tau
}_{x}-\tau
_{x}|\lesssim_P\sqrt{v(\bar{p})/n}.
\]
\end{theorem}

The convergence is uniform over the set $\mathcal{C}_{\bar{p}}$ under the
condition that $v(\bar{p})=o(n)$, which allows for $v(\bar{p})$ to grow,
that is, for $\bar{p}$ to diminish, as a function of the sample size.
That is
of interest to achieve convergence in the whole space as $n$ grows, and for
increasing-dimension frameworks as proposed by Huber \cite{Huber}.
Theorem \ref{Thm:RateFronteir} allows for the estimation of extreme
partial quantile indices as long as they have a reasonable probability
of comparison.

This result highlights the difficulty of estimating properly the
quantile $%
\tau_{x}$ of points for which comparable points are rare. Intuitively,
if $%
p_x\leq1/n$ there is a nonnegligible probability that our sample might
miss $\mathcal{C}(x)$ completely, since
\[
P\bigl( X_i \notin\mathcal{C}(x), i=1,\ldots,n\bigr)= ( 1 - p_x
)^n \geq
\biggl( 1 - \frac{1}{n} \biggr)^n \geq\frac{1}{3},
\]
which creates ambiguity regarding the choice of the partial quantile
index of $x$.

Within $\mathcal{C}_{\bar{p}}$, the estimation of
the probability of comparison $p_x$ holds uniformly directly from E.2.
However, it is typical for this to hold uniformly over $\mathcal{S}$
in many cases of interest.


For $\tau\in(0,1)$, the natural sample analog of partial quantile
surfaces (\ref{Def:Eq:PQsurfaces}) is
given by
%
\begin{equation} \label{Def:SampleSurface}\qquad
\widehat{\mathcal{Q}}(\tau)=\bigl\{ x\in\mathcal{S}\dvtx%
\mathbb{P}_{n}\bigl(X_i\succcurlyeq x|\mathcal{C}(x)\bigr)\geq(1-\tau),
\mathbb{P}_{n}\bigl(X_i\preccurlyeq x|\mathcal{C}(x)\bigr)\geq\tau%
\bigr\}.
\end{equation}
From Theorem \ref{Thm:RateFronteir} it follows that if $x \in\mathcal
{Q}(\tau)$ and $p_x \geq\bar{p}$, $x \in\mathcal{\widehat Q}(\tau
')$, where $|\tau-\tau'|\lesssim_P \sqrt{v(\bar{p})/n}$.

\subsection{Rate for partial quantile points}

Next, we turn to the estimation of partial quantile points. We are also
interested in deriving rates uniformly over a set of quantile indices. We
will consider uniform estimation over a~compact set $\mathcal
{U}\subset
(0,1) $. Note that, by definition, for any $\tau\in\mathcal{U}$ we
have $%
p_{\tau}\geq\wp$. Intuitively, this ensures that observations are likely
to be on the comparable set of partial quantile points as long as $\wp
$ is
not too small. We consider the following estimator:
%
\begin{eqnarray}\label{Def:Qestimation}
&& \hat{x}_{\tau}\in\mathop{\arg\max}_{x\in\mathcal{S}} \hat
{p}_{x}\nonumber\\
&&\qquad\mbox{s.t. }\mathbb{P}_{n}(X_i\succcurlyeq x)\geq (1-\tau)\cdot\hat p_x
- \epsilon_n, \\
&&\qquad\phantom{\mbox{s.t. }}\mathbb{P}_{n}(X_i\preccurlyeq x)\geq \tau \cdot\hat p_x -
\epsilon_n,\phantom{(1-)}\nonumber
\end{eqnarray}
where $\epsilon_{n}$ is a slack parameter that goes to zero (see
Comment \ref{Remark:En} below). We denote the optimal value in (\ref
{Def:Qestimation}) by
\[
\hat{p}_{\tau}=\hat{p}_{\hat{x}_{\tau}}=\mathbb{P}_{n}(\mathcal
{C}(\hat{x}%
_{\tau})).
\]
\begin{remark}\label{Remark:En}
The introduction of $\epsilon_{n}$ aims to ensure that the feasible
set in (\ref{Def:Qestimation}) is nonempty uniformly over $\tau\in
\mathcal{U}$ with high probability. It suffices to choose $\epsilon
_n$ to bound discontinuities of functions in $\mathcal{T}$ associated
with partial quantile points, namely
\begin{eqnarray*}
\epsilon_n^D &:=& 2 \sup_{\tau\in\mathcal{U}} \sup_{x_\tau\in
\mathcal{Q}^*(\tau)} \limsup_{x\to x_\tau} |\Pn( X_i \preccurlyeq
x ) - \Pn( X_i \preccurlyeq x_\tau)| \\
&&\hspace*{93.5pt}{}\vee|\Pn( X_i \succcurlyeq x )
- \Pn( X_i \succcurlyeq x_\tau)|\vee| \hat p_{x} - \hat p_{x_\tau}|.
\end{eqnarray*}
In the convex cone partial order described in Example \ref{Ex:Conic},
if $X$ is an $\RR^d$-valued random variable with no point mass, with
probability one it follows that $\epsilon_n^0 \leq2d/n$. In the case
of discrete spaces like Example \ref{Ex:AcyclicGRAPHpreference}, we
can take $\epsilon_n = 0$ for $n$ sufficiently large. In more general
cases, it also suffices to choose $\epsilon_n$ to majorize $\epsilon
_n^{D'} := \sup_{x \in \mathcal{S}, p_{x}\geq\wp}|\hat{\tau
}_{x}-\tau
_{x}|$. Under E.1 and E.2, Theorem \ref{Thm:RateFronteir} ensures that
$\epsilon_n^{D'} \lesssim_P \sqrt{v(\wp)/n}$. The latter simplifies
the analysis considerably and does not affect the final rate of
convergence of the estimator, but could introduce a $\sqrt{n}$-bias in
the partial quantile index of the estimator of the partial quantile
point (see Theorem \ref{Thm:Rates1} and Corollary \ref
{Cor:AsymptoticTauPQP} below). We explicitly allow for either choice in
Theorem \ref{Thm:Rates1} since it automatically leads to practical
choices of $\epsilon_n$ in cases of interest, including
Example~\ref{Ex:Conic}.
\end{remark}

In contrast to the estimation of partial quantile indices, where the
convergence is independent of the underlying space $\mS$, the
estimation in (\ref{Def:Qestimation}) brings forth the need to work
with a metric to measure the distance in $\mathcal{S}$
between the estimated and true parameters. It must be noted that the
choice of metric might be application dependent. A possible choice of
metric that relies completely on the partial order to avoid the
geometry of $\mathcal{S}$ is given by\looseness=1
%
\begin{equation} \label{Def:Metricd0}
d_{0}(w,z)=P( \{X\succcurlyeq w\}\bigtriangleup\{X\succcurlyeq
z\}) +P( \{X\preccurlyeq w\}\bigtriangleup\{X\preccurlyeq
z\}),
\end{equation}
where $A\bigtriangleup B=(A\cap B^{c})\cup(B\cap A^{c})$ denotes the
symmetric difference between two sets. A typical choice of metric in
many applications when $\mS=\RR^d$, which is connected to the
geometry, is given by the $\ell_2$-norm $d(w,z) = \|w-z\|$. Moreover,
some identification condition with respect to the particular metric
needs to hold, in our case E.3.

In the analysis of the rate of convergence, one needs to account for
nonstandard issues: the underlying parameter might not be unique, the
empirical criterion function lacks continuity, a restricted
identification condition, and the constraints in (\ref
{Def:Qestimation}) define a random set. For instance, the lack of continuity
of indicator functions will lead to $\phi_n(r)\lesssim
(r^{1/2}+n^{-1/4})\sqrt{\log n}$ in
many cases of interest and would not allow for the usual $\sqrt
{n}$-rate in
general. Examples of nonstandard rates of convergence are given in Kim and
Pollard \cite{KimPollard1990} and van der Vaart and Wellner~\cite{vdV-W}.
Moreover, for each quantile $\tau\in(0,1)$, the identification condition
holds only within $\mathcal{Q}(\tau)$ instead of over the entire
space $\mathcal{S}$. That can lead to a slower rate of convergence
since the partial
quantile surface $\mathcal{Q}(\tau)$ is unknown and needs to be
replaced by a~parameter set that is random.
\begin{theorem}[(Uniform rate for partial quantile points)]
\label{Thm:Rates1} Consider a compact set of quantiles $\mathcal
{U}\subset
(0,1)$ and let $\epsilon_n \geq\epsilon_n^D \wedge\epsilon
_n^{D'}$. Assume that conditions \textup{E.1--E.5} hold for $\mathcal{U}$ and some
metric $d(\cdot,\cdot)$. Then, provided that $v(\wp) = o(n\wp^2)$,
we have
\[
\sup_{\tau\in\mathcal{U}} \inf_{x_{\tau}\in\mathcal{Q}^{\ast
}(\tau
)}d(x_{\tau},\hat{x}_{\tau})\lesssim_P\biggl( \frac{v(\wp
)}{n}+\frac{\epsilon_n^2}{\wp^2}\biggr)^{
{1/2}\wedge{\gamma/(2\alpha)}}\vee r_{n}^{-1},
\]
where
\[
r_{n}^{\alpha}\phi_n(1/r_{n})\leq\sqrt{n}.
\]
\end{theorem}

In the typical case of $\phi_{n}(r)\lesssim(r^{1/2}+n^{-1/4})\sqrt
{\log n}$, if $\gamma/\alpha=1/2$, we have an $n^{1/4}$-rate of
convergence, and if $\gamma/\alpha=1$ we have an $(n/\log
n)^{1/3}$-rate of convergence. Under
mild regularity conditions, the logarithmic term can
be removed in the later case if we are interested on a
single quantile index recovering a $n^{1/3}$-rate of convergence, as in
\cite{KimPollard1990}.
However, it is instructive to revisit Theorem
\ref{Thm:Rates1} in the case of a complete order, for which it turns
out that Theorem \ref{Thm:Rates1} implies a $\sqrt{n}$-rate of convergence.
\begin{corollary}
\label{Cor:CompleteOrder} Under \textup{E.1, E.2} and \textup{E.4(ii)}, if the binary relation
is a~complete ordering, for a compact set $\mathcal{U} \subset(0,1)$
and $\epsilon_n := \epsilon_n^D \wedge\epsilon_n^{D'}$, we
have\looseness=-1
\[
\sup_{\tau\in\mathcal{U}}  \inf_{x_\tau\in\mathcal{Q}^*(\tau)}
d(x_\tau,\hat x_\tau) \lesssim_P \sqrt{v(1)/n }.
\]
\end{corollary}

The presence of a complete order resolves the issues with the
restricted identification condition and discontinuity of the criterion
function since the criterion function becomes constant, namely $\hat
p_x = p_x = 1$ for all $x\in\mathcal{S}$. Also, in this case, the
multiplicity of partial quantiles is the same multiplicity as in the
univariate quantile under the natural ordering, $%
\mathcal{Q}^{\ast}(\tau)=\mathcal{Q}(\tau)$.

Finally, we note that in discrete spaces $\mS$ with $|\mS|<\infty$,
like Example \ref{Ex:AcyclicGRAPHpreference}, for $n$ sufficiently
large, with high probability we perfectly recover the partial quantile
points associated with most indices [a consequence of Lemma \ref
{Lemma:GRAPH-assump} and the metric $d(x,y) = 1\{x \neq y\}$].

\subsection{Asymptotic distributions}\label{Sec:Asymp}

In this section, we discuss the derivation of asymptotic distributions of
quantities defined in this paper.
%
\begin{theorem}[(Asymptotic distribution of partial quantile indices)]
\label{Thm:Inference} Let $\bar{p}>0$ be fixed, and assume that
conditions \textup{E.1},
\textup{E.2} and \textup{E.6} hold. Then, if $v(\bar{p}) = o(n)$, for any $x \in
\mathcal{C}_{\bar p}$
\[
\sqrt{n}(\hat{\tau}_{x}-\tau_{x})\rightsquigarrow N\biggl( 0,\frac
{\tau
_{x}(1-\tau_{x})}{p_{x}}\biggr) .
\]
Moreover, 
the process $\beta_n(x) = \sqrt{n}(\hat{\tau}_{x}-\tau_{x})$
indexed by $\mathcal{C}_{\bar p}$ converges weakly in $\ell^{\infty
}(\mathcal{C}_{\bar p})$ to a bounded, mean zero Gaussian process
$G_P$ indexed by $\mathcal{C}_{\bar p}$ with covariance function given by
\begin{eqnarray*}
\Omega_{z,y} &=& \tau_{z}\tau_{y}  \biggl( \frac{%
P(X \preccurlyeq z\cap X \preccurlyeq y)}{P(X \preccurlyeq z)P(X
\preccurlyeq y)}+\frac{P(\mathcal{C}(z)\cap\mathcal{C}(y))}{p_zp_y}\\
&&\hspace*{24pt}{}-\frac{P(\mathcal{C}(z)\cap X \preccurlyeq y)}{%
p_zP(X \preccurlyeq y)}-\frac{P(X \preccurlyeq z\cap\mathcal{C}(y))}{P(X
\preccurlyeq z)p_y}\biggr)
\end{eqnarray*}
for any $z, y \in\mathcal{C}_{\bar p}$.
%
\end{theorem}

Theorem \ref{Thm:Inference} characterizes the empirical process
associated with the estimation of partial quantile indices. Moreover,
it allows us to make inference on the unknown partial quantile index
associated with the estimated partial quantile point process.
\begin{corollary}\label{Cor:AsymptoticTauPQP}
Assume that the conditions of Theorem \ref{Thm:Rates1} and \textup{E.6} hold.
Then, uniformly over $\tau\in\mathcal{U}$, we have
\[
\sqrt{n}(\tau_{\hat x_\tau}-\tau) = G_P(\hat x_\tau) + o_P(1)
+\sqrt{n}(\hat\tau_{\hat x_\tau} - \tau),
\]
where $\sqrt{n}(\hat\tau_{\hat x_\tau}-\tau)$ is observed.
\end{corollary}

We note that the quantity $\sqrt{n}(\hat\tau_{\hat x_\tau} - \tau
)$ is observed in the estimation, so Corollary~\ref
{Cor:AsymptoticTauPQP} can be used for inference.
In particular, if $P(X\succcurlyeq x)$, $P(X\preccurlyeq x)$ and $p_x$
are continuous in $x$, we have $\sqrt{n}|\hat\tau_{\hat x_\tau} -
\tau| = O_P(\epsilon_n\sqrt{n}/\wp)$. In that case, if $\epsilon_n
= o(\wp/\sqrt{n})$, it establishes that the partial quantile index of
the estimated partial quantile point is $\sqrt{n}$-consistent.

Finally, we turn to the estimation of the partial quantile comparability
that aims to characterize the overall comparability of points. We consider
the estimator given by
%
\begin{equation} \label{Def:HATp*}
\hat{\wp}=\min_{\tau\in\mathcal{U}}\hat{p}_{\tau},
\end{equation}
where $\mathcal{U} \subset(0,1)$ is a compact set sufficiently large. The
next result studies the property of the estimator. It is interesting to note
that one can estimate this quantity at a $\sqrt{n}$-rate under mild
regularity conditions.

We use the following notation. For $\tau\in(0,1)$, let
\[
Z_P(\tau)=\sup_{x_{\tau}\in\mathcal{Q}^{\ast}(\tau)}Z_P(\mathcal
{C}(x_{\tau
})),
\]
where $Z_P$ is a Gaussian process defined as in E.6.
\begin{theorem}[(Asymptotic distribution of partial quantile comparability)]
\label{Thm:p*} Consider a compact set of quantiles $\mathcal
{U}\subset
(0,1)$, let $\epsilon_n \geq\epsilon_n^D \wedge\epsilon_n^{D'}$,
$\epsilon_n^2 = o(n^{-1/2})$, and assume $v(\wp)=o(n\wp^2)$ and that
\textup{E.1--E.6} hold. Assume that the
function $\tau\mapsto p_{\tau}$ is twice continuously differentiable
with a
unique minimum, that is, $\wp=p_{\tau^{\ast}}$ for a unique $\tau
^{\ast}\in\operatorname{int}\, \mathcal{U}$. Then
\[
\sqrt{n}(\hat{\wp}-\wp)=o_P(1)+Z_P(\tau^{\ast}).
\]
\end{theorem}

Theorem \ref{Thm:p*} shows that we have a Gaussian limit for $\sqrt
{n}(\hat\wp- \wp)$
only if the set $\mathcal{Q}^*(\tau^*)$ is single-valued. Otherwise
we should expect non-Gaussian limits. Similar findings of non-Gaussian
limits within generalizations of quantiles have been found in \cite
{EinmahlMason1992}; see Section \ref{sec:add} for a detailed
discussion.

\section{Additional issues}\label{sec:add}

In this section, we discuss several other relevant issues. First, we
discuss robustness to outliers. Next, we study monotonicity properties
of the underlying partial quantiles and their sample counterparts. We
provide conditions under which partial quantile indices and
probabilities of comparison characterize completely the underlying
probability distribution. Then we\vspace*{1pt} establish that under independence and
$(\RR^d,\geq)$, there is a concentration of measure for partial
quantile indices and points. We also develop dispersion measures based
on partial quantiles.
Computational tractability of computing partial
quantiles of a random variable with known probability distribution is
then considered. Finally, we have a detailed comparison with the
generalized quantile processes developed in \cite{EinmahlMason1992}.

\subsection{Robustness to outliers}

Next, we investigate robustness to outlier properties of partial
quantile indices and probabilities of comparison. To do that, we
consider the influence function of these functions. Let $F$ denote the
distribution of $X$ and $F_\varepsilon$ denote a contaminated
distribution by $y \in\mathcal{S}$,
\[
F_\varepsilon= \varepsilon\delta_y + (1-\varepsilon)F.
\]
Viewing the quantities as functions of the probability distribution, we
have $\tau_x(F) = \tau_x$ and $p_x(F)=p_x$. Thus, $\tau
_x(F_\varepsilon)$ and $p_x(F_\varepsilon)$ are the partial quantile
index and probability of comparison associated with $x$ for the
contaminated distribution. Recall that the influence function of a
function $\theta(\cdot)$ at $F$ and~$y$ is defined as
\[
IF_\theta(y,F) = \lim_{\varepsilon\to0} \frac{\theta
(F_\varepsilon)-\theta(F)}{\varepsilon}.
\]
The following result follow (whose proof follows from direct calculation).
\begin{lemma}[(Influence functions)]\label{Lemma:IF}The influence
function for partial quantile indices and probabilities of comparisons
are given by
\[
IF_{\tau_x}(y,F) =
\frac{1\{y \preccurlyeq x\} - \tau_x1\{y\preccurlyeq x \cup
y\succcurlyeq x\}}{p_x}
\]
and
\[
IF_{p_x}(y,F) =
1\{y\preccurlyeq x \cup y\succcurlyeq x\} - p_x.
\]
\end{lemma}

As in the case of univariate quantiles, the influence functions do not
depend on the exact ``place'' of $y$. They only depend on whether $y$
precedes $x$, $y$ is incomparable to $x$, or $x$ precedes $y$. Thus, an
outlier cannot impact probabilities of comparison much nor partial
quantile indices if $p_x $ is far from zero.\vfill\eject

Note that partial quantile points are defined based on $p_x$ and $\tau
_x$. Nonetheless, the influential function associated with partial
quantile points is not defined in the generality of the paper. In
particular, we cannot take differences between elements of $\mathcal
{S}$ unless additional structure is imposed. One could generalize the
influence function to $\lim_{\varepsilon\to0} d(x_\tau(F), x_\tau
(F_\varepsilon))/\varepsilon$ for some metric $d$ defined in
$\mathcal{S}$. However, extending the notion of the influence function
is outside the scope of this work.

\subsection{Characterization properties}

One important question is whether the partial quantile quantities
characterize the underlying probability distribution, as univariate
quantiles do in the univariate case. The answer relies on the richness
of the partial order.

A family of sets $\mathcal{E}$ is said to be a \textit{determining class}
if for any two probabilities measures $\mu, \nu$ such that $\mu(E) =
\nu(E)$ for all $E \in\mathcal{E}$, we have $\mu= \nu$.
Reference \cite{DudleyRAP} contains properties and definitions of determining
classes which is a well studied topic in probability theory \cite
{BagchiSitaram1982,Sitaram1983,Sitaram1984}. The classic example of a
determining class for probabilities measures is $\{ x + \RR^d_- \dvtx x
\in\RR^d\}$.

By definition of probabilities of comparison and partial quantile
indices, we have the identity
\[
p_x \tau_x = P(X \preccurlyeq x ).
\]
Thus, if the family of sets $\{ X \preccurlyeq x \}$, $x\in\mS$, is a
\textit{determining class}, the probabilities of comparison and partial
quantile indices characterize the underlying measure.
\begin{theorem}
If the family of sets $\mathcal{M}(\preccurlyeq) = \{ \{ y \in\mS\dvtx
y \preccurlyeq x\} \dvtx  x \in\mS\}$ is a~determining class, then
partial quantile indices and probabilities of comparison uniquely
determines the probability distribution.
\end{theorem}

Below we show that partial orders described in Examples \ref{Ex:Conic}
and \ref{Ex:AcyclicGRAPHpreference} lead to partial quantiles that
characterize the probability measure.
\begin{lemma}\label{Lemma:CHARACTERIZATIONcone}
If $y \preccurlyeq x$ only if $x-y \in K$ where $K$ is a proper convex
cone, as in Example \ref{Ex:Conic}, we have that $\mathcal
{M}(\preccurlyeq)$ is a determining class.
\end{lemma}
\begin{lemma}\label{Lemma:CHARACTERIZATIONgraph}
If the partial order is given by an acyclic directed graph, as in
Example \ref{Ex:AcyclicGRAPHpreference}, we have that $\mathcal
{M}(\preccurlyeq)$ is a determining class.
\end{lemma}

Recall that a binary relation is said to be antisymmetric if $x
\succcurlyeq y$ and $y \succcurlyeq x$ implies that $x=y$. In general,
it follows that antisymmetry is a necessary condition for the
probability measure to be characterized by the partial quantiles.
Otherwise, any transfer of probability mass within indifferent points
$x\sim y$ would not change probabilities of comparison and partial
quantile indices. Partial orders are antisymmetric by definition.\

\subsection{Monotonicity and partial quantiles}

Recall that for univariate quantiles with the natural ordering, estimated
quantiles are nondecreasing. In this section, we consider monotonicity
properties with respect to the
partial order of the estimated partial quantile surfaces and points.
Similar to the standard univariate quantile case, such properties are
valuable for
interpretation and applicability of the partial quantile concept.

We start with a positive result for the estimation of partial quantile
surfaces. The
following result states that the transitivity in the partial order
translates into monotonicity of the estimated partial quantile indices.
Theorem \ref{Thm:MonotoneQuantileIndex} below is analogous to
Proposition \ref{Prop:Transitive} but deals with estimated partial
quantile indices
instead of the true partial quantile indices.
\begin{theorem}
\label{Thm:MonotoneQuantileIndex} Assume that the binary relation is
transitive. Then, if $x \succcurlyeq y$ we have $\hat\tau_x \geq\hat
\tau_y $.
\end{theorem}

Next, we turn to partial quantile points where monotonicity is more delicate.
In this section our interest lies in cases for which the true partial
quantile points are partial-monotone, that is,
%
\begin{equation} \label{Def:Monotone}
x_{\tau}\succcurlyeq x_{\tau^{\prime}} \qquad\mbox{if } \tau\geq
\tau
^{\prime}.
\end{equation}
In particular, under transitivity, this implies that $x_{\tau}$ is
unique for every $\tau\in(0,1)$. In general, the true partial
quantile points might not be partial-monotone with respect to the
partial order (e.g., Example \ref{Ex:NonUniq}).

However, even if the true partial quantile points are partial-monotone
in the sense of (\ref{Def:Monotone}), the estimated partial quantile
points might violate this partial-monotonicity due to estimation
error.\footnote{This can be observed in Figure \ref
{Fig:UniformQP_Ptau} in Section \ref{Sec:Ex}, where the partial
quantile points for the uniform distribution over the unit square are
estimated. A close inspection of Figure \ref{Fig:UniformQP_Ptau} shows
that $\hat{x}_{0.35}=(0.39,0.44)$ and $\hat{x}_{0.4}=(0.47,0.42)$,
which violates
the partial-monotonicity condition (\ref{Def:Monotone}) although the
true partial quantile points satisfy (\ref{Def:Monotone}), as can be
seen from Example \ref{Ex:Square} in Section \ref{Sec:Ex}.} A similar
lack of
monotonicity is observed in quantile regression when
conditional quantile curves are being estimated, see
Koenker \cite{K2005}. The result of this section is
motivated by techniques recently developed to correct the
lack of monotonicity of estimated conditional quantile
curves in Chernozhukov, Fern\'{a}ndez-Val and Galichon \cite
{CFG2007,CFG2009} and Neocleous and Portnoy~\cite%
{NeocleousPortnoy2008}.

Unlike the quantile index result mentioned above that makes no
assumption in the space, additional structure is needed on the pair
$(\mathcal{S}%
,\succcurlyeq)$. Based on the partial order, define the operations
$\bigvee$ and $\bigwedge$, which denote the least upper bound and the greatest
lower bound, respectively, of any two points in $\mathcal{S}$ (these are
also referred to as the ``join'' and the ``meet''). We assume that
$(\mathcal{S},\succcurlyeq)$ is a \textit{lattice
space}, that is, $\mathcal{S}$ is closed under $\wedge$ and $\vee$.
For example, $( {\mathbb{R}}^d, \geq)$ is a lattice space under the operations
\[
x\bigwedge y = (x_1 \wedge y_1,\ldots,x_d\wedge y_d)
\quad\mbox{and}\quad
x\bigvee y
= (x_1\vee y_1,\ldots,x_d\vee y_d).
\]

Given an initial estimator $\{ \hat x_\tau\dvtx \tau\in(0,1)\}$, we define
its majorant and minorant as
%
\begin{equation} \label{Def:MM}
\hat x^\wedge_\tau= \bigwedge_{\tau^{\prime}\geq\tau,  \tau
^{\prime}\in
(0,1) } \hat x_{\tau^{\prime}} \quad\mbox{and}\quad \hat x^\vee
_\tau=
\bigvee_{\tau^{\prime}\leq\tau, \tau^{\prime}\in(0,1)} \hat
x_{\tau^{\prime}} .
\end{equation}
Note that by construction, $\hat x^\wedge_\tau$ and $\hat x^\vee
_\tau$ are
partial-monotones. They can be thought as upper and lower envelopes
constructed based on the initial estimator. Also note that if $\hat
x_\tau$ is partial-monotone, then we would have $\hat x_\tau= \hat
x^\wedge_\tau= \hat x^\vee_\tau$.

\subsubsection{Rearrangement and the case $({\mathbb{R}}^{d},\geq)$}

Due to its importance in applications, we carry over a monotonization scheme
for the case of $\mathcal{S}={\mathbb{R}}^{d}$ with the partial order being
induced by the convex cone $K={\mathbb{R}}_{+}^{d}$. The particular
structure of the cone is such that $K={\mathbb{R}}_{+}\times\cdots
\times{\mathbb{R}}_{+}$ is the cartesian product of the natural order.

A possible monotonization scheme is given by a componentwise
rearrangement, namely
\[
\hat x^r_{\tau,j} = \inf_y \biggl\{ y \in{\mathbb{R}}\dvtx \int_0^1
1\{ \hat
x_{u,j} \leq y \} \,du \geq\tau\biggr\},\qquad  j = 1,\ldots,d.
\]
Note that $\hat x^r_{\tau}$ is such that $\hat x^\wedge_\tau\leq
\hat
x^r_\tau\leq\hat x^\vee_\tau$. We have the following result.
\begin{theorem}
\label{Thm:Rearrangement} Assume that $x_{\tau}$ is partial-monotone. Then,
for any $\kappa\geq1$,
\[
\int_{0}^{1}\sum_{j=1}^{d}|\hat{x}_{\tau,j}^{r}-x_{\tau
,j}|^{\kappa
}\,d\tau\leq\int_{0}^{1}\sum_{j=1}^{d}|\hat{x}_{\tau,j}-x_{\tau
,j}|^{\kappa
}\,d\tau
\]
with probability one.
\end{theorem}

Chernozhukov, Fern\'{a}ndez-Val and Galichon \cite{CFG2009} had
previously derived this improvement in the estimation by using
rearrangement in the estimation of monotone functions (of which
univariate conditional quantiles are a~particular case).

The usefulness of Theorem \ref{Thm:Rearrangement} is twofold. On the one
hand, it states that we always improve in terms of the $L_{\kappa}$-norm
with respect to the original estimator. On the other hand, it allows us to
check if the partial-monotone assumption is valid.
\begin{corollary}
\label{Corollary:Rearrangement} Assuming that $x_{\tau}$ is
partial-monotone, for any $\kappa\geq1$ we have
\[
\int_{0}^{1}\sum_{j=1}^{d}|\hat{x}_{\tau,j}^{r}-\hat{x}_{\tau
,j}|^{\kappa
}\,d\tau\leq2^{\kappa}\int_{0}^{1}\sum_{j=1}^{d}|\hat{x}_{\tau
,j}-x_{\tau,j}|^{\kappa}\,d\tau.
\]
Consequently, if
\[
\Biggl( \int_{0}^{1}\sum_{j=1}^{d}|\hat{x}_{\tau,j}^{r}-\hat
{x}_{\tau,j}|^{\kappa}\,d\tau\Biggr) ^{1/\kappa}>2\sup_{\tau\in
(0,1)}\Vert\hat{x}%
_{\tau}-x_{\tau}\Vert_{\kappa},
\]
$x_{\tau}$ is not partial-monotone.
\end{corollary}

Note that if conditions E.3 and E.4 are satisfied with $d(x,y)=\Vert
x-y\Vert_{\kappa}=({\sum_{j=1}^{d}}|x_{j}-y_{j}|^{\kappa})^{1/\kappa}$,
the right-hand side of the expression above can be bounded by the rate of
convergence of Theorem \ref{Thm:Rates1}. Therefore, although Corollary
\ref%
{Corollary:Rearrangement} is not a formal statistical test, it can provide
evidence for the lack of partial-monotonicity of partial quantile points
since we can compute the $L_{\kappa}$ distance between $\hat{x}_{\tau}^{r}$
and~$\hat{x}_{\tau}$. The lack of partial-monotonicity of partial quantile
points can arise due to nonuniqueness of partial quantile points. (In
general, it can also arise if the binary relation is not transitive.)

\subsection{Independence, natural ordering and concentration of measure}

Note that in general, even if the components are independent, partial
quantiles can reflect a dependence created by the partial order.
However, if the partial order is given by the componentwise natural
order, some independence carries over.\vspace*{1pt} The next result specializes to
the case where $(\mS,\succcurlyeq)$ is $(\RR^d,\geq)$ and $X$ is an
${\mathbb{R}}^{d}$-valued random variable whose components are
independent with no point mass. In the following, let $q_X(\tau
)=(q_{X_1}(\tau),q_{X_2}(\tau),\ldots, q_{X_d}(\tau))'$ denote the
vector whose components are the $\tau$-quantiles of the components of~$X$.
\begin{theorem}[(Independence, concentration of measure and partial
quantile points)]\label{Thm:Ind}
Consider an ${\mathbb{R}}^{d}$-valued random variable $X$ with no
point mass and the natural partial order $\geq$. If the components of
$X$ are independent, then the partial
quantile points (\ref{Def:Beta-Prob}) satisfy
\begin{eqnarray}
x_{\tau}=q_X\biggl( \frac{\tau^{1/d}}{\tau^{1/d}+(1-\tau
)^{1/d}}\biggr)
\quad\mbox{and}\quad p_{\tau}=\frac{1}{( \tau^{1/d}+(1-\tau
)^{1/d}) ^{d}}\nonumber\\
&&\eqntext{\mbox{for all }  \tau\in(0,1).}
\end{eqnarray}
In particular, we have $x_{0.5}=q_X(0.5)$, and for any $\ell_\kappa
$-norm we have
\[
\|x_\tau- x_{0.5}\|_\kappa\leq\| q_X(\tau) - q_X(0.5)\|_\kappa\qquad
\mbox{for all } \tau\in(0,1).
\]
\end{theorem}

Theorem \ref{Thm:Ind} leads to $x_{0.5} =
(q_{X_1}(0.5),q_{X_2}(0.5),\ldots, q_{X_d}(0.5))'$, the vector with
componentwise medians,
which is intuitively reasonable in terms of the geometry. Moreover, we
observe that for $d\geq1$,
\[
\biggl\vert\frac{\tau^{1/d}}{\tau^{1/d}+(1-\tau)^{1/d}}-\frac
{1}{2}%
\biggr\vert\leq\biggl\vert\tau-\frac{1}{2}\biggr\vert,
\]
so that under independence, partial quantiles are always closer to the
median than univariate
quantiles. Therefore, partial quantiles exhibit a concentration of
measure phenomenon under independence and this partial order. However,
the case of $%
\tau=0.5$ also leads to $\wp=1/2^{d-1}$, which decreases exponentially
fast in the dimension $d$. In contrast, as $\tau$ becomes extreme
(i.e., $%
\tau$ converges to zero or one), $p_{\tau}$ approaches one. The simplicity
of the $d=1$ case follows from the fact that all points are comparable. We
typically lose this advantage as soon as $d>1$, and the degree to which
increases in $d$ make comparisons less likely depends on the partial order,
the probability distribution, and the value of $\tau$. This
illustrates a ``concentration of measure phenomenon'' and a ``curse of
dimensionality for comparisons.''
\begin{remark}[{[Impact of correlations under $(\RR^d,\geq)$]}]
Under $(\RR^d,\geq)$, if the components of $X$ are positively
correlated, the probabilities of comparison tend to be larger than
under independence. However, under negative correlation, the
probabilities of comparison tend to be smaller than under independence.
These reflect cases in which the distributions are more or less aligned
with the partial order.
\end{remark}
\begin{remark}[(Perfect positive correlation)]
In the case $(\RR^d,\geq)$, if a~(strictly) monotone transformation
of the components of $X$ are perfectly positively correlated, we have
$x_\tau= q_X(\tau)$ and $p_\tau= 1$ for every $\tau\in(0,1)$.
This is a trivial case in which multivariate partial quantiles collapse
into the univariate quantiles. Not surprising, the concentration of
measure statement is satisfied with equality.
\end{remark}

Next, we turn to partial quantile indices which also exhibit a
concentration of measure under independence.
\begin{theorem}[(Independence, concentration of measure and partial
quantile indices)]\label{Thm:IndInd}
Consider a ${\mathbb{R}}^{d}$-valued random variable $X$ with no point
mass and the natural partial order $\geq$. If the components of $X$
are independent, then the partial quantile
indices (\ref{DefAlt:tau_x}) satisfy
\[
P( \tau_X \leq\tau) = P\Biggl(\sum_{j=1}^d Z_j \leq\log\biggl( \frac{\tau
}{1-\tau} \biggr) \Biggr),
\]
where $Z_j$ are independent logistic random variables with zero mean,
and variance $\pi^2/3$.

In particular, we have that $P(\tau_X \geq1/2) = 1/2$ and that $\tau
_X$ concentrates on extreme quantiles with respect to the dimension.
Namely, for any positive number~$C$,\vspace{-6pt}
\[
P( |\tau_X - 0.5| \leq0.5 - Cd^{-1/2} ) \leq1/C.\vspace{-4pt}
\]
\end{theorem}

Theorem \ref{Thm:IndInd} yields a concentration of measure for partial
quantile indices under independence. As the dimension grows, a
realization of the random variable is more likely to have an extreme
partial quantile index. Equivalently, a realization of the random
variable is likely to belong to a partial quantile surface $\mathcal
{Q}(\tau)$ for $\tau$ close to zero or one. This has close
connections to the concentration of measure for a uniform distribution
over the $d$-dimensional unit cube, where most of the mass concentrates
on corners. In our case, corners correspond to the extremes zero or
one.\vspace{-3pt}
\begin{remark}[{[$\mathcal{Q}(\tau)$ as a partially-efficient frontier]}]
The notion of a~partial quantile surface can be connected with that of
an efficient frontier. A point $x\in\mathcal{S}$ is said to be in the
efficient frontier of $E$ with respect to a partial order if there is
no point $x^{\prime}\in E$ that dominates $x$ in terms of the partial
order. The definition of partial quantile surfaces allows us to
generalize the concept of efficient frontiers for random variables. In
this case, the support of the possible realizations of $X$ plays the
role of the set $E$.
We can interpret the partial quantile surfaces $\mathcal{Q}(\tau)$ as
partially-efficient frontiers parametrized by $\tau$, the probability
of drawing a preceding point
conditional on it being a~comparable point. Partially-efficient frontiers
for high values of $\tau$ are likely to be of particular interest. It
might be quite difficult to reach a point on the
efficient frontier but much easier to reach a point on a partially-efficient
frontier with $\tau$ close to but not equal to one (as shown by
Theorem \ref{Thm:IndInd} under independence). In such cases, the
partially-efficient frontier notion might be quite appealing. In
particular, if the support of $X$ is $\RR^d$, partially-efficient
frontiers are meaningful while the efficient frontier is empty.
\end{remark}

\subsection{Partial quantile regions}

One common use of univariate quantiles is to provide measures of
dispersion. In this section, we propose an approach to build such
measures of dispersion based on the partial quantiles. Traditionally, a
measure of dispersion would be centered on the median and expanded to
extreme quantiles. In the univariate case, for instance, Serfling~\cite
{Serfling2002} advocates the interval
%
\begin{equation} \label{Def:IntervalUni}
I(\kappa)=\biggl[ q\biggl( \frac{1-\kappa}{2}\biggr) ,q\biggl(
\frac{1+\kappa
}{2}\biggr) \biggr],\qquad \kappa\in[0,1],
\end{equation}
to measure the dispersion of a random variable. With $\kappa=0$,
$I(\kappa) $ is the~me\-dian, and as $\kappa$ increases from zero to
one we obtain an interval with probability at least~$\kappa$.

In the extension to the multivariate case, we shift from
``interval'' to ``region.'' Moreover, in order to use partial
quantiles, we need to specify not only the quantiles but also the minimum
probability of comparison in which we are interested. We define the
partial quantile region
of levels $\theta\in[0,1]$ and $\eta\in[0,1]$ as
%
\begin{eqnarray}\label{Def:IntervalPQ}
\mathcal{R}(\theta,\eta) &=&  \biggl\{ x\in\mathcal{S}
\dvtx
P\bigl(X\preccurlyeq
x|\mathcal{C}(x)\bigr)\geq\frac{1-\theta}{2},\nonumber\\[-8pt]\\[-8pt]
&&\hspace*{5.8pt}P\bigl(X\succcurlyeq x|\mathcal{C}(x)\bigr)\geq\frac{1-\theta}{2},
p_{x}\geq(1-\eta)\cdot p_{\tau_{x}} \biggr\} .\nonumber
\end{eqnarray}

These regions consist of points that are ``typical,'' that is,
nonextreme partial quantiles with respect to the given partial order,
which are more comparable to other points. Thus, partial quantile
regions can help characterize dispersion around typical and comparable points.

The family of sets $\mathcal{R}$ is such that $\mathcal{R}(\theta
,\eta)\subseteq\mathcal{R}(\theta',\eta')$ whenever $\theta\leq
\theta'$ and $\eta\leq\eta'$. By definition, $\mathcal{R}(\theta
,0)$ contains only the partial quantile points
for indices $\tau\in$ $[(1-\theta)/2,(1+\theta)/2]$. On the other
hand, $\mathcal{R}%
(\theta,1)$ contains all the partial quantile surfaces for indices
$\tau\in
[(1-\theta)/2,(1+\theta)/2]$. Note that if we do not constrain
the probability
of comparisons, we would obtain unbounded regions in some situations.
In the
univariate case with the natural order (i.e., a complete order holds),
we recover (\ref{Def:IntervalUni}) since $p_{x}=1$ for every $x\in
{\mathbb{R}}$.

In order to endow the partial quantile region with some probability
coverage, we fix a nondecreasing function $g\dvtx[0,1]\rightarrow[
0,1]$ such that $g(0)=0$ and $g(1)=1$. (A
simple rule would be to set $\eta=\theta$.)
Define
\[
\theta^{\ast}_\kappa=\inf\bigl\{\theta\dvtx P\bigl(X\in\mathcal{R}(\theta
,g(\theta))\bigr)\geq\kappa\bigr\},
\]
and let the dispersion region
\[
\mathcal{R}(\kappa)=\mathcal{R}(\theta^{\ast}_\kappa,g(\theta
^{\ast}_\kappa)).
\]
Therefore, the family $\{\mathcal{R}(\kappa)\dvtx  \kappa\in[
0,1]\}$
satisfies the following properties:
\begin{longlist}
\item \textit{Nested property.} This family of sets is nested,
$\mathcal{R}(0) = \mathcal{Q}^*(0.5)$ and $\mathcal{R}(1) =
\mathcal{S}$;
\item \textit{Coverage property.} $\mathcal{R}(\kappa)$ is the
smallest set in the family with probability at least $\kappa$;
\item \textit{Ordering property.} Any element $x \in\mathcal
{R}(\kappa)$ satisfies $|\tau_x - 0.5| \leq\theta^*_\kappa/2$;
\item \textit{Comparability property.} Any element $x \in\mathcal
{R}(\kappa)$ satisfies $p_x \geq
(1-g(\theta^*_\kappa))p_{\tau_x}$.
\end{longlist}
\begin{remark}
With respect to the estimation of (\ref{Def:IntervalPQ}), results in
Section %
\ref{sec:estimation} can be directly applied to estimate $\mathcal
{R}(\theta,\eta)$
uniformly on $\theta\in[0,1-\varepsilon]$ and $\eta\in
[0,1-\varepsilon]$, where $\varepsilon>0$ is fixed or goes to zero
sufficiently slowly.
\end{remark}

\subsection{Efficient computation}\label{Sec:Comp}

In this section, we turn our attention to the question of whether the
computation of the partial quantiles (\ref{Def:Beta-Prob}) can be
performed efficiently. The notion of efficiency we use is the one in
the computational complexity literature, that is, that it can be
computed in polynomial time with the ``size'' of the problem (usually
the dimension of $\mathcal{S}$; see \cite
{BlumCuckerShubSmale1997,GLS1984,MR1995}).

Such a question is usually tied to regularity conditions on the
relevant objects (in this case, on the probability distribution and on
the partial order) and on the representation of the relevant objects.
For example, the partial order could be given only by an
\textit{oracle}: for every two points in $\mathcal{S}$, the \textit{oracle}
returns the better point or reports that the points are incomparable.
Alternatively, it could
have an explicit format that allows us to exploit additional structure
(a similar idea holds for the representation of the probability
distribution of the random variable).

A simple result that pertains to the case when $\mS$ has a finite
number of elements.
\begin{lemma}
\label{Lemma:Finite} Assume that the cardinality of $\mathcal{S}$ is finite,
that we can compute $P(\{x\})$ for every $x\in\mathcal{S}$, and that
we can
evaluate the partial order for any pair of points in $\mathcal{S} $.
Then we
can compute all the partial quantiles in at most $O(|\mathcal{S}|^2)$
operations.
\end{lemma}

Lemma \ref{Lemma:Finite} explicitly evaluates all points in $\mS$.
Therefore, it might be problematic to rely on it when the cardinality
of $\mathcal{S}$ is large. Moreover, we emphasize that Lemma \ref
{Lemma:Finite} does not
provide any information regarding the case where $\mathcal{S}$ is not
finite. A simple discretization of $\mS\subset\RR^d$ would typically
suffer from the curse of dimensionality (e.g., computational
requirements would be larger than $1/\varepsilon^d$). It is not
surprising that the general case cannot be computed efficiently.
\begin{example}
\label{Ex:BadComp} Let $\mathcal{S}=[0,1]^{d}$ be the unit cube, and assume
that the binary relation is such that $x$ and $y$ are incomparable for
all $%
x,y $ different from an unknown point $x^{\ast}\in\mathcal{S}$ for which
$P(X\succcurlyeq x^{\ast}|\mathcal{C}%
(x))=P(X\preccurlyeq x^{\ast}|\mathcal{C}(x))=1/2$. With no additional
information, it is not possible to approximate $x^{\ast}$ efficiently with
any deterministic method. On the other hand, probabilistic methods have an
exponentially small chance of ever being close to $x^{\ast}$. (This
computational problem is equivalent to maximizing a discontinuous
function over the unit cube.)
\end{example}

Note that Example \ref{Ex:BadComp} is an extreme and, arguably,
uninteresting case. There are many interesting cases for which
additional structure is available and can be explored. Here we will
provide sufficient regularity/representation conditions on the
probability distribution and on the partial order to allow efficient
computation of partial quantiles that require the maximization of the
probability of drawing a comparable
point over a subset of $\mathcal{S}$. These conditions cover many
relevant cases.

Our analysis relies on the following two regularity conditions, one for
the probability distribution and another for the partial order:

\begin{longlist}[C.2.]
\item[C.1.] \textit{Condition on the probability density function.\quad}
Let $\mS= \RR^d$ and let the probability density function $f$ of the
random variable $X$ be
log-concave. That is, for every $x,y\in\mathcal{S}$ and $\lambda\in
[0,1]$, we have
\[
f\bigl(\lambda x+(1-\lambda)y\bigr)\geq f(x)^{\lambda}f(y)^{1-\lambda}.
\]

\item[C.2.] \textit{Condition on the partial order.\quad} For
every $x,
y\in\mathcal{S}$, we have
%
\begin{equation} \label{POconvex}
x\succcurlyeq y  \quad\mbox{only if}\quad x-y\in K,
\end{equation}
where $K$ is a convex cone with nonempty interior.
\end{longlist}

In particular, condition C.1, log-concavity of $f$ over $\mathcal{S}$,
implies that $\mathcal{S}$~is convex. Moreover, a log-concave density
function is unimodal, a useful property to achieve computational
tractability. This is needed because of the representation model we
will be
using. Following the literature on computational complexity for Monte Carlo
Markov Chains (see Vempala \cite{VempalaSurvey} for a survey), we assume
that we can evaluate the density function $f$ at any given point.
Nonetheless, the class of log-concave density functions covers many
cases of
interest, including Gaussian and uniform distributions over convex
sets. As illustrated by Example \ref{Ex:BadComp}, the restriction to
log-concave distributions alone is not sufficient to ensure good
computational properties. Condition C.2 provides sufficient regularity
conditions. The partial orders allowed in (\ref{POconvex}) cover many
cases of practical interest, with $K$ being equal to the nonnegative
orthant or the cone of semi-definite positive matrices.
Now we can state a key equivalence lemma for partial quantile points
under these regularity conditions. It allows to replace the function
$p_x$ by a~variable $p \in[0,1]$ in the formulation of partial
quantile points under C.1 and C.2 which simplifies the optimization
problem considerable.
\begin{lemma}
\label{Lemma:LogConcaveConvexCone} Assume that conditions \textup{C.1} and \textup{C.2} hold.
Then the optimization problem formulation in (\ref{Def:Beta-Prob}) is
equivalent to the following optimization problem:
%
\begin{eqnarray}\label{Prob:Aux}
&& (p_{\tau},x_{\tau})\in\mathop{\arg\max}_{p,x} p \nonumber\\
&&\qquad\mbox{s.t. } P(X\succcurlyeq x)\geq(1-\tau)p,
\nonumber\\[-8pt]\\[-8pt]
&&\qquad\mbox{\phantom{\mbox{s.t. }}} P(X\preccurlyeq x)\geq\tau p, \nonumber\\
&&\qquad\mbox{\phantom{\mbox{s.t. }}} x\in\mathcal{S}, 0\leq
p\leq1.\nonumber
\end{eqnarray}
\end{lemma}

An important consequence of Lemma \ref{Lemma:LogConcaveConvexCone},
due to the log-concavity assumption, is that by a simple change of
variable $p=\exp(v)$, (\ref{Prob:Aux}) can be recast as a convex
programming problem. We will be interested in computing an $\varepsilon
$-approximate solution, that is, a point $x_{\tau}^{\varepsilon}$
such that $|\tau_{x_{\tau}^{\varepsilon}}-\tau|\leq\varepsilon$
and $%
p_{x_{\tau}^{\varepsilon}}\geq p_{\tau}(1-\varepsilon)$.

It is helpful to first consider the case that a membership oracle to
evaluate $P(X\succcurlyeq x)$ and $P(X\preccurlyeq x)$ is available. In
that case, because of Lemma \ref{Lemma:LogConcaveConvexCone}, we can
directly use random walks and simulating annealing proposed in Kalai
and Vempala \cite{KalaiVempala2006} and Lov\'asz and Vempala \cite
{LV2006} to compute an approximate maximizer. Table \ref{Table:ALG}
displays the efficient algorithm.
%
\begin{table}
\caption{The hit-and-run method is a random walk that takes as input a
covariance matrix $T_i$,\break an initial point $(V_i^k,X_i^k)$, a
probability density function $g_i$, and a membership oracle\break for
a~convex set $H(\bar{p})$. The output is a random point whose
distribution is approximately according to $g_i$ restricted to~$H(\bar
{p})$. The simulating annealing procedure changes the power to which
the objective function is raised, $g_i(v,x) = \exp(a_iv_i)$, so that
the probability mass concentrates on the maximum (starting from near
uniform). The final output is a point $X^* \in H(\bar{p})$ such that
with probability $1-\delta$, $p_{X^*}\geq(1-\varepsilon)p_\tau$.
The optimization algorithm\break is based on Kalai and Vempala \cite
{KalaiVempala2006} and Lov\'asz and Vempala \cite{LV2006}}
\label{Table:ALG}
\begin{tabular*}{\tablewidth}{@{\extracolsep{\fill}}rl@{}}
\hline
\multicolumn{2}{@{}c@{}}{\textbf{Optimization algorithm}}\\
\hline
Step 0. & Let $\bar{p} < p_\tau$, $\delta\in(0,1)$, set $m =
\lceil\sqrt{d}\ln\frac{2p_\tau(d + \ln(1/\delta))}{\bar
{p}\varepsilon} \rceil$, $k = \lceil c_o d \log^5 d \rceil$
and \\
& $  a_i = \frac{\bar{p}}{p_\tau} ( 1 + \frac
{1}{\sqrt{d}})^i \mbox{ and } g_i(v,x)=\exp(a_iv)\mbox{,
 for }  i=1,\ldots,m. $\\
Step 1. & Let $(V_0^1,X_0^1),\ldots,(V_0^k,X_0^k)$ be independent
uniform random points from \\[2pt]
& $  H(\bar{p}) := \left\{ (v,x) \in\RR\times\mS:
\begin{array}{l}
\log P( X \succcurlyeq x ) \geq\log(1-\tau) + v,\\
\log P( X \preccurlyeq x) \geq\log\tau+ v, \\
\log\bar{p} \leq v \leq0\\
\end{array}
\right\}$\\[14pt]
&and let $T_0$ be their empirical covariance matrix.\\
Step 2. & For $i=1,\ldots,m$ do the following:\\
&Get independent random samples $(V_i^1,X_i^1),\ldots,(V_i,^k,X_i^k)$
from $g_i$ on $H(\bar{p})$, \\
&using hit-and-run with covariance matrix $T_i$, starting from
$(V_{i-1}^1,X_{i-1}^1),\ldots$,\\
&$(V_{i-1}^k,X_{i-1}^k)$, respectively. Set $T_{i+1}$ to be the empirical covariance matrix of\\
&$X_i^1,\ldots,X_i^k$.\\
Step 3. & Output $\max_{j=1,\ldots,k} p_{X_m^j}$ and the maximizer
point $X^*$.\\
\hline
\end{tabular*}
\end{table}

In the case that only a membership oracle for the probability density
function $f$ is available, we can efficiently approximate $P(X
\preccurlyeq x)$ and $P(X \succcurlyeq x)$ by a factor of
$1+\varepsilon$ again by random walks and simulating annealing as
proposed in Lov\'asz and Vempala \cite{LV2006}. This can be used in
the above algorithm to construct the following result.
\begin{theorem}
\label{Thm:Comp} Assume that conditions \textup{C.1} and \textup{C.2} hold. If we have
a~membership oracle to evaluate the probability density function and to
evaluate the partial order, then for every precision $\varepsilon>0$,
with probability $1-\delta$ we can compute an $\varepsilon$-solution
for a $\tau$-partial quantile polynomially in $d$, $\ln(1/\delta)$
and $1/(p_\tau\varepsilon)$.
\end{theorem}

Theorem \ref{Thm:Comp} establishes that conditions C.1 and C.2 are
sufficient for the existence of an efficient probabilistic method to
approximate partial quantile points.\looseness=-1

\subsection{Comparison with generalized quantile processes}

At this point, it is clarifying to discuss relations with the
interesting work of Einmahl and Mason \cite{EinmahlMason1992}.
These authors proposed a broad class of generalized quantile processes
%
\begin{equation}\label{def:U}
U(\tau)=\min\{\lambda(A)\dvtx P(A)\geq\tau,  A\in\mathbb{A}\}
\end{equation}
for $\tau\in(0,1)$, where $\lambda$ is a continuous function
(usually the volume function) and $\mathbb{A}$ is a chosen family of
sets. Formulation
(\ref{def:U}) does not cover the proposed approach. In particular, the
family of sets in (\ref{def:U}) is nested in $\tau$. One important
difference is the incorporation of a partial order structure which
raises issues of incomparability between points, leading to the use of
conditional probabilities.
Moreover, the focus of \cite{EinmahlMason1992} is on the ${\mathbb
{R}}$-valued process $\{U(\tau)\dvtx\tau
\in(0,1)\}$. In this work, in additional to the process $\{p_{\tau
}\dvtx\tau\in(0,1)\}$, we are interested in other processes such as $%
\{x_{\tau}\dvtx\tau\in(0,1)\}$ and $\{\tau_{x}\dvtx x\in\mathcal{S}\}$, which
are, respectively, $\mathcal{S}$-valued and indexed by $\mathcal{S}$.

The generalized quantile process $U\dvtx(0,1)\rightarrow{\mathbb{R}%
}$ as defined in (\ref{def:U}) is estimated by
%
\begin{equation}\label{Def:Un}
U_{n}(\tau)=\inf\{\lambda(A)\dvtx\mathbb{P}_{n}(A)\geq\tau,A\in
\mathbb{A}\}.
\end{equation}
Einmahl and Mason \cite{EinmahlMason1992} establish an asymptotic
approximation for the
process $\tau\mapsto\sqrt{n}(U_{n}(\tau)-U(\tau))$. However, their
analysis does not apply to partial quantiles. For instance, partial
quantiles are built upon conditional probabilities induced by the
partial order instead of the original probabilities. (This is also very
different from that of Polonik and
Yao \cite{PolonikYao2002}, for which the conditioning is fixed within the
maximization.) In addition, note that (\ref{Def:Un}) automatically implies
that $U_{n}(s)\leq U_{n}(t)$ for $s\leq t$, which is likely to fail in our
case. Their analysis relies on a regularity condition that requires $U$ to
be strictly increasing. Regarding their assumptions, they also impose
E.1, E.2, E.4 and E.6. Note that condition E.5 does not appear in
Einmahl and Mason \cite{EinmahlMason1992} because the objective
function is deterministic.

In our context, we would like to estimate the mapping
$\tau\mapsto p_\tau$ by its sample counterpart
$\tau\mapsto\hat p_\tau$. However, the monotonicity
assumption cannot be invoked in general. In fact, it does
not hold in many cases of interest or under independence
as shown in Theorem \ref{Thm:Ind}. Moreover, our
estimated partial quantiles involve an objective function that is data
dependent, $\hat{p}_{x}=%
\mathbb{P}_{n}(\mathcal{C}(x))$, and not a fixed value as the
objective function in (\ref{Def:Un}). In general, we will not be able
to uniformly estimate the entire function at a $\sqrt{n}$-rate due to
the weaker identification condition, which seems to introduce a bias
even if the $\epsilon_{n}$ term is zero. As in \cite
{EinmahlMason1992} for the process $\sqrt{n}(U_n(\tau)-U(\tau))$,
one should expect possibly non-Gaussian limits for $\sqrt{n}(\hat
{p}_{\tau}-p_{\tau})$ since the partial quantile points might be
nonunique. Since Einmahl and Mason \cite{EinmahlMason1992} are
interested in $U$, they did not study the convergence property of the
points (sets $A\in\mathbb{A}$ in their framework) that achieve the
maximum, as Theorem~\ref{Thm:Rates1} does. Also, there are no analogs
of partial quantile indices in \cite{EinmahlMason1992}.

Finally, note that it is potentially
interesting to apply the machinery of the generalized quantile process of
Einmahl and Mason \cite{EinmahlMason1992} with $\lambda(A)=\mathrm
{volume}%
(A)$ and $\mathbb{A}=\{\mathcal{R}(\kappa)\dvtx  \kappa\in[
0,1]\}$,
since the sets in $\mathbb{A}$ are nested. However, unlike in \cite
{EinmahlMason1992}, the sets in $\mathbb{A}$ are unknown a priori and also
need to be estimated.

\section{Illustrative examples}\label{Sec:Ex}

The following examples illustrate our definitions in different
settings, thereby illustrating some possible characteristics of partial
quantiles. Our intention is to provide some intuition regarding the
behavior of $\tau_x$, $\mathcal{Q}(\tau)$, $x_{\tau}$, $p_{x}$,
$p_{\tau}$ and $\wp$ in a variety of cases and to show that the
interaction between the partial order and the probability distribution
plays a key role.
\begin{example}[(Unit square in ${\mathbb{R}}^{2}$)]
\label{Ex:Square} Let $X\sim\operatorname{Uniform}([0,1]^{2})$, with $%
a\succcurlyeq b$ only if $a\geq b$ componentwise. Note that
\[
P(X\succcurlyeq x)=(1-x_{1})(1-x_{2}), \qquad P(X\preccurlyeq
x)=x_{1}x_{2}
\]
and
\[
p_x=1-x_{1}-x_{2}+2x_{1}x_{2}
\]
characterize the partial quantile indices for every $x\in[0,1]^{2}$.
It follows that to maximize $p_x$ for $x\in\mathcal{Q}(\tau)$, the partial
quantile points are on the diagonal $x_{1}=x_{2}$ and are given by
\[
x_{\tau}=\frac{\tau^{1/2}}{\tau^{1/2}+(1-\tau)
^{1/2}}\pmatrix{1\cr1}\qquad
\mbox{with } p_{\tau}=\frac{1}{1+2\sqrt{\tau(1-\tau)}}.
\]
Figure \ref{Fig:UniformSquare} illustrates the partial quantile
%
\begin{figure}
\begin{tabular}{@{}cc@{}}

\includegraphics{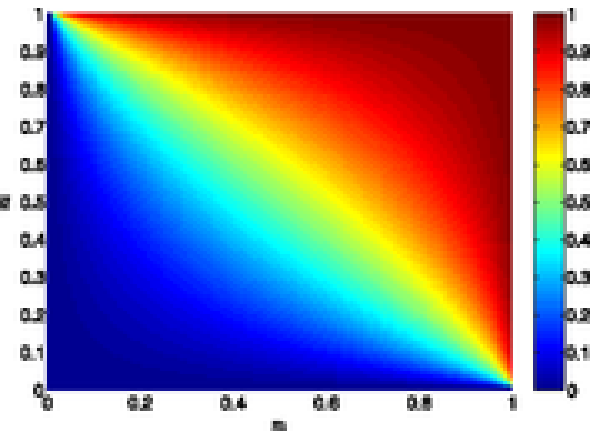}
 &\includegraphics{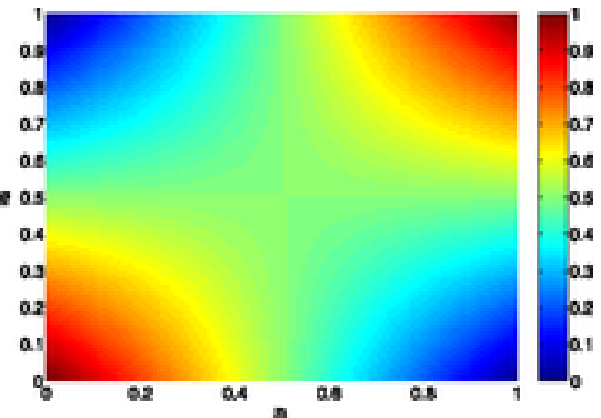}\\
(a) & (b)
\end{tabular}
\caption{\textup{(a)} Partial quantile indices and \textup{(b)}
probabilities of comparison
for $x\in$ $[0,1]^{2}$ in Example~\protect\ref{Ex:Square}.}
\label{Fig:UniformSquare}
\end{figure}
indices $%
\tau_{x}$ and $p_{x}$ for each $x\in$ $[0,1]^{2}$. The shapes of the
partial quantile surfaces can be inferred from the color bands of partial
quantile indices, with each band containing $\mathcal{Q}(\tau)$ for an
interval of values of $\tau$. The symmetry leads to the partial quantiles
being on the diagonal, and we can see from the graph of values of
$p_{x}$ on
the diagonal that $p_{\tau}\rightarrow1$ as $\tau\rightarrow0$ or $1$
and is minimized at the partial median $x_{0.5}=(1/2,1/2)$, with
$\wp=1/2$.
\end{example}

Since partial quantiles generalize univariate quantiles under the
natural ordering, we must
inherit some of its features. For example, multiplicity is possible.
However, we note that in a multidimensional setting with the additional
freedom of a partial order, the set of $\tau$-partial quantiles for a given
$\tau$ does not need to be convex. Multiplicity and nonconvexity of the
set of $\tau$-partial quantiles for a given $\tau$ are illustrated by the
next example, which can be thought of as a mixture of two populations.
In the univariate case, mixtures, just as any other distributions,
always lead to convex collections of quantiles.
\begin{example}[(Nonuniqueness)]
\label{Ex:NonUniq} Consider the random variable
\[
X\sim\operatorname{Uniform}\bigl((-1,1)\times(1,3)\cup(1,3)\times(-1,1)\bigr)
\]
with $a\succcurlyeq b$ only if $a\geq b$ componentwise. In this case, no
points in the square $(-1,1)\times(1,3)$ can be compared with any
point in
the square $(1,3)\times(-1,1)$. This situation leads to nonuniqueness of
the partial quantiles. For $\tau\in(0,1)$, we have
\[
x_{\tau}\in\left\{\pmatrix{-1+2\tau\cr1+2\tau}, %
\pmatrix{1+2\tau\cr-1+2\tau} \right\} \quad\mbox{and}\quad p_{\tau}=
\frac{(1-\tau)^{2}+\tau^{2}}{2}.
\]
Here $\wp= 1/4$ and $p_{\tau}\leq1/2$ for every $\tau\in(0,1)$, because
the two squares are not in alignment with the partial order. See Figure
\ref{Fig:NU} for the representation. Moreover, the set of $\tau$-partial
%
\begin{figure}
\begin{tabular}{@{}cc@{}}

\includegraphics{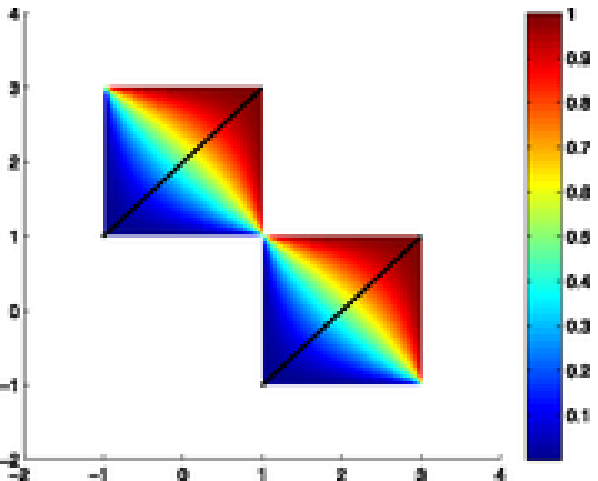}
 &\includegraphics{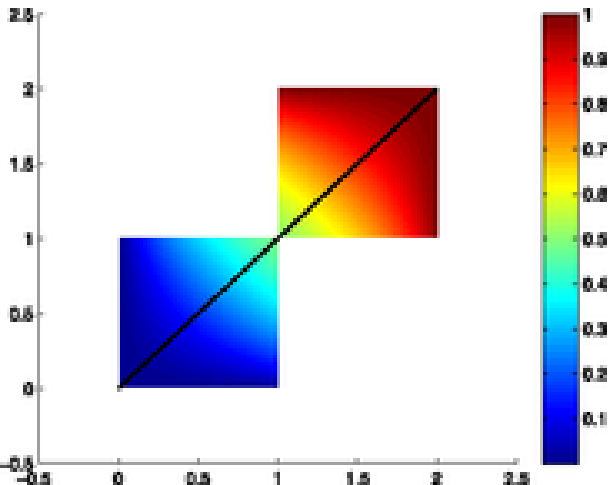}\\
(a) & (b)
\end{tabular}
\caption{\textup{(a)} The potential nonuniqueness of the partial quantiles arising from
the partial order (Example \protect\ref{Ex:NonUniq}). \textup{(b)}
The case of the
partial order being aligned with the probability distribution
(Example~\protect\ref{Ex:Uniq}).}
\label{Fig:NU}
\end{figure}
quantiles for a given $\tau$ is not convex. For example, the set of
$\tau$%
-partial quantiles for $\tau=1/2$ is $\{(0,2)',(2,0)'\}$. The intuitive
geometric notion of a spatial median would report the point
$(1,1)^{\prime }, $ which is not a partial quantile because it is not
comparable with any point in the support of the distribution and thus
having $p_{(1,1)} = 0$.
\end{example}

In the next example, which also involves a mixture of two populations,
the probability distribution is better aligned with the partial order.
\begin{example}[(Aligned distribution and partial order)]
\label{Ex:Uniq} Consider the random variable
\[
X\sim\operatorname{Uniform}([0,1]^{2}\cup[1,2]^{2})
\]
with $a\succcurlyeq b$ only if $a\geq b$ componentwise. The
probabilities of
the events $\{X\succcurlyeq x\}$ and $\{X\preccurlyeq x\}$ are
\begin{eqnarray*}
P(X\succcurlyeq x)&=&\frac{1+(1-x_{1})(1-x_{2})}{2},\\
P(X\preccurlyeq x)&=&\frac{%
x_{1}x_{2}}{2} \qquad\mbox{for } x\in[0,1]^{2}
\end{eqnarray*}
and
\begin{eqnarray*}
P(X\succcurlyeq x)&=&\frac{(2-x_{1})(2-x_{2})}{2},\\
P(X\preccurlyeq x)&=&\frac{%
1+(x_{1}-1)(x_{2}-1)}{2} \qquad\mbox{for } x\in[1,2]^{2}.
\end{eqnarray*}
The partial quantiles can be computed explicitly:
\begin{eqnarray*}
x_{\tau}&=&\frac{\sqrt{1+4( {1}/({2\tau
})-1) }-1}{%
2( {1}/({2\tau})-1) }\pmatrix{1\cr1}
\qquad\mbox{for } \tau<1/2,\\
x_{\tau}&=&\pmatrix{1\cr1}
\qquad\mbox{for }  \tau=1/2
\end{eqnarray*}
and
\[
 x_{\tau}=\biggl( 2-\frac{\sqrt{1+4(
{1}/({%
2(1-\tau)})-1) }-1}{2( {1}/({2(1-\tau)})-1)
}\biggr)
\pmatrix{1\cr1}  \qquad\mbox{for } \tau>1/2.%
\]
Note that in contrast to Example \ref{Ex:NonUniq}, we have $\wp= 3/4$ in
this case since the ordering is somewhat aligned with the distribution [see
Figure \ref{Fig:NU}(b)].
\end{example}

Examples \ref{Ex:NonUniq} and \ref{Ex:Uniq} show the impact the
alignment of
the probability distribution with the partial order can have on the partial
quantiles and on~$p_{x_{\tau}}$. This alignment is good in Example
\ref{Ex:Uniq}, and the partial quantiles are on the main diagonal. Any
point $%
x\in\mathcal{Q}(\tau)$ for some $\tau$ will have a lower $p_{x}$
than $%
x_{\tau}$, the member of $\mathcal{Q}(\tau)$ on the main diagonal. Here
the maximization of the probability of drawing a comparable point leads to
partial quantiles that are consistent with what we might expect. In
Example %
\ref{Ex:NonUniq}, on the other hand, the maximization of the
probability of
drawing a comparable point leads to two partial quantiles for each
value of $%
\tau$. Each of these two partial quantiles seems reasonable in the context
of the square that it is in. Since the two squares are not in alignment with
the partial order, however, the two $\tau$-partial quantiles for a
given $%
\tau$ are disconnected. Results like this are to be expected with such a
lack of alignment. This is analogous to trying to identify a mode with a
bimodal distribution having widely separated modes. 

There are extreme cases in which the probability distribution is not aligned
at all with the partial order, as illustrated by Example \ref{Ex:NonComp}.
%
%
\begin{example}[(Noncomparable)]
\label{Ex:NonComp} Let $X\sim\operatorname{Uniform}(\Delta^{d-1})$, where
$d>1$,
\[
\Delta^{d-1}=\Biggl\{x\in{\mathbb{R}}^{d}\dvtx x\geq
0,\sum_{j=1}^{d}x_{j}=1\Biggr\}
\]
is the $(d-1)$-dimensional simplex, and $a\succcurlyeq b$ only if
$a\geq b$
componentwise. In this case, no two points can be compared. Therefore, we
have $p_{x}=0$ and $P(X\succcurlyeq x|\mathcal{C}(x))=P(X\preccurlyeq
x|%
\mathcal{C}(x))=1$ for all $x\in\Delta^{d-1}$. Definition \ref%
{Def:PartialSurface} yields $\mathcal{Q}^*(\tau) = \mathcal{Q}(\tau
)=\Delta
^{d-1}$ for all $\tau\in(0,1)$ and $\wp= 0$.
\end{example}

Although Example \ref{Ex:NonComp} might suggest a departure from the
traditional quantile definition, it deals with the somewhat extreme
case in
which no points are comparable. This situation is in sharp contrast
with the
complete order that we are accustomed to in the univariate case.
Nonetheless, it provides a~meaningful illustration of a situation in which
no point is better than any other if we rely only on the partial order.
This situation is analogous to trying to compare points on a
Pareto-efficient set, or an efficient frontier, where the points on the
frontier dominate
other points below and to the left of the frontier but the partial
order does not allow us to say that any point on the efficient
frontier is better than any other.

Next, we consider the case of a complete order in detail, as described
earlier. Note that many complete orders are not partial orders since
antisymmetry might fail. Nonetheless, all the quantities proposed here
can be defined analogously.
\begin{example}[(Complete order)]
\label{Ex:h} Suppose that the binary relation $\preccurlyeq$ can be
represented by a
real-valued measurable function, that is, $x\succcurlyeq y$ if and only
if $%
u(x)\geq u(y)$ for some $u\dvtx\mathcal{S}\rightarrow{\mathbb{R}}$. This
is a well-behaved case in which we have a complete order in $\mathcal
{S}$. Therefore, we have
\[
P(X\succcurlyeq x_{\tau})=P\bigl(u(X)\geq u(x_{\tau})\bigr)\geq(1-\tau) %
\quad\mbox{and}\quad P\bigl(u(X)\leq u(x_{\tau})\bigr)\geq\tau.
\]
Consider the (standard) quantile curve $q_{u(X)}\dvtx(0,1)\to\RR$ of the
random variable $u(X)$. Then $p_x = p_\tau= \wp= 1$, $\tau_x =
q_{u(X)}^{-1}(u(x))$, $\mathcal{Q}(\tau) = u^{-1}(q_{u(X)}(\tau))$
and $\mathcal{Q}^*(\tau) = \mathcal{Q}(\tau)$.
\end{example}
%

The situation described in Example \ref{Ex:h} is encountered, for
example, in decision analysis when the consequences in a
decision-making problem are multidimensional in nature and $u$ might be
represented by a payoff or utility function (e.g., Keeney and Raiffa
\cite{KeeneyRaiffa1976}). We emphasize that the reparametrization
allows us to reduce to the standard univariate case, but the partial
quantiles in the original space $\mathcal{S}
$ would be given by the preimage of the function $u$ and could have an
arbitrary geometry even if we have an interval (possibly a point) in
terms of $u$.


In the following example, a random set is the random element of
interest in
the appropriate space under the inclusion ordering (see Molchanov \cite
{Molchanov2005} for precise definitions).
\begin{example}[(Interval covering)]
\label{Ex:Interval} Let $\mathcal{S}$ be the set of all closed
intervals on $%
[0,1]$, and let $X$ be a closed random interval,
\[
X=[\xi_{1},\xi_{2}],\qquad  \xi_{j}\sim\operatorname{Uniform}%
([0,1])  \qquad\mbox{for }  j=1,2.
\]
The partial order is given by $a\succcurlyeq b$ only if $b\subset a$.
Let $%
x=[x_{1},x_{2}]\subset[0,1]$ be an interval. Then we have
\[
P(X\succcurlyeq x)=2x_{1}(1-x_{2}) \quad\mbox{and}\quad P(X\preccurlyeq
x)=|x_{2}-x_{1}|^{2},
\]
which characterize the partial quantile surfaces. Using Anderson's
lemma, and letting $a(\tau) = \sqrt{2(1-\tau)/\tau}$, one
can show that partial quantiles are achieved on symmetric intervals centered
at $1/2$ and given by
\[
x_{\tau}=\biggl[ \frac{1}{2}-\frac{1}{2+2a(\tau)},\frac
{1}{2}+\frac{1}{2+2a(\tau)}\biggr]
\]
and\vspace{-5pt}
\[
p_\tau=\biggl( \frac{1}{1+a(\tau)}%
\biggr)^{2}+2\biggl( \frac{1}{2}-\frac{1}{2+2a(\tau)}\biggr) ^{2}.
\]
\end{example}

\begin{figure}

\includegraphics{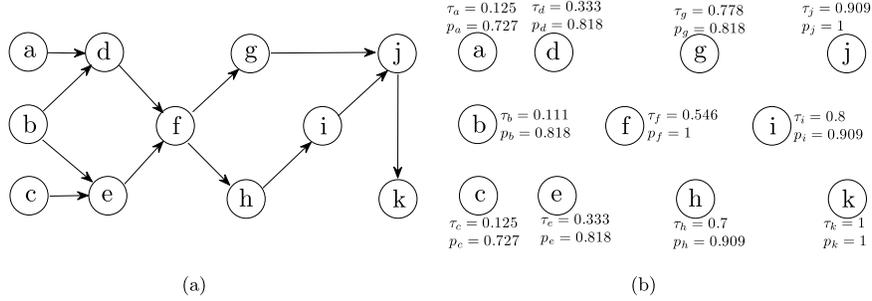}
\vspace{-3pt}
\caption{\textup{(a)} Acyclic directed graph with $x \preccurlyeq y$
if there
is a path from $x$ to $y$. \textup{(b)} Displays partial quantile
indices and
probabilities of comparisons.}\label{Fig:GRAPHexample}
\vspace{-8pt}
\end{figure}

Next, we consider an example of a discrete set $\mathcal{S}$.
\begin{example}[(Partial order based on acyclic directed graphs)]\label
{Ex:AcyclicGRAPH}
Let $X$ be a uniform random variable on $\mathcal{S}=\{
a,b,c,d,e,f,g,h,i,j,k\}$. The partial order relation is given by an
acyclic directed graph, as in Figure \ref{Fig:GRAPHexample}(a), and
$x\preccurlyeq y$ if there is a path from $x$ to $y$ in the graph.
Figure \ref{Fig:GRAPHexample}(b) illustrates how the partial order
relation impacts the partial quantile indices and probabilities of
comparison. Note also that $P( X \preccurlyeq f ) \geq0.5$ and $P(X
\succcurlyeq f)\geq0.5$, making~$f$ the partial median.\vspace{-2pt}
\end{example}
%
%
\begin{figure}[b]
\vspace{-5pt}
\includegraphics{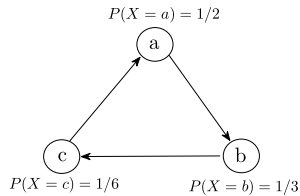}
\vspace{-3pt}
\caption{The cyclic directed graph with $x
\preccurlyeq y$ if there is an arc from $x$ to $y$. The cycle indicates
that the binary relation is not transitive. Moreover, there are no
extreme partial quantiles in this example.}\label{Fig:GRAPHcyclic}
\vspace{-3pt}
\end{figure}

We conclude the examples with a binary relation that is not
transitive.\vspace{-2pt}
\begin{example}[(Nontransitive binary relation)]\label{Ex:CyclicGRAPH}
Let $X$ be a random variable with values in $\mathcal{S}=\{a,b,c\}$,
$P(X=a)=1/2$, $P(X=b)=1/3$ and $P(X=c)=1/6$. The binary\vadjust{\goodbreak} relation is
given by a directed graph, as in Figure \ref{Fig:GRAPHcyclic}, and
$x\preccurlyeq y$ if there is an arc from $x$ to $y$ in the graph. The
cycle in the graph indicates that the binary relation is not
transitive. We note that in this particular example, there are no
extreme partial quantiles. That is, the partial quantile surfaces are
$\mathcal{Q}(\tau) = \varnothing$ for $\tau$ sufficiently close to
$0$ or $1$.
\end{example}

\subsection{Illustration of estimation: The unit square example}

In order to illustrate previous results and statements from Sections
\ref{sec:pq}, \ref{sec:estimation}, \ref{Sec:Asymp} and \ref
{sec:add}, we consider
Example~\ref{Ex:Square} in detail. In this case, $\mathcal{S}=[0,1]^{2}$,
the probability distribution $P$ is the uniform distribution on $[0,1]^{2}$,
and the partial order is given by the $a\succcurlyeq b$ only if $a\geq
b$ (i.e., $a_1 \geq b_1$ and $a_2\geq b_2$), which is
a conic order with $K={\mathbb{R}}_{+}^{2}$. For convenience, we
denote the dimension of $\mathcal{S}$ be $d=2$.

%
\begin{figure}[b]
\begin{tabular}{@{}cc@{}}

\includegraphics{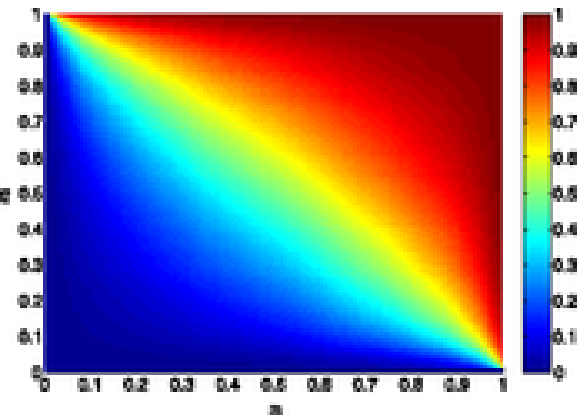}
 &\includegraphics{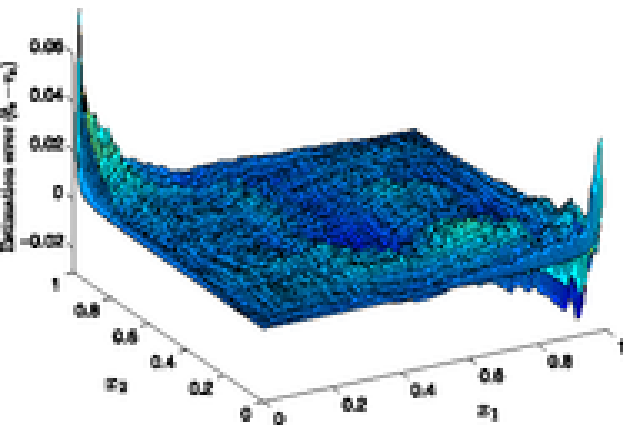}\\
(a) & (b)
\end{tabular}
\caption{\textup{(a)} Estimated partial quantile indices and \textup
{(b)} the difference
between the estimated and true partial quantile indices for uniform samples
on the unit square.}
\label{Fig:UniformSurfaces}
\end{figure}

The class of sets $\mathcal{T}=\{\mathcal{C}(x),\{y\in\mathcal{S}%
\dvtx y\preccurlyeq x\}, \{y\in\mathcal{S}%
\dvtx y\succcurlyeq x\}\dvtx x\in\mathcal{S}\}$ is a VC class of sets whose VC
dimension is of the order $d$, so we have $v(\mathcal{T})\lesssim d$. We
consider the metric to be the usual euclidian norm $d(x,y)=\Vert
x-y\Vert$.
From Theorem~\ref{Thm:Ind}, we have $\wp=1/2^{d-1}$.

Condition E.2 holds with $v(\bar{p})\lesssim d/\bar{p}^2$. Condition
E.3 for $\tau\in
(0,1)$ holds with $\alpha=2$ and $c=1/2^{d}$ (note that for $\tau\in
\{0,1\}$ we would have $\alpha=1$). Condition~E.4 holds with $\gamma=2$
for $\tau=0.5$ and $\gamma=1$ otherwise. Condition E.5 holds with
$\phi
_{n}(r)\lesssim(r^{1/2}+n^{-1/4})\sqrt{\log n}$ by applying maximal
inequalities (the $\log n$ term can be dropped if we are interested in
a single quantile). Finally, condition E.6 holds by an uniform central limit
theorem over $\mathcal{T}$ (see Dudley \cite{Dudley2000}, Theorem
3.7.2, or
van der Vaart and Wellner \cite{vdV-W}, Theorem~2.5.2).

\begin{figure}
\begin{tabular}{@{}cc@{}}

\includegraphics{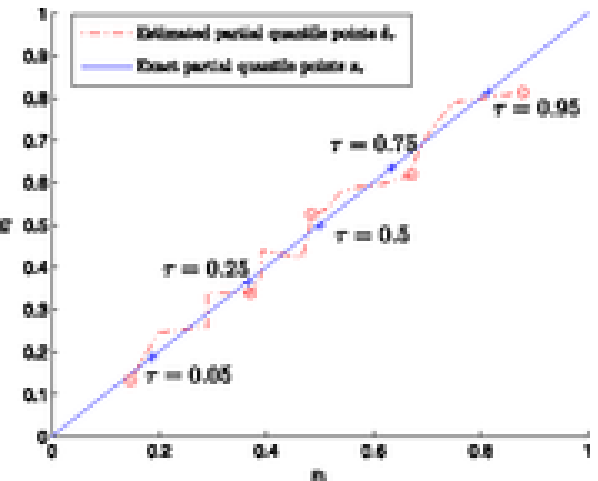}
 &\includegraphics{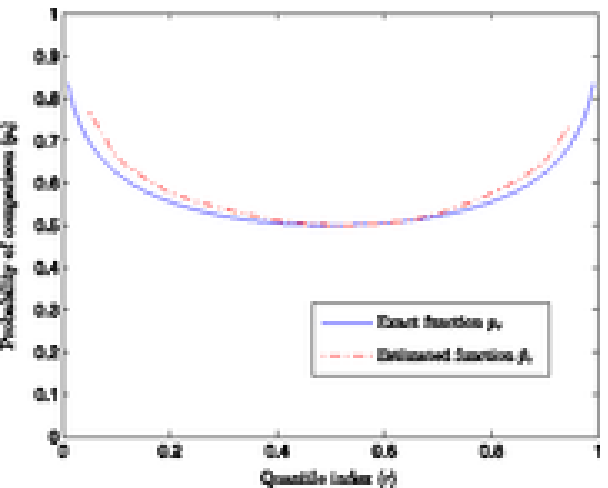}\\
(a) & (b)
\end{tabular}
\caption{\textup{(a)} True and estimated partial quantiles and
\textup{(b)} true and estimated
$p_{\protect\tau}$ as a function of $\protect\tau$ for uniform
samples on
the unit square.}
\label{Fig:UniformQP_Ptau}
\end{figure}

In Figures \ref{Fig:UniformSurfaces} and \ref{Fig:UniformQP_Ptau}, we display
the estimated partial quantile indices and points for the case of $d=2$ with
a sample size of $n=5\mbox{,}000$. Note that the graph of the estimated partial
quantile indices in Figure \ref{Fig:UniformSurfaces} looks very
similar to
the graph of the true partial quantile indices in Figure \ref%
{Fig:UniformSquare}. The difference between the true and estimated
values is
also shown in Figure \ref{Fig:UniformSurfaces}. In light of Theorem
\ref%
{Thm:RateFronteir}, the partial quantile surface is estimated uniformly over
$\mathcal{C}_{\bar p}$ at an $n^{1/2}$-rate of convergence if $\bar
{p}$ is fixed. We
see from the difference between the true and estimated values in Figure
\ref%
{Fig:UniformSurfaces} that the convergence is slower at the top left and
bottom right corners, which correspond to points with small
probabilities of
comparison $p_{x}$.

Although the exact partial quantiles fall on the $x_{1}=x_{2}$
diagonal, we
can see from the few quantiles labeled in Figure \ref{Fig:UniformQP_Ptau}
that they are not evenly spaced along the diagonal. Instead, they are closer
together for $\tau$ near 0.5 and more spread out as $\tau\rightarrow
$ 0
or 1. Moreover, the exact and estimated values of $p_{\tau}$ are smaller
for $\tau$ near 0.5 (the minimum value of the exact $p_{\tau}$ is $%
p_{0.5}=0.5$) and grow larger as $\tau\rightarrow$ 0 or 1. The estimated
quantiles in Figure \ref{Fig:UniformQP_Ptau} are close to but not
equal to
the true quantiles. Also, there is a slight violation of monotonicity
in the
estimated quantiles, a point we will expand upon later.

%

If we are interested in computing partial quantiles only for the case
of $%
\mathcal{U}=\{1/2\}$, we can take $\gamma=2$, which yields a
$n^{1/3}$-rate of
convergence by Theorem~\ref{Thm:Rates1}. Note that for $\mathcal{U}%
=\{0,1\}$ we have $\gamma=1$ and $\alpha=1$, which also leads us to a
$n^{1/3}$-rate
of convergence by Theorem \ref{Thm:Rates1}. On the other hand, if
we are interested in computing for a nondegenerate interval $\mathcal{U}$
of quantiles, we have that $\gamma=1$, which leads to an
$n^{1/4}$-rate of
convergence.

\begin{figure}

\includegraphics{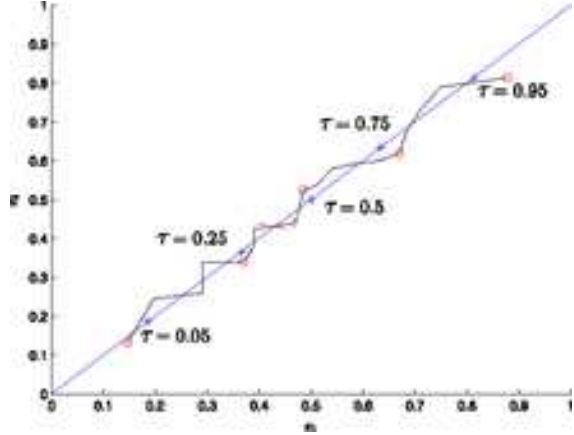}

\caption{The componentwise rearrangement procedure applied to the estimated
partial quantiles from Figure \protect\ref{Fig:UniformQP_Ptau}.}
\label{Fig:Monotone}
\end{figure}

Figure \ref{Fig:Monotone} illustrates the application of the rearrangement
procedure proposed here to the estimated partial quantiles in Figure
\ref%
{Fig:UniformQP_Ptau}, which violated monotonicity for $\tau\in[
0.35,0.40]$. The rearrangement results in estimated partial quantiles that
coincide with the original estimates except for $\tau\in[
0.35,0.40], $ where they are modified to eliminate the violation of
monotonicity.

%

\begin{figure}[b]
\begin{tabular}{@{}cc@{}}

\includegraphics{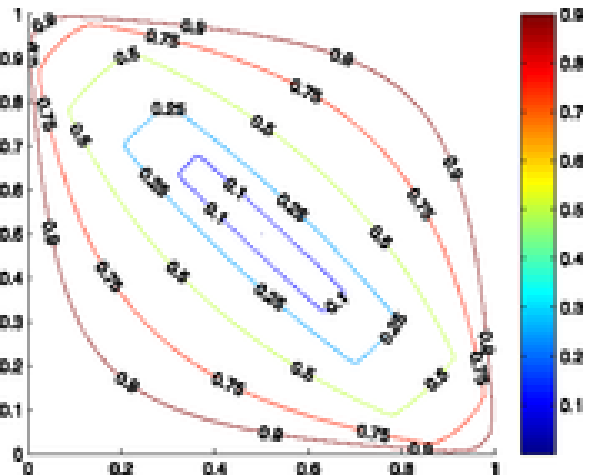}
 &\includegraphics{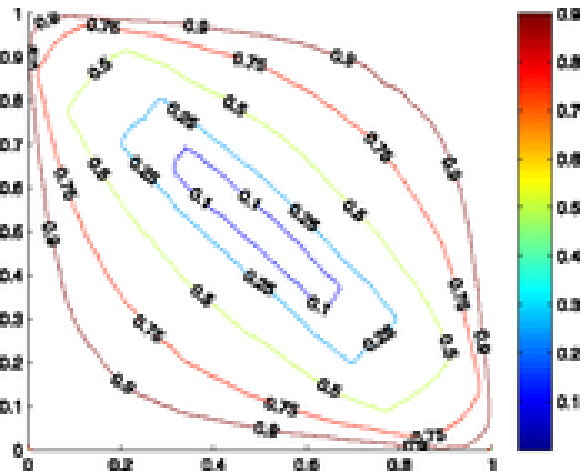}\\
(a) & (b)
\end{tabular}
\caption{\textup{(a)} True and \textup{(b)} estimated dispersion
regions $\mathcal{R}(\theta,\theta)$ for
Example \protect\ref{Ex:Square}, with the boundaries of the regions labeled
by $\theta$.}
\label{Fig:DispersionUniformSquare}
\end{figure}

Exact and estimated dispersion regions with $\eta=g(\theta)=\theta$
for Example~\ref{Ex:Square} are shown in Figure \ref{Fig:DispersionUniformSquare},
corresponding to the exact and estimated partial quantile indices given in
Figures \ref{Fig:UniformSquare} and \ref{Fig:UniformSurfaces}. The
%
dispersion regions seem intuitively reasonable, and the estimated regions
are quite similar to the exact regions. The dispersion regions for high
values of $\theta$ extend out toward $(0,1)'$ and $(1,0)'$, to regions
where the
probabilities of comparison are low. 

\section{Applications}\label{Sec:Applications}




In this section, we use the concept of partial quantiles in two empirical
applications, one involving the intake of dietary components and the other
involving the performance of mutual funds. Our goal is not to do a detailed,
full-scale analysis in each case, but to briefly illustrate the use of
partial quantiles and show some of the capabilities of the concepts and
measures discussed here. In particular, partial quantiles provide useful
graphical and quantitative summaries of the data.

\subsection{Intake nutrients within diets}

Quantitative information regarding the intake distribution of several
dietary components (e.g., calcium, iron, protein, Vitamin A and
Vitamin C)
has been collected by the U.S. Department of Agriculture (USDA) through
periodic surveys. This information is used to formulate food assistance
programs, consumer education efforts, and food regulatory activities. One
important concept in analyzing food consumption data is the usual intake,
defined as the long-run average of daily intakes of dietary components by
individuals. Nusser et al. \cite{NCDF1996}
propose an
approach that assumes the existence of a transformation of the data such
that both the original distribution and measurement errors are normally
distributed. Among other relevant statistics, they estimate the
quantiles of
several dietary components, focusing on each component separately.

For simplicity, we consider only two dietary components, daily intakes of
iron (in milligrams) and protein (in grams), in our analysis. The
partial order is the componentwise natural order. Partial quantiles are
relevant in this situation because not all pairs of diets (as
summarized by their usual intakes) are necessarily comparable in the
sense that we can say that one of
the pair is ``better'' than the other. If one diet has more iron and the other
has more protein, for example, they are not comparable. We recognize that
this partial order rule may not hold for all values of the intakes. At
extremely high levels of a component, it may be undesirable to increase the
intake yet further, but we will assume that the partial order holds within
the range of the data. Another factor that can be relevant is that intakes
of different dietary components are not independent. With this partial
order, for example, a positive correlation between iron intakes and
protein intakes is more in alignment with the partial order and will
lead to higher probabilities of comparison than a negative correlation.
Therefore, understanding this dependence can be important in designing
policies such as those mentioned above. Moreover, the invariance of
partial quantiles under
order-preserving transformations is important since different components
tend to have different scales.

%
\begin{figure}

\includegraphics{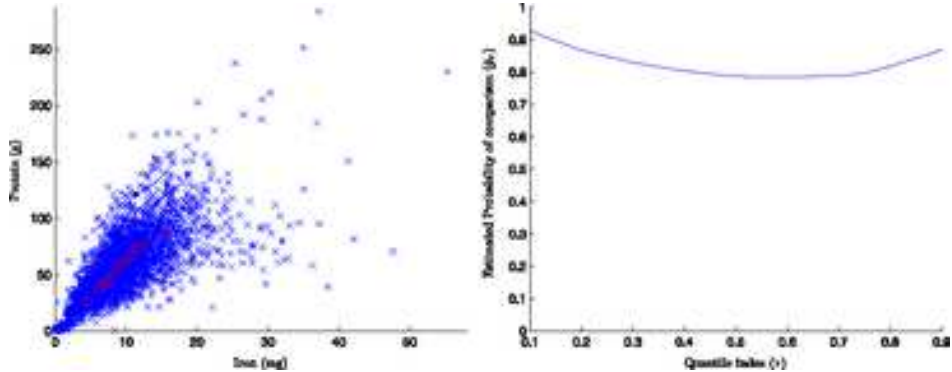}

\caption{\textup{(a)} Data (scatter diagram) and partial quantiles,
and \textup{(b)} estimated
probabilities of comparison $\hat{p}_{\protect\tau}$ for the
(multidimensional) iron and protein levels in food intakes.}
\label{Fig:Nutrients1}
\end{figure}

\begin{figure}
\begin{tabular}{@{}c@{\hspace*{6pt}}c@{}}

\includegraphics{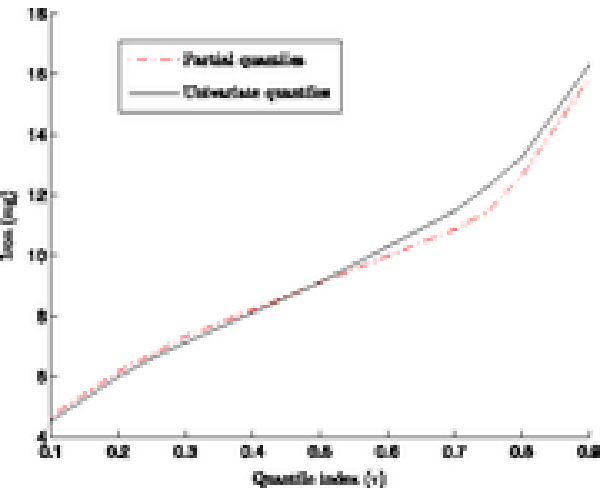}
 &\includegraphics{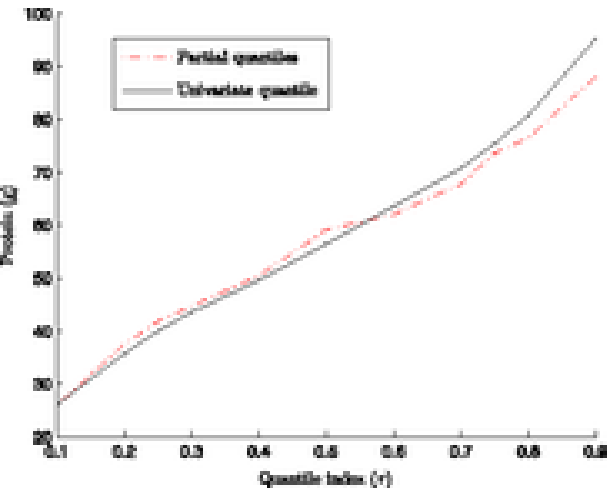}\\
(a) & (b)
\end{tabular}
\caption{Estimated partial quantiles and univariate quantiles for
intakes of
\textup{(a)} iron and \textup{(b)} protein.}
\label{Fig:Q-PQ}
\end{figure}

The data we use are a subset of the data from the 1985 Continuing
Survey of
Food Intakes by Individuals (CSFII) \cite{USDA87}, a data source
used in
\cite{NCDF1996}. A~scatter diagram of the data is given in Figure \ref
{Fig:Nutrients1}, which indicates that the data are quite well-aligned with
the partial order. The estimated partial quantiles shown on
this scatter diagram are monotonically increasing (in terms of the
partial order) in $\tau$. We would expect to see some diets that are not
comparable. Different people may tend to emphasize different types of foods,
with different mixes of nutrients, in their diets. Nonetheless,\vadjust{\goodbreak} the data
indicate that all of the estimated partial quantiles $\hat{x}_{\tau}$ are
comparable with more than $78\%$ of the sampled diets, as can be seen from
Figure~\ref{Fig:Nutrients1}. This suggests that partial quantiles can be
interpreted very similarly to the usual univariate quantiles. For example,
when deriving policies/activities/programs, the decision maker can consider
the $0.5$-partial quantile to be a reasonable representation of the
``median'' individual. Table \ref{Table:Nutrients} and
Figure \ref{Fig:Q-PQ} display comparisons of estimated univariate quantiles and
partial quantiles. In this case, the partial quantiles are slightly more
concentrated around central values than are the univariate quantiles. This
reflects the intuitive notion that it is too extreme to interpret a
componentwise univariate quantile as its multidimensional counterpart. We
note that the univariate quantiles in Table \ref{Table:Nutrients} differ
from those for the same nutrients in \cite{NCDF1996} because we
present the
standard sample quantiles, whereas a measurement error model and assumptions
of normality are used to generate estimated quantiles in~\cite{NCDF1996}.

%
\begin{table}[b]
\caption{Comparison between estimated univariate quantiles and partial
quantiles for iron and protein~intakes}
\label{Table:Nutrients}
\begin{tabular*}{\tablewidth}{@{\extracolsep{\fill}}ld{2.2}d{2.2}d{2.2}d{2.2}@{}}
\hline
\multicolumn{1}{@{}l}{\textbf{Quantile}} & \multicolumn{2}{c}{\textbf{Univariate quantile}} &
\multicolumn{2}{c@{}}{\textbf{Partial quantile}} \\[-4pt]
& \multicolumn{2}{c}{\hrulefill} & \multicolumn{2}{c@{}}{\hrulefill}\\
\multicolumn{1}{@{}l}{\textbf{Index} \textbf{(}$\bolds\tau$\textbf{)}}
& \multicolumn{1}{c}{\textbf{Iron (mg)}} & \multicolumn{1}{c}{\textbf{Protein (g)}}
& \multicolumn{1}{c}{\textbf{Iron (mg)}} & \multicolumn{1}{c@{}}{\textbf{Protein (g)}} \\
\hline
0.1 & 4.51 & 25.95 & 4.69 & 25.97 \\
0.2 & 5.99 & 35.62 & 6.16 & 37.51 \\
0.25 & 6.61 & 39.89 & 6.74 & 41.83 \\
[4pt]
0.3 & 7.12 & 43.53 & 7.33 & 44.72 \\
0.4 & 8.11 & 49.63 & 8.21 & 50.49 \\
0.5 & 9.12 & 56.48 & 9.09 & 59.14 \\
0.6 & 10.29 & 63.61 & 9.97 & 62.03 \\
0.7 & 11.47 & 70.81 & 10.85 & 67.80 \\
[4pt]
0.75 & 12.30 & 75.50 & 11.44 & 73.57 \\
0.8 & 13.25 & 80.82 & 12.61 & 76.45 \\
0.9 & 16.30 & 95.34 & 15.84 & 87.99 \\
\hline
\end{tabular*}
\end{table}
%

%

%
\begin{figure}[b]
\begin{tabular}{@{}c@{\hspace*{6pt}}c@{}}

\includegraphics{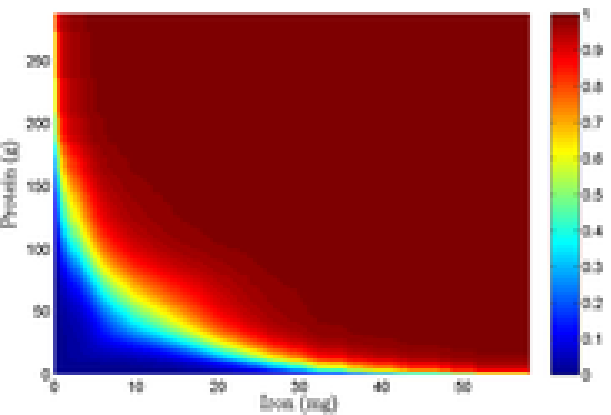}
 &\includegraphics{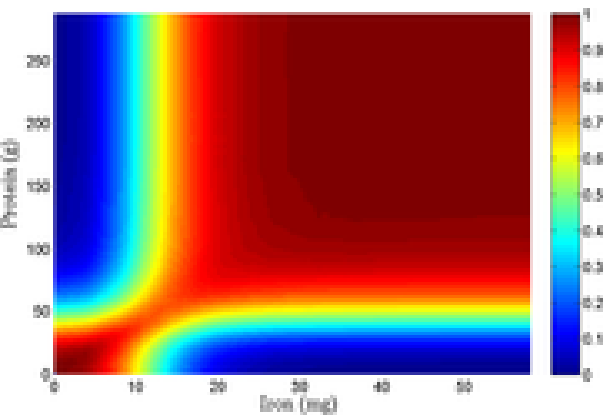}\\
(a) & (b)
\end{tabular}
\caption{\textup{(a)} Estimated partial quantile indices and \textup
{(b)} estimated probabilities of
comparison for levels of iron and protein in food intakes.}
\label{Fig:Nutrients2}
\end{figure}

Figure \ref{Fig:Nutrients2} gives more details, showing the estimated
partial quantile indices $\tau_{x}$ and the probabilities of
comparison $%
p_{x}$ for all $x$. The borders between colors indicating the partial
quantile indices capture the shape of the ``quality'' of the
diets in a comparative sense and show that the partial quantile surfaces
appear convex for these data. For example, a subject with levels of
iron and
protein of $(17.894,87.995)$ will be on the $0.95$ partial quantile surface
among diets that are comparable with her diet, since her diet is on the
upper right-hand border of the light red band in Figure \ref
{Fig:Nutrients2}%
(a). This border can be thought of as a partially efficient frontier of the
intake of iron and protein at a $95\%$ level in this application since any
diets on that border are better than $95\%$ of the comparable diets.
Moreover, this partial quantile surface allows us to consider comparative
statics of the changes needed to stay at the same partial quantile
level but
with higher probabilities of comparison. Note that the graph of the
probabilities of comparison is roughly symmetric, with $p_{x}$
%
%
\begin{figure}ka

\includegraphics{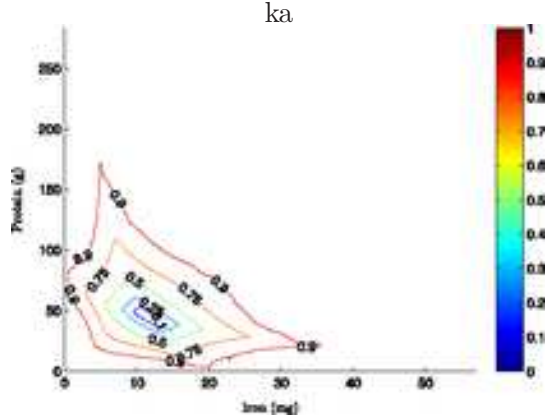}

\caption{The dispersion measure $\mathcal{R}(\theta,\theta)$ based
on estimated
partial quantiles for the (multidimensional) iron and protein levels in food
intakes. The boundaries of the regions are labeled by~$\theta$.}
\label{Fig:Nutrients03}
\end{figure}
decreasing as
we move away from the rough ``axis of symmetry'' along a
particular partial quantile surface. This is consistent with the
location of
the partial quantiles in Figure \ref{Fig:Nutrients1}. Figure \ref%
{Fig:Nutrients03} provides yet additional information by showing the regions
$\mathcal{R}(\theta,\theta)$ from the dispersion measure in (\ref
{Def:IntervalPQ})
for selected values of $\theta$.

\subsection{Evaluating investment funds}

Next, we consider evaluating the performance of investment funds. Several
indices have been considered toward this end in the Finance literature. A
central approach is to regress the return of the fund ($R_F$) above the
return on the risk free asset ($r$) against the return of the market ($R_M$)
above the return on the risk free asset
\[
(R_F-r)=\alpha+\beta(R_{M}-r),
\]
which arises from a standard CAPM model (e.g., \cite{Sharpe1964}). The
exposure with respect to $\beta$ should not be rewarded, and higher values
of the intercept $\alpha$, the risk adjusted return (i.e., the expected
return on the fund when the market yields a return of zero) should be
rewarded.

An emerging literature within finance advocates that in addition to the
risk-adjusted return, \textit{market timing} should also be rewarded
(see \cite%
{HM1981,JK1986,Wermers2000,Andrade2008} and the references therein). The
difference between returns on the market and returns on the fund can be
broken down by whether they are positive or negative to capture market
timing \cite{Andrade2008}:
%
\begin{equation}\label{Eq:DefDELTA}
(R_{F}-r)=\alpha+\beta^{+}\max\{R_{M}-r,0\}+\beta^{-}\min\{
R_{M}-r,0\}.
\end{equation}
Note that $\max\{R_{M}-r,0\}\geq0$ and $\min\{R_{M}-r,0\}\leq0$; a better
performance would have $\beta^{+}$ positive\vadjust{\goodbreak} (the more positive the better)
and $\beta^{-}$ negative (the more negative the better). Therefore, in the
model (\ref{Eq:DefDELTA}), the quantity $\Delta:=\beta^{+}-\beta^{-}$
captures the market timing ability of the fund.
Once again, the partial order that we will use for the pair $(\alpha
,\Delta
)$ is the componentwise natural order.

We use the data used by Andrade in \cite{Andrade2008}. Figure \ref%
{Fig:Finance} shows the data, the estimated partial quantiles, and the
associated probabilities of comparison. Since the partial order is not
complete, we expect to have funds that are noncomparable. In contrast
to the previous application, the data are not well-aligned with the
partial order. It appears that $\alpha$ and $\Delta$ have a strong
negative correlation. As a result, the estimated values for the
probabilities of comparison $p_{\tau
} $ are very small, always below $0.20$ and with $\hat{\wp}=0.00651$.

%
\begin{figure}[t]

\includegraphics{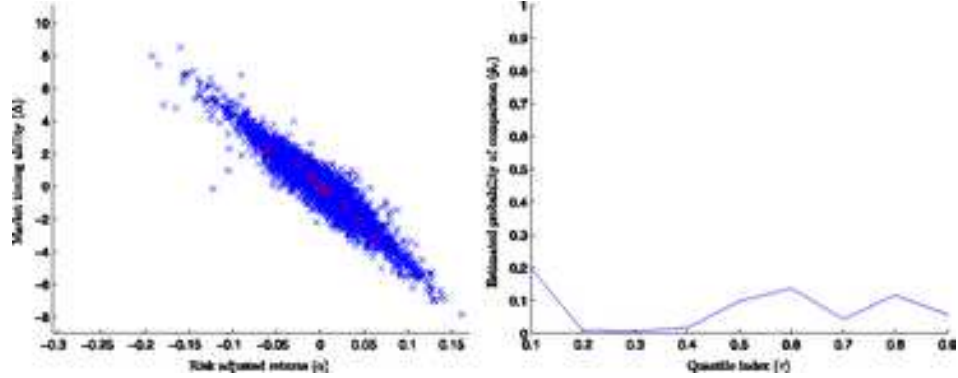}

\caption{Data, estimated partial quantiles, and estimated
probabilities of
comparison for the performance of investment funds.}
\label{Fig:Finance}
\end{figure}

%
\begin{figure}[b]
\begin{tabular}{@{}c@{\hspace*{6pt}}c@{}}

\includegraphics{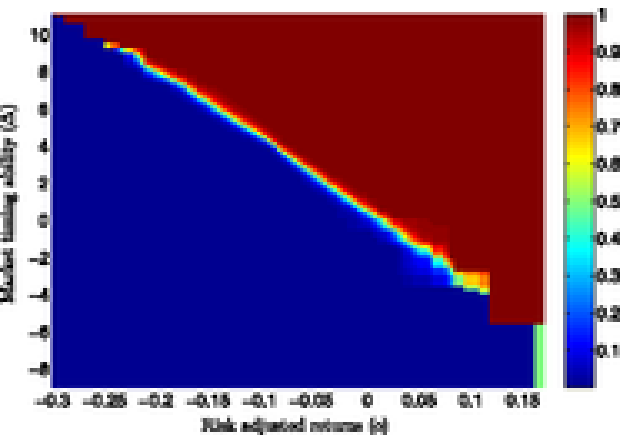}
 &\includegraphics{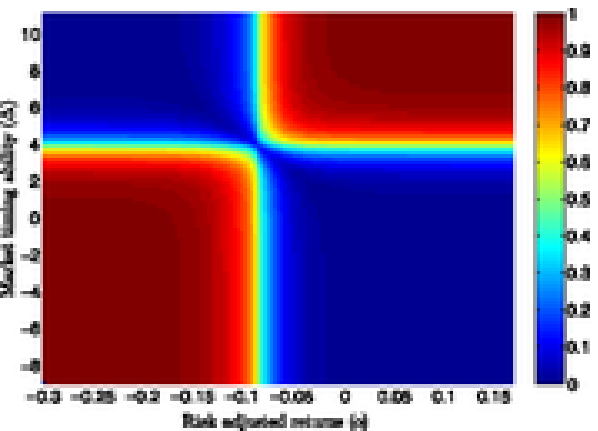}\\
(a) & (b)
\end{tabular}
\caption{\textup{(a)} Estimated partial quantile indices and \textup
{(b)} estimated
probabilities of comparison for the performance of investment funds.}
\label{Fig:FinanceII}
\end{figure}

%
\begin{figure}

\includegraphics{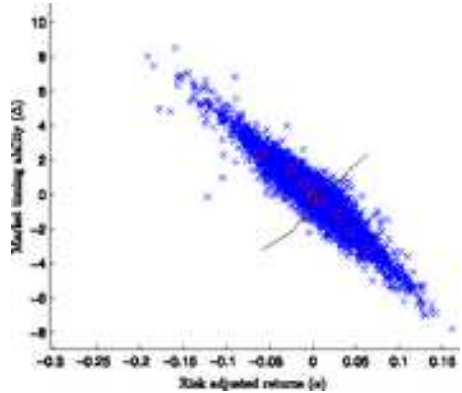}

\caption{The componentwise rearrangement procedure applied to the estimated
partial quantile points for the performance of investment funds. The
difference is $2.141$.}
\label{Fig:FinanceIII}
\end{figure}

Figure \ref{Fig:FinanceII}(a) shows that the partial quantile surfaces for
different values of~$\tau$ are quite close to each other and, except for
extreme values of $\tau$, follow a~pattern that is linear with a negative
slope. This narrow band passes through a region with probabilities of
comparison quite low everywhere, consistent with the above observation
regarding Figure \ref{Fig:Finance}. Therefore, small random variation can
cause potentially large shifts in partial quantile indices. As a
result, the
estimated partial quantiles are not monotonic. When we apply the
rearrangement procedure from Section \ref{sec:add}, we get the results shown
in Figure \ref{Fig:FinanceIII}. The rearranged partial quantiles are
monotonic, but note that many fall outside the support of the data. Moreover,
the $\ell_{2}(\mathcal{U})$ distance between the rearranged and the
original estimator of the partial quantile point process is $2.141$
within the range of $\tau\in(0.1,0.9)$. These observations provide
strong evidence that the true partial quantiles are not
partial-monotone in the sense of (\ref{Def:Monotone}).

How can we interpret the results for this evaluation of investment
funds? We suggest that the results provide some evidence that most (if
not all) of the
funds may actually be optimizing their choices and (up to random
fluctuation) performing on the efficient frontier. Therefore, their
performance is not dominated by many other funds, and when it is, the
differences in performance are slight and seem consistent with random
variation. Similarly, their performance does not dominate many other firms.
This lack of much domination in the data set would explain the low
probabilities of comparability. Since funds have different targets for the
ideal trade-off between risk and return, we should not be surprised to
observe many points on or near different portions of the efficient frontier
in the data, and the data seem to be consistent with this expectation. To
some extent, this is very similar in spirit to Example \ref{Ex:NonComp},
where no point is comparable with any other point.

\subsection{Tobacco and health knowledge scale (THKS)}

We consider the Television School and Family Smoking Prevention
Cessation Project (TVSFP) study (Flay et al. \cite{Flayetal1988} and
Gibbons and Hedeker \cite{GibbonsHedeker1997}), which was designed to
test the effects of a school-based social resistance classroom
curriculum and a media (television) intervention program in terms of
tobacco use prevention and cessation. We refer the reader to \cite
{GibbonsHedeker1997} for the details of the experiment, and we report
the data collected in Table \ref{Table:THKS}.

%
\begin{table}
\tablewidth=250pt
\caption{Tobacco and health knowledge scale postintervention results
subgroups frequencies (and percentages) \protect\cite{GibbonsHedeker1997}}
\label{Table:THKS}
\begin{tabular*}{\tablewidth}{@{\extracolsep{\fill}}lcccc@{}}
\hline
\multicolumn{2}{@{}c}{\textbf{Subgroup}} &
\multicolumn{2}{c}{\textbf{THKS score}} & \\ [-4pt]
\multicolumn{2}{@{}c}{\hrulefill} & \multicolumn{2}{c}{\hrulefill} &\\
\multicolumn{1}{@{}l}{\textbf{CC}} & \multicolumn{1}{c}{\textbf{TV}} &
\multicolumn{1}{c}{\textbf{Pass}} & \multicolumn{1}{c}{\textbf{Fail}}
& \multicolumn{1}{c@{}}{\textbf{Total}} \\
\hline
No & No &  175 & 246 & \phantom{0,}412 \\
& & (41.6) & (58.6) & \\
No & Yes & 201 & 215 & \phantom{0,}416 \\
& & (48.3) & (51.7) & \\
Yes & No & 240 & 140 & \phantom{0,}380 \\
& & (63.2) & (36.8) & \\
Yes & Yes & 231 & 152 & \phantom{0,}383 \\
& & (60.3) & (39.7) &
\\[4pt]
\multicolumn{2}{c}{Total} & 847 & 753 & 1\mbox{,}600\\
& & (52.9) & (47.1) & \\
\hline
\end{tabular*}
\end{table}

\begin{figure}[b]

\includegraphics{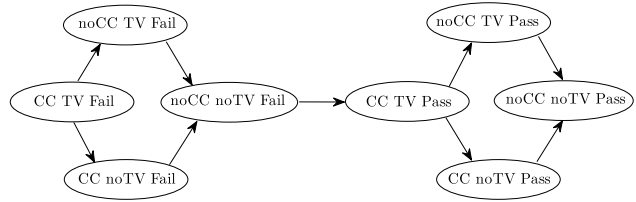}

\caption{The partial order represented by an acyclic directed graph.
We have that $a \preccurlyeq b$ if there is a directed path from $a$ to $b$.}
\label{Fig:GraphSMOKE}
\end{figure}

The partial order of the policy maker is to obtain a ``Pass'' over
``Fail'' regardless of the subgroup. For the same result of the THKS,
given cost and political considerations, it is preferred not to have
used social resistance classroom curriculum (CC) or a media
(television) intervention (TV). However, the subgroup with no CC and TV
is not comparable to CC and no TV. The partial order is summarized by
the acyclic directed graph in Figure \ref{Fig:GraphSMOKE}.

%

\begin{figure}

\includegraphics{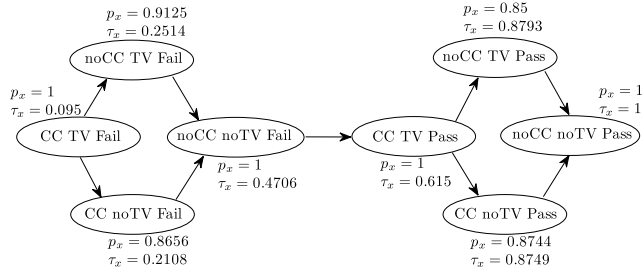}

\caption{The figure displays partial quantile indices and
probabilities of comparisons. According to the partial order of the
policy maker we have $P(X\succcurlyeq\mbox{``CC TV Pass''})\geq1/2$
and $P(X\preccurlyeq\mbox{``CC TV Pass''})\geq1/2$ making ``CC TV
Pass'' the (partial) median.}
\label{Fig:GraphSMOKEwPQ}
\end{figure}

Based on the data of Table \ref{Table:THKS} and the partial order
described in Figure~\ref{Fig:GraphSMOKE}, we compute the partial
quantile indices and probabilities of comparison, see Figure~\ref
{Fig:GraphSMOKEwPQ}.

%

In this application we note the high values of the probability of
comparisons. That makes the interpretation of partial quantiles very
similar to traditional quantiles. In particular, the outcome ``CC TV
Pass'' is such that $P(X\succcurlyeq\mbox{``CC TV Pass''})\geq1/2$ and
$P(X\preccurlyeq\mbox{``CC TV Pass''})\geq1/2$ making ``CC TV Pass''
the (partial) median.

\section{Conclusions}\label{Sec:Conclusion}

We propose a new generalization of quantiles to the multivariate case
based on a given partial order. An important feature of our definition
is that it
is based only on the probability distribution and on the partial order,
which might or not on the geometry of the underlying space. It leads to
a concept that has several desirable properties, including robustness
to outliers and equivatiance/invariance under transformations that
preserve the partial order. Several issues regarding estimation and
computability are investigated and discussed. In particular, rates of
convergence are derived, as are asymptotic distributions of many
quantities, and efficient computation is shown for an important
subclass of distributions and partial orders.

The partial order is the additional structure exploited in this work.
It is clear that partial quantiles depend crucially on the choice of
the partial order. Therefore, their interpretation will also depend
heavily on the partial order. We advocate that the choice of the
partial order is application dependent. Thus, the relevance of these
concepts for a particular application is linked with how meaningful the
partial order is for that application. An alternative approach would be
to choose the partial order to achieve partial quantiles with a desired
property. For instance, one might want partial quantiles with high
probabilities of comparison (which can be achieved with any binary
relation that is a complete order), or partial quantiles that
characterize the probability distribution (which can be achieved if the
partial order induces a determining class), etc. Although these types
of goals can be achieved by the appropriate choice of a partial order,
it is very important for the partial order to make sense in the context
of the specific application because the interpretation of all the
concepts will be tied with that partial order.

Many extensions of the concept of partial quantiles are possible. For
instance, the idea of embedding the partial quantile notion within a
regression framework is of interest, as in
\cite{Chakraborty2003,Chaudhuri1996,ChenGooijer2007,HPS2009,KB78}.
Another possibility is to study the pattern of partial
quantile surfaces conditional on covariates, since partial quantile
surfaces also provide a meaningful generalization of the concept of an
efficient frontier.

Censored models have a wide range of applications and have attracted
considerable interest due to\vadjust{\goodbreak} their connection with quantiles
observed by Powell \cite
{Powell1984,Powell1986a,Powell1986b,Powell1991} and
others \cite%
{NeweyPowell1990,Portnoy2003,NBP2006,BuchinskyHahn1998,Womersley1986}.
However, typical data exhibit censoring in more than one variable. Due to
the equivariance under order-preserving transformations, the proposed
generalization of quantiles is suitable to be applied to censored
multidimensional data.

Moreover, another motivation to consider partial orders, or more
general preferences, is the connection with the literature of decision
theory. For example,
the identification of axioms on the preferences that allow for
statistical inference, computational tractability, etc., is of
interest. Similarly, the identification of classes of decision problems
for which partial quantiles play an important role in optimal
strategies would be very valuable.
Although the pursuit of these extensions is outside the scope of this
paper, we believe that they provide questions of interest for future research.
\vspace{6pt}

\begin{appendix}
\section{\texorpdfstring{Section \lowercase{\protect\ref{sec:pq}} proofs}{Appendix A: Section 2 proofs}}

\begin{pf*}{Proof of Proposition \ref{Prop:Invariance}}
This follows from the equivalence between the events $\{
h(X)\succcurlyeq
h(Y)\}$ and $\{X\succcurlyeq Y\}$, and the events $\{h(X)\succ h(Y)\}$
and $\{X\succ Y\}$.
\end{pf*}
\begin{pf*}{Proof of Proposition \ref{Prop:Symmetry}}
If $m$ is an invariance mapping, it follows that $\mathcal{C}(%
m(x))=m(\mathcal{C}(x))$ and $X\succcurlyeq m(x)=m(X\succcurlyeq x)$.
Therefore,
\begin{eqnarray*}
\Pp\bigl(X\succcurlyeq m(x)|\mathcal{C}(m(x))\bigr)
& = &  \frac{%
\Pp(X\succcurlyeq m(x))}{\Pp(  \mathcal{C}(m(x))  )}=\frac{\Pp(
m(X\succcurlyeq x))}{\Pp(m(\mathcal{C}(x)))}\\
& = & \frac{\Pp(X\succcurlyeq x)}{\Pp(\mathcal{C}%
(x))}=\Pp\bigl(X\succcurlyeq x|\mathcal{C}(x)\bigr).
\end{eqnarray*}
This implies that if $x\in\mathcal{Q}(\tau)$, then $m(x)\in\mathcal
{Q}(\tau)$, and
if $x$ is a $\tau$-partial quantile, so is $m(x)$.
\end{pf*}
\begin{pf*}{Proof of Proposition \ref{Prop:Transitive}}
Since the binary relation is transitive, $\{X\succ x\}\subseteq\{
X\succ
x^{\prime}\}$ and $\{X\preccurlyeq x\}\supseteq\{X\preccurlyeq
x^{\prime
}\}$, so that $P(X\succ x^{\prime})\geq P(X\succ x)\geq0$ and $%
P(X\preccurlyeq x)\geq P(X\preccurlyeq x^{\prime})\geq0$. Therefore,
\begin{eqnarray*}
\tau_{x} & = &  \Pp\bigl(X\preccurlyeq
x|\mathcal{C}(x)\bigr)\\
&=&\frac{\Pp%
(X\preccurlyeq x)}{\Pp(X\preccurlyeq x)+\Pp(X\succ x)}
\geq\frac{\Pp(X\preccurlyeq x^{\prime})}{\Pp%
(X\preccurlyeq x^{\prime})+\Pp(X\succ x)}\\
&\geq&\frac{\Pp%
(X\preccurlyeq x^{\prime})}{\Pp(X\preccurlyeq x^{\prime})+\Pp%
(X\succ x^{\prime})}=\tau_{x^{\prime}}.%
\end{eqnarray*}
\upqed\end{pf*}

\section{\texorpdfstring{Section \lowercase{\protect\ref{sec:estimation}} proofs}{Appendix B: Section 3 proofs}}

\begin{pf*}{Proof of Lemma \ref{Lemma:LCPO-assump}}
We can assume that $X$ has a compact support to ensure that integrals
are well\vadjust{\goodbreak} defined (and standard approximation arguments yields the full
result, or we are establishing probabilistic bounds and the compact set
is chosen to control the probability).

Since $K$ is a convex set, the associated class of functions $\mathcal
{T}$ is measurable and $\partial K$ has zero Lebesgue measure by Lemma
2.4.3 in Dudley \cite{Dudley2000}. Moreover, $\mathcal{T}$ is a VC
class of sets with VC index at most $3d+4$. Therefore, condition E.2
holds with $v(\bar p) = (3d+4)/\bar p^2$.

Let $\sigma_0$ denote the surface measure on $\partial K$. To
establish E.5, let
\[
\mu:= \sup_{x \in\RR^d} \int_{\partial(-K\cup
K)} f(x+y) \,d\sigma_0(y) < \infty,
\]
since the support of $X$ is
compact. Next, note that $d(x,y) = \|x-y\| \geq E[|1\{ X \in\mathcal
{C}(x)\} - 1\{ X \in\mathcal{C}(y)\}|^2]/\mu$. Then E.5 holds with
$\phi_n(r)\lesssim( \sqrt{\mu r} + n^{-1/4} ) \sqrt{\log n}$ by
Theorem 2.14.17 of van der Vaart and Wellner \cite{vdV-W}. If
$\mathcal{U}$ is a singleton, we can improve the bound to $\phi
_n(r)\lesssim\sqrt{\mu r} + n^{-1/4}$ using arguments in Kim and
Pollard~\cite{KimPollard1990}.

Since $\mathcal{T}$ is a VC class and $K$ is a convex set which
ensures enough measurability, E.6 holds by Theorem of 2.6.8 in van der
Vaart and Wellner \cite{vdV-W}.

To establish E.3, building upon Section 5 in Kim and Pollard \cite
{KimPollard1990}, note that
\begin{eqnarray*}
\nabla\tau_x &=& \frac{1}{p_x} \int_{\partial(-K)}
f(x+y)n_{(-K)}(y) \,d\sigma_0(y) \\
&&{}- \frac{\tau_x}{p_x} \int_{\partial
(-K\cup K)} f(x+y)n_{(-K\cup K)}(y)\,d\sigma_0(y)
\end{eqnarray*}
and
\[
\nabla^2p_x = \int_{\partial(-K\cup K)} \nabla
f(x+y) n_{(-K\cup K)}(y)' \,d\sigma_0(y),
\]
where $n_A(y)$ is the outward pointing unit vector normal to $\partial
A$ at $y$. Letting $ B_1 = \partial(-K\cup K)\cap\partial(-K)$, $B_2
= \partial(-K\cup K)\setminus\partial(-K) \subset\partial K$, we
have
\begin{eqnarray*}
p_x\nabla\tau_x & = & \int_{\partial(-K)\setminus
B_1} f(x+y)n_{(-K)}(y) \,d\sigma_0(y) \\
&&{}+ (1-\tau_x)\int_{B_1}
f(x+y)n_{(-K)}(y) \,d\sigma_0(y) \\
&&{} - \tau_x \int_{B_2} f(x+y)n_{(-K\cup K)}(y)\,d\sigma_0(y)\\
&=& \int_{\partial(-K)}\bigl( 1\{y\in B_1^c\} + (1-\tau_x)1\{y \in B_1\} +
\tau_x1\{y\in-B_2\}\bigr)\\
&&\hspace*{28.1pt}{}\times f(x+y)n_{(-K)}(y)\,d\sigma_0(y).
\end{eqnarray*}
Since $-K$ is a convex cone with nonempty interior, the normal vectors
cannot be (positively) linearly dependent. Therefore, we have $\nabla
\tau_x \neq0$ for any $x$ in the interior of the support of the
random variable $X$. Therefore, $\mathcal{Q}(\tau) = \tau
_x^{-1}(\tau)$ is a continuously differentiable hypersurface for every
$\tau\in(0,1)$ by the Global Implicit Function theorem. The
smoothness of $p_x$ and $\mathcal{Q}(\tau)$ yields condition E.3 with
$\alpha= 2$ for all $\tau\in\mathcal{U}$.

Also, $p_x = \int_{-K\cup K} f(x+y) \,dy$ and $\tau_x = (1/p_x) \int
_{-K} f(x+y) \,dy$ are twice differentiable functions. Therefore, $p_\tau
$ is Lipschitz for $\tau\in\mathcal{U}$ since $\mathcal{U}\subset
(0,1)$ is compact and $\wp> 0$ under our conditions. Thus, condition
E.4(i) is satisfied with $\gamma= 1$. Moreover, continuity of $p_x$
and $\tau_x$ also implies that the mapping $\mathcal{Q}^*(\tau)$ is
upper-semi continuous.
\end{pf*}
\begin{pf*}{Proof of Lemma \ref{Lemma:GRAPH-assump}}
The bound on E.2 follows from the union bound. Condition E.3 follows from the
finite cardinality of $\mS$ since for $x\in\mathcal{Q}(\tau
)\setminus\mathcal{Q}^*(\tau)$ we have $p_x < p_\tau$ and for $x
\in\mathcal{Q}^{\ast}(\tau)$ we have $p_x = p_\tau$. Take $c =
\min_{\tau\in\mathcal{U}} p_\tau- \max_{x \in\mathcal{Q}(\tau
)\setminus\mathcal{Q}^*(\tau)} p_x > 0$ since $\mathcal{U}$ is
compact. Condition E.5 follows similarly to~E.2, noting that for $d(x,y)<1$ we
have $x=y$. Condition E.6 follows trivially.
Finally, E.4 follows by noting that $p_\tau$ and $x_\tau$ are
piecewise constant mappings with a finite number of jumps. Thus, if
$\mathcal{U}$ does not include the indices corresponding to these
jumps, E.4 holds trivially.
\end{pf*}
\begin{pf*}{Proof of Theorem \ref{Thm:RateFronteir}}
For convenience, let $W_{x}=\{X\preccurlyeq x\}$.
Then, for all $x \in\mathcal{S}$ such that $p_x \geq\bar{p}$ we
have, by condition E.2,
\begin{eqnarray*}
|\hat\tau_x - \tau_x| & = & \biggl|\frac{\Pn(W_x) - \Pp
(W_x)}{p_x} + \Pn(W_x) \biggl( \frac{1}{\hat p_x} - \frac
{1}{p_x}\biggr) \biggr| \\
&=& \biggl|\frac{\Pn(W_x) - \Pp(W_x)}{p_x} +
\hat\tau_x \biggl( \frac{ p_x - \hat p_x}{p_x} \biggr)\biggr| \\
& \leq& \biggl|\frac{\Pn(W_x) - \Pp(W_x)}{p_x}\biggr| + \hat\tau
_x \biggl| \frac{ p_x - \hat p_x}{p_x} \biggr| \\
&\lesssim&\!\mbox{$_P$\,} \sqrt
{v(\bar{p})/n}.
\end{eqnarray*}
\upqed\end{pf*}
\begin{lemma}[(Technical lemma)]\label{Lemma:Tech}
Let $0 < \epsilon_1 \vee\epsilon_2 < \epsilon_3 < 1/2$ and $f, g,
h\dvtx\break[0, 1]\to[0,1]$, such that for all $t\in[0,1]$,
%
\begin{eqnarray}\label{Eq:ContJump}
\limsup_{t^k\to t} f(t^k) &\leq&
f(t) + \epsilon_1,\qquad \limsup_{t^k\to t} g(t^k) \leq g(t) +
\epsilon_1\quad
\mbox{and} \nonumber\\[-8pt]\\[-8pt]
\liminf_{t^k\to t} h(t^k) &\geq& h(t) - \epsilon_1.\nonumber
\end{eqnarray}
Moreover, assume that $\epsilon_2 < \epsilon_3\min_{t\in
[0,1]}h(t)$, and for every $t \in[\epsilon_3, 1-\epsilon_3]$:
\begin{longlist}
\item $|f(t) - th(t)|\leq\epsilon_2$,
\item $|g(t)- (1-t)h(t)| \leq
\epsilon_2$ and
\item $f(t) + g(t) \geq h(t)$.
\end{longlist}
Then, for every $\tau\in(3\epsilon_3, 1-3\varepsilon_3)$ there is
$\bar{t}$ such that $f(\bar{t}) \geq\tau h(\bar{t})-2\epsilon_1$
and $g(\bar{t}) \geq(1-\tau)h(\bar{t}) - 2\epsilon_1$.
\end{lemma}
\begin{pf}
Let $\bar{t} = \sup_{t \in[\epsilon_3, 1-\epsilon_3]} t \dvtx g(t)
\geq(1-\tau)h(t)$.
We have that $g(2\epsilon_3) \geq(1-2\epsilon_3)h(2\epsilon
_3)-\epsilon_2 = (1-\tau)h(2\epsilon_3) + (\tau-2\epsilon
_3)h(2\epsilon_3)-\epsilon_2 \geq(1-\tau)h(2\epsilon_3)$ by the
assumption on $\epsilon_2$ and $\tau$. Similarly, $g(1-2\epsilon_3)
\leq2\epsilon_3h(1-2\epsilon_3)-\epsilon_2<(1-\tau)h(2\epsilon
_3)$. So $\bar{t}\in[2\epsilon_3,1-2\epsilon_3]$.

Moreover, the condition (\ref{Eq:ContJump}) on $g$ and $h$ implies
that $g(\bar{t}) \geq(1-\tau)h(\bar{t})-2\epsilon_1$ and, by the
definition of $\bar{t}$, $g(\bar{t}+\mu) < (1 - \tau)h(\bar{t}+\mu
)$ for every $\mu> 0$. Thus, $f(\bar{t}+\mu) > \tau h(\bar{t}+\mu
)$ for every $\mu> 0$ by (iii). In turn, condition (\ref
{Eq:ContJump}) for $f$ and $h$ yields $f(\bar{t}) \geq\tau h(\bar
{t})-2\epsilon_1$, which establishes the result.
\end{pf}
\begin{pf*}{Proof of Theorem \ref{Thm:Rates1}}
The proof proceeds in steps. Step 1 establishes feasibility of a
``near'' partial quantile point. Step 2 derives the main arguments.
Step~3 concludes the proof.\vspace*{8pt}

\textit{Step} 1. \textit{Feasibility of near partial quantile
point}.\quad
Note that for any point $x$ that is feasible for (\ref
{Def:Qestimation}) we have $|\tau- \hat\tau_x| \leq\epsilon_n/\hat
p_x$. Moreover, by Theorem~\ref{Thm:RateFronteir}, if also $p_x \geq
\wp$, we have $|\hat\tau_{x} - \tau_{x}| \lesssim_P \sqrt{v(\wp
)/n}$, so that $|\tau- \tau_{x}| \lesssim_P u_n := \sqrt{v(\wp)/n}
+ \epsilon_n/\wp$.

Assume that $\epsilon_n \geq\epsilon_n^D$. Pick an arbitrary $x_\tau
\in\mathcal{Q}^*(\tau)$. By condition E.4, there is a continuous
path of quantile points, $\mathcal{P} = \{ x_{\tau'} \dvtx \tau' \in
(0,1) \}$, that passes through $x_\tau$.
Let $\epsilon_1 = \epsilon_n/2$, $\epsilon_2 = \sqrt{v(\wp)/n}$
and $\epsilon_3 = (1/6)\min_{u\in\mathcal{U}} u \wedge(1-u)$, so
that $f(t)= \Pn( X \preccurlyeq x_t )$,
$g(t)= \Pn( X \succcurlyeq x_t )$, and $h(t) = \hat p_{x_t}$ satisfies
condition (\ref{Eq:ContJump}), (i) and (ii) by Theorem \ref
{Thm:RateFronteir} and (iii) by definition. By Lemma~\ref
{Lemma:Tech}, there exists $x_{\tau^*} \in\mathcal{P}$ that is
feasible for (\ref{Def:Qestimation}). Since $p_{\tau^*} \geq\wp$,
we have $|\tau- \tau^*| \lesssim_P u_n$.
On the other hand, if $\epsilon_n \geq\epsilon_n^{D'}$, $x_\tau\in
\mathcal{Q}^*(\tau)$ is itself feasible with high probability. We can
take $x_{\tau^*}=x_\tau$ and the relation $|\tau- \tau^*| \lesssim
_P u_n$ would still hold.\vspace*{8pt}

\textit{Step} 2. \textit{Main argument.}\quad We will derive the rate of convergence by bounding
\[
p_{\tau_{\hat x_\tau}} - p_{\hat x_\tau} = p_{\tau_{\hat x_\tau}}
- p_\tau+ p_{x_\tau} - p_{\hat x_\tau}
\]
from above using E.5 and the optimality of $\hat x_\tau$, and from
below using the restricted identification condition E.2.

To establish the upper bound first note that by optimality of $\hat
x_\tau$, we have $\hat p_{x_{\tau^*}} \leq\hat p_{\hat x_{\tau}} $
and using E.5,
\begin{eqnarray*}
p_{x_{\tau}} - p_{\hat x_{\tau}} & \lesssim_P & \phi_n(d(\hat
x_\tau,x_\tau))/\sqrt{n} + \hat p_{x_{\tau}} -\hat p_{\hat x_{\tau
}}\\
& \lesssim_P & \phi_n(d(\hat x_\tau,x_\tau))/\sqrt{n} + \hat
p_{x_{\tau}} -\hat p_{x_{\tau^*}}.
\end{eqnarray*}
Applying E.5 one more time, and using that $|\tau^*-\tau|\lesssim_P
u_n$ so that $d(x_{\tau^*},\break x_\tau) \lesssim_P u_n$ and $p_{x_{\tau
}} - p_{x_{\tau^*}}\lesssim_P u_n^\gamma$,
\begin{eqnarray*}
p_{x_{\tau}} - p_{\hat x_{\tau}} & \lesssim_P & \phi_n(d(\hat
x_\tau,x_\tau))/\sqrt{n} + \phi_n(d(x_\tau,x_{\tau^*}))/\sqrt{n}
+ p_{x_{\tau}} - p_{x_{\tau^*}}\\
& \lesssim_P & \phi_n(d(\hat x_\tau,x_\tau))/\sqrt{n} + \phi
_n(u_n)/\sqrt{n} + u_n^\gamma.
\end{eqnarray*}
Also, since $|\tau_{\hat x_\tau} - \tau|\lesssim_P u_n$, by E.4,
$p_{\tau_{\hat x_\tau}} - p_\tau\lesssim_P u_n^\gamma$.

Note that if $d(\hat x_\tau,x_\tau) \lesssim_P u_n^\gamma$ we are
done. Therefore, the relations above yields
\[
p_{\tau_{\hat x_\tau}} - p_{\hat x_\tau} \lesssim_P \phi_n(d(\hat
x_\tau, x_\tau))/\sqrt{n} + u_n^\gamma.
\]
By E.3, we can minorate the left-hand side and obtain
\[
c \wedge\inf_{z \in\mathcal{Q}^*(\tau_{\hat x_\tau})} d(\hat
x_\tau, z)^\alpha\lesssim_P \phi_n(d(\hat x_\tau, x_\tau))/\sqrt
{n} + u_n^\gamma.
\]
Since the argument holds for all $x_\tau\in\mathcal{Q}^*(\tau)$, we have
\[
c \wedge\inf_{z \in\mathcal{Q}^*(\tau_{\hat x_\tau})} d(\hat
x_\tau, z)^\alpha\lesssim_P \phi_n\Bigl(\inf_{x_\tau\in\mathcal
{Q}^*(\tau)}d(\hat x_\tau, x_\tau)\Bigr)\big/\sqrt{n} + u_n^\gamma.
\]
Next note that the minimum in the left-hand side cannot be $c$ as $n$
grows [since $\phi_n(d(\hat x_\tau,x_\tau))$ can be bounded by
$2\sqrt{v(\wp/2)} = o(n^{1/2})$ by Theorem \ref{Thm:RateFronteir}].\vspace*{8pt}

\textit{Step} 3. \textit{Conclusion of the proof.}\quad Using that $\alpha\geq1$ by
E.3, E.4(ii), and the last relation in Step 2,
\begin{eqnarray*}
&& \inf_{x_{\tau}\in\mathcal{Q}^*(\tau)}d(\hat
{x}_{\tau},x_{\tau}) \\
&&\qquad \leq \inf_{x_{\tau}\in
\mathcal{Q}^*(\tau), z \in\mathcal{Q}^*(\tau_{\hat x_{\tau}})}
d(\hat{x}_{\tau},z) + d(z,x_{\tau}) \\
&&\qquad \lesssim  \inf_{z \in\mathcal{Q}^*(\tau_{\hat
x_{\tau}})} d(\hat{x}_{\tau},z) + |\tau- \tau_{\hat x_{\tau}}|\\
&&\qquad\lesssim_P \phi_n^{1/\alpha}\Bigl(\inf_{x_\tau
\in\mathcal{Q}^*(\tau)} d(\hat x_\tau,x_\tau)\Bigr)\big/n^{1/2\alpha
} + u_n^{\gamma/\alpha} + u_n
\\
&&\qquad\lesssim_P    u_n \vee u_n^{\gamma/\alpha}
\vee\phi_n^{1/\alpha}\Bigl(\inf_{x_\tau\in\mathcal{Q}^*(\tau)}
d(\hat x_\tau,x_\tau)\Bigr)\big/n^{1/2\alpha}.
\end{eqnarray*}
The rate result follows as in \cite{vdV-W}.
\end{pf*}
\begin{pf*}{Proof of Corollary \ref{Cor:CompleteOrder}}
Since the order is complete, $p_x = \hat p_x = 1$ for every $x \in
\mathcal{S}$. In particular, condition E.5 is satisfied with $\phi
_n(r) = 0$, E.4 with $\gamma= \alpha$ and E.3 with any positive
$\alpha$ since $\mathcal{Q}^*(\tau) = \mathcal{Q}(\tau)$.
In this case $\epsilon_n := \epsilon_n^D \wedge\epsilon_n^{D'} \leq
\epsilon_n^{D'} \lesssim_P \sqrt{v(1)/n}$.
\end{pf*}
\begin{pf*}{Proof of Theorem \ref{Thm:Inference}}
For convenience, let $\mathcal{W}_{x}=\{ X \preccurlyeq x\}$ and
$\mathcal{C}_{x}=%
\mathcal{C}(x)$. By E.6 we have $\sqrt{n}(\mathbb{P}_{n}(\mathcal
{W}_{x})-%
\Pp(\mathcal{W}_{x}))\rightsquigarrow N(0,\Pp(\mathcal{W}%
_{x})(1-\Pp(\mathcal{W}_{x})))$ and $\sqrt{n}(\hat p_x
-p_x)\rightsquigarrow N(0,p_x(1-p_x))$.

Moreover, we have%
\begin{eqnarray*}
\hat{\tau}_{x}-\tau_{x}
&=&  \frac{\mathbb{P}_{n}(\mathcal
{W}_{x})}{%
\hat p_x}-\frac{\Pp(\mathcal{W}_{x})}{p_x}
= \frac{\mathbb{P}_{n}(\mathcal{W}_{x})}{\hat p_x}%
-\frac{\mathbb{P}_{n}(\mathcal{W}_{x})}{p_x}+\frac{%
\mathbb{P}_{n}(\mathcal{W}_{x})-\Pp(\mathcal{W}_{x})}{p_x} \\
&=&  \mathbb{P}_{n}(\mathcal{W}_{x})\biggl( \frac{1}{\hat
p_x}-\frac{1}{p_x}\biggr) +\frac{\mathbb{P}_{n}(%
\mathcal{W}_{x})-\Pp(\mathcal{W}_{x})}{p_x}\\
&=& \frac{\mathbb{P}_{n}(\mathcal{W}_{x})}{\hat p_x}%
\frac{p_x-\hat p_x}{p_x}+\frac{\mathbb{P}_{n}(\mathcal{W}_{x})-\Pp
(\mathcal{%
W}_{x})}{p_x} \\
&=& -\hat{\tau}_{x}\frac{\hat p_x-p_x}{p_x}+\frac{\mathbb
{P}_{n}(\mathcal{%
W}_{x})-\Pp(\mathcal{W}_{x})}{p_x}\\
&=& (\tau_{x}-\hat{\tau}_{x})\frac{\hat p_x-%
p_x}{p_x}-\tau_{x}\frac{%
\hat p_x-p_x}{p_x}+\frac{\mathbb{P}_{n}(\mathcal{W}_{x})-\Pp
(\mathcal{W%
}_{x})}{p_x}.%
\end{eqnarray*}
By Condition E.2, $|\frac{\hat p_x-p_x}{p_x}| \lesssim_P
\sqrt{v(\bar{p})/n} = o_P(1)$, so that
\begin{eqnarray*}
\bigl(1+o_P(1)\bigr)p_x(\hat{\tau}_{x}-\tau_{x}) & = &
-\tau_{x}( \hat p_x-p_x) +\mathbb{P}_{n}(\mathcal
{W}_{x})-\Pp(\mathcal{W}_{x}) \\
& = & \frac{1}{\sqrt{n}}\mathbb{G}_{n}( 1\{\mathcal{W}_{x}\}
-\tau
_{x}1\{\mathcal{C}_{x}\}) .%
\end{eqnarray*}
Therefore, we have
$
p_x\sqrt{n}(\hat{\tau}_{x}-\tau_{x})=_P \mathbb{G}%
_{n}( 1\{\mathcal{W}_{x}\}-\tau_{x}1\{\mathcal{C}_{x}\}
) .
$ That converges to a zero mean normal distribution with variance
\begin{eqnarray*}
\mathbf{E}[ (1\{\mathcal{W}_{x}\}-\tau_{x}1\{\mathcal{C}_{x}\}
)^{2}%
] & = & \Pp(\mathcal{W}_{x})+\tau_{x}^{2}p_x-2\tau_{x}\Pp
(\mathcal{W}_{x}) \\
& = & \Pp(\mathcal{W}_{x})(1-\tau_{x})+\tau_{x}\bigl(\tau_{x}p_x-\Pp
(\mathcal{W}_{x})\bigr) \\
& = & \Pp(\mathcal{W}_{x})(1-\tau_{x})
\end{eqnarray*}
using $\mathcal{W}_{x}\subseteq\mathcal{C}_{x}$ and $\tau_{x}=\Pp(%
\mathcal{W}_{x})/p_x$. Finally, we get
\[
\sqrt{n}(\hat{\tau}_{x}-\tau_{x})\rightsquigarrow N\biggl( 0,\frac
{\tau
_{x}(1-\tau_{x})}{p_x}\biggr) .
\]

Note that within $\mathcal{C}_{\bar{p}}$, all the functions are
bounded by $2/\bar{p}$
with high probability for large enough sample size. Therefore, a
multidimensional central limit theorem applies and the covariance structure
of a pair $x,y\in\mathcal{S}$ is given by
\[
\Omega_{x,y}=\mathbf{E}\biggl[ \frac{(1\{\mathcal{W}_{x}\}-\tau
_{x}1\{%
\mathcal{C}_{x}\})}{p_x}\frac{(1\{\mathcal{W}%
_{y}\}-\tau_{y}1\{\mathcal{C}_{y}\})}{p_y}\biggr] .
\]
After simplification, we obtain
\begin{eqnarray*}
\Omega_{x,y} & = &  \frac{\Pp(\mathcal{W}_{x}\cap
\mathcal{W}_{y})}{\Pp(\mathcal{C}_{x})\Pp(\mathcal{C}_{y})}%
-\tau_{x}\frac{\Pp(\mathcal{C}_{x}\cap\mathcal{W}_{y})}{\Pp(%
\mathcal{C}_{x})\Pp(\mathcal{C}_{y})}-\tau_{y}\frac{\Pp(%
\mathcal{C}_{y}\cap\mathcal{W}_{x})}{\Pp(\mathcal{C}_{x})\Pp(%
\mathcal{C}_{y})}+\tau_{x}\tau_{y}\frac{\Pp(\mathcal{C}_{x}\cap
\mathcal{C}_{y})}{\Pp(\mathcal{C}_{x})\Pp(\mathcal{C}_{y})} \\
& = &  \tau_{x}\tau_{y}\biggl( \frac{\Pp(\mathcal{W}%
_{x}\cap\mathcal{W}_{y})}{\Pp(\mathcal{W}_{x})\Pp(\mathcal{W}%
_{y})}-\frac{\Pp(\mathcal{C}_{x}\cap\mathcal{W}_{y})}{p_x\Pp
(\mathcal{W}_{y})}-\frac{\Pp(\mathcal{W}%
_{x}\cap\mathcal{C}_{y})}{\Pp(\mathcal{W}_{x})p_y}+\frac{\Pp
(\mathcal{C}_{x}\cap\mathcal{C}_{y})}{p_xp_y}\biggr) .%
\end{eqnarray*}
Finally, asymptotic equicontinuity of $\beta_n(x)$ follows directly
from the asymptotic equicontinuity of $\alpha_n(x)$ implied by E.6 and
$\bar{p}>0$ being fixed.
\end{pf*}
\begin{pf*}{Proof of Corollary \ref{Cor:AsymptoticTauPQP}}
The proof of the second result builds upon arguments in \cite
{EinmahlMason1992,BucchianicoEinmahlMushkudiani2001}. Based on Theorem
\ref{Thm:Inference}, we have that for $\mathcal{C}_{\wp/2} = \{ x
\in\mathcal{S}, p_x \geq\wp/2\}$, the process $\beta_n(x):=\sqrt
{n}(\hat\tau_x - \tau_x)$
converges weakly in $\ell^{\infty}(\mathcal{C}_{\wp/2})$ to a
bounded, mean zero Gaussian process $G_P$. By the
Skorohod--Dudley--Wichura representation theorem, there exists a
probability space
$(\widetilde\Omega, \mathcal{\widetilde A}, \widetilde P)$ carrying
versions $\widetilde G_P$ and $\widetilde\beta_n$ of $G_P$ and $\beta
_n$ such that
$\sup_{x\in\mathcal{C}_{\wp/2}}| \widetilde\beta_n(x)-\widetilde
G_P(x)| \to0$ as $n$ grows. Next, note that for all $\tau\in\mathcal
{U}$, $\hat x_{\tau} \in\mathcal{C}_{\wp/2}$ provided that $\sqrt
{v(\wp)/n} = o(\wp)$. Thus,
\[
\sqrt{n}(\tau_{\hat x_\tau} - \tau) = -\widetilde\beta_n(\hat
x_\tau) + \sqrt{n}( \hat\tau_{\hat x_{\tau}} - \tau) = o(1) +
\widetilde G_P(\hat x_\tau) + \sqrt{n}( \hat\tau_{\hat x_{\tau}} -
\tau).\quad
\]
\upqed\end{pf*}
\begin{pf*}{Proof of Theorem \ref{Thm:p*}}
Let $\tau^*$ and $\hat\tau^*$ be such that $\wp= p_{\tau^*}$ and
$\hat\wp= \hat p_{\hat\tau^*}$. Thus, we have $\hat x_{\hat\tau
^*}$ and $\hat x_{\tau^*}$ satisfying $\hat p_{\hat\tau^*} = \hat
p_{\hat x_{\hat\tau^*}}$ and $ \hat p_{\tau^*} = \hat p_{\hat
x_{\tau^*}}$. Moreover, let $u_n := \sqrt{v(\wp)/n} + \epsilon
_n/\wp\lesssim n^{-1/2}$ by assumption.

First, note that since $\hat\wp\leq\hat p_{\tau^*}$, and $p_{\tau
_{\hat x_{\tau^*}}} \geq p_{\hat x_{\tau^*}}$, we have, by E.5,
\begin{eqnarray*}
\hat\wp- \wp& \leq& \hat p_{\tau^*} - p_{\tau^*} = \hat p_{\hat
x_{\tau^*}} - p_{x_{\tau^*}} \\
& = & \hat p_{\hat x_{\tau^*}} - p_{\hat x_{\tau^*}} - ( \hat
p_{x_{\tau^*}} - p_{x_{\tau^*}}) + p_{\hat x_{\tau^*}} - p_{x_{\tau
^*}} + \hat p_{x_{\tau^*}} - p_{x_{\tau^*}} \\
& \lesssim&\!\mbox{$_P$\,}  \phi_n(d(\hat x_{\tau^*}, x_{\tau^*})) /\sqrt{n} +
p_{\tau_{\hat x_{\tau^*}}} - p_{x_{\tau^*}} + \hat p_{x_{\tau^*}} -
p_{x_{\tau^*}}.
\end{eqnarray*}
Note also that by Step 1 in the proof of Theorem \ref{Thm:Rates1} we
have $|\tau_{\hat x_{\tau^*}} - \tau^*| \lesssim_P u_n$. Moreover,
$p_\tau$ is locally quadratic around $\tau^*$. Therefore,
\[
\hat\wp- \wp \lesssim_P  \phi_n(d(\hat x_{\tau^*}, x_{\tau
^*})) /\sqrt{n} + u_n^2 + \hat p_{x_{\tau^*}} - p_{x_{\tau^*}}.
\]
Since it holds for any $x_{\tau^*}\in\mathcal{Q}^*(\tau^*)$,
%
\begin{eqnarray}\label{Eq:UBp*}
\hat\wp- \wp& \lesssim_P & \phi_n\Bigl(\inf_{x_{\tau^*}\in\mathcal
{Q}^*(\tau^*)}d(\hat x_{\tau^*}, x_{\tau^*})\Bigr)\big/\sqrt{n}\nonumber\\
&&{} + u_n +
\max_{x_{\tau^*}\in\mathcal{Q}^*(\tau^*)}\{\hat p_{x_{\tau^*}} -
p_{x_{\tau^*}}\}\\
& \lesssim_P & o(n^{-1/2}) + \max_{x_{\tau^*}\in\mathcal{Q}^*(\tau
^*)}\{\hat p_{x_{\tau^*}} - p_{x_{\tau^*}}\}\nonumber
\end{eqnarray}
since $u_n^2 = o(n^{-1/2})$, and $\inf_{x_{\tau^*}\in\mathcal
{Q}^*(\tau^*)}d(\hat x_{ \tau^*}, x_{ \tau^*})=o_P(1)$ by Theorem
\ref{Thm:Rates1}.

Next, by Step 1 in the proof of Theorem \ref{Thm:Rates1}, for every
$x_{\hat\tau^*}$ there is a partial quantile point $x_{\bar\tau}$,
$d(x_{\hat\tau^*},x_{\bar\tau})\lesssim u_n$ that is feasible for
(\ref{Def:Qestimation}) with $\hat\tau^*$. Thus, $\hat p_{\hat
x_{\hat\tau^*}} \geq\hat p_{x_{\bar\tau}}$.
Using this inequality, E.5, and that $p_{\hat\tau^*} \geq p_{\tau
^*}$ by definition (\ref{Def:p*}),
%
\begin{eqnarray}\label{Eq:LBp*}
\hat\wp- \wp&\geq& \hat p_{x_{\bar\tau}} - p_{x_{\tau^*}}\nonumber\\
& = & \hat p_{x_{\bar\tau}} - p_{ x_{\bar\tau}} - ( \hat p_{
x_{\tau^*}} - p_{x_{\tau^*}} ) + p_{x_{\bar\tau}} - p_{x_{\tau^*}}
+ \hat p_{x_{\tau^*}} - p_{x_{\tau^*}}\nonumber\\[-8pt]\\[-8pt]
& \gtrsim&\!\mbox{$_P$\,} - \phi_n( d(x_{\bar\tau}, x_{\tau^*}) ) /\sqrt{n} +
p_{x_{\hat\tau^*}} - p_{x_{\tau^*}} + \hat p_{x_{\tau^*}} -
p_{x_{\tau^*}}\nonumber\\
& \geq& - \phi_n( d(x_{\bar\tau}, x_{\tau^*})) /\sqrt{n} + \hat
p_{x_{\tau^*}} - p_{x_{\tau^*}},\nonumber
\end{eqnarray}
where $x_{\bar\tau}$ was chosen to be close to $x_{\tau^*}$, namely
$d(x_{\bar\tau},x_{\tau^*}) \leq d(x_{\bar\tau},x_{\hat\tau^*})
+d(x_{\hat\tau^*},x_{\tau^*}) \lesssim_P u_n + |\hat\tau^*-\tau^*|$.
Therefore, (\ref{Eq:LBp*}) holds for any $x_{\tau^*} \in\mathcal
{Q}^*(\tau^*)$ and
$d(x_{\hat\tau^*}, x_{\tau^*}) \lesssim|\hat\tau^* - \tau^*| =
o_P(1)$ by Lemma \ref{Lemma:Aux01} below. Thus,
%
\begin{equation}\label{Eq:LBp*2}
\hat\wp- \wp \geq - o_P(n^{-1/2}) + \max_{x_{\tau^*}\in
\mathcal{Q}^*(\tau^*)}\{ \hat p_{x_{\tau^*}} - p_{x_{\tau^*}}\}.
\end{equation}
Combining (\ref{Eq:LBp*2}) and (\ref{Eq:UBp*}), we obtain $ \sqrt
{n}(\hat\wp- \wp) = o_P(1) + Z_P(\tau^*)$.
\end{pf*}
\begin{lemma}\label{Lemma:Aux01}
Under the assumptions of Theorem \ref{Thm:Rates1}, and that $\tau
\mapsto p_\tau$ is a twice differentiable function, let $\wp= p_{\tau
^*}$ and $\hat\wp= \hat p_{\hat\tau^*}$. Then $|\hat\tau^* - \tau
^*| = o_P(1)$.
\end{lemma}
\begin{pf}
Consider the twice differentiable function $\tau\mapsto p_\tau$.
Since $p_{\tau^*}$ is its strict minimum at the interior of $\mathcal
{U}$, we have
$p_{\tau} - p_{\tau^*} \gtrsim|\tau- \tau^*|^2$ for
$\tau\in\mathcal{U}$. 

By Step 1 of the proof of Theorem \ref{Thm:Rates1}, for every $\tau
\in\mathcal{U}$ we have that there is an $x_{\bar\tau}$ that is
feasible and $|\bar\tau- \tau| \lesssim_P u_n=o_P(1)$. Thus,
\[
\hat p_\tau= \hat p_{\hat x_\tau} \geq\hat p_{x_{\bar\tau}}
\gtrsim_P p_{\bar\tau} - \sqrt{v(\wp)/n} \gtrsim_P p_\tau- \sqrt
{v(\wp)/n} - u_n^\gamma= o_P(1) + p_\tau.
\]
Similarly, since $|\tau_{x_{\bar\tau}}-\tau|\lesssim_P u_n$,
\begin{eqnarray*}
\hat p_\tau&\lesssim_P& p_{\hat x_\tau} + \sqrt{v(\wp/2)/n} \lesssim
_P p_{\tau_{x_{\bar\tau}}} + \sqrt{v(\wp/2)/n} \\
&\lesssim_P& p_{\tau} + \sqrt{v(\wp/2)/n} + u_n^\gamma= o_P(1) + p_\tau.
\end{eqnarray*}
Therefore, using that $\hat p_{\hat\tau^*}\leq\hat p_{\tau^*}$,
\[
|\hat\tau- \tau^*|^2 \leq p_{\hat\tau^*} - p_{\tau^*} = o_P(1) +
v - \hat p_{\tau^*} = o_P(1).
\]
\upqed\end{pf}

\vspace*{-8pt}

\section{\texorpdfstring{Section \lowercase{\protect\ref{sec:add}} proofs}{Apendix C: Section 4 proofs}}

\vspace*{-8pt}

\begin{pf*}{Proof of Lemma \ref{Lemma:CHARACTERIZATIONcone}}
This follows if $\operatorname{support}   \hat1_K = \RR^d$, where $\hat
1_K$ is the Fourier transform of the indicator function of the set $K$,
see \cite{BagchiSitaram1982}, Proposition~3.1. (We proceed as in
Proposition~3.2 in \cite{BagchiSitaram1982} with the necessary
modifications.)\vspace*{8pt}

\textit{Step} 1. Let $0 \neq f \in L^1(\RR^d)$ such that $\operatorname{support}
f \subseteq K$, $\hat f(w) = \int_{\RR^d} e^{-iw'x}\times f(x)\,dx=\int
_{K} e^{-iw'x}f(x)\,dx$, and $K^o = \{ y \in\RR^d \dvtx y'x \leq0
$ for all $x \in K\}$ denote the polar cone of $K$. Define the
regions (in the complex space $\mathbb{C}^d$)
\[
H = \{ z \in\mathbb{C}^d \dvtx  \operatorname{Im}(z) \in K^o \}
\quad\mbox{and}\quad
H_0 = \{ z \in\mathbb{C}^d \dvtx \operatorname{Im}(z) \in\operatorname{int}  K^o \}.
\]
It follows from the definition that $\hat f$ can be extended to a
bounded function~$g$ in the region $H$ [because $K$ is a proper convex
cone, for any $w \in H_0$ and $x \in K$, we have $\operatorname{Re}(-iw'x) \leq
0$]. Moreover, $g$ is analytic in $H_0$ and continuous in $H$.
Therefore, $\hat f$ is the restriction of the bounded analytic function
$g$ on the boundary of $H$ \cite{Bochner1937}. Consequently,\vspace*{1pt} $\hat f$
cannot be identically zero on an open subset of $\RR^d$ (which would
imply that $\hat f = 0$ and, thus, $f=0$), equivalently, $\operatorname{support}
\hat f = \RR^d$.\vspace*{8pt}

\textit{Step} 2. Next, we consider $1_K$ which is a nonzero bounded Borel
function which is not in $L^1(\RR^d)$. By contradiction, assume that
$\hat1_K$ vanishes on a~nonempty open set $U$ of $\RR^d$, that is,
$(\operatorname{support}  \hat1_K)\cap U = \varnothing$. Let $x_0$ and
$\varepsilon>0$ such that $B(x_0,2\varepsilon) \subset U$.

Let $0\neq h_1 \in L^1(\RR^d)$ such that $\hat h_1$ is a $C^{\infty}$
function and $\operatorname{support} \hat h_1 \subset B(0,\varepsilon)$. Then
\[
\mathrm{support} \widehat{(h_1 \cdot1_K)} = \operatorname{support}(\hat h_1 *
\hat1_K) \subseteq\operatorname{support}  \hat h_1 + \operatorname{support}  \hat
1_K \subset\RR^d \setminus B(x_0,\varepsilon),
\]
where ``$*$'' denotes the convolution operator.

On the other hand, $h_1 \cdot1_K \in L^1(\RR^d)$ with $\operatorname{support}
(h_1\cdot1_K) \subset K$. Therefore, by Step 1, $h_1\cdot1_K = 0$
almost everywhere on $\RR^d$. In turn, $\hat h_1$ is a $C^\infty
$-function of compact support, so $h_1$ is the restriction of an entire
function to $\RR^d$, and hence $h_1(x)\neq0$ almost everywhere in
$\RR^d$. Thus, $1_K$ is zero almost everywhere which give us a
contradiction since $K$ is a proper convex cone.
\end{pf*}
\begin{pf*}{Proof of Lemma \ref{Lemma:CHARACTERIZATIONgraph}}
Without loss of generality, we can consider only connected graphs
(otherwise we proceed with each connected component separately).
We provide an algorithm.

For each node, we have $\tau_x p_x = P(X \preccurlyeq x)$. If there is
no incoming arc on~$x$, we have that $P(X \preccurlyeq x) = P( X = x)$.
For a general node $x$, if we already computed $P(X=y)$ for all $y\neq
x$, $y \preccurlyeq x$, then we have $P(X=x) = \tau_xp_x - \sum
_{y\neq x, y \preccurlyeq x} P(X=y)$. Otherwise, ``backtrack'' to
consider a $y\neq x$, $y \preccurlyeq x$ for which $P(X=y)$ is not
known. Since there are no cycles, we can only ``backtrack'' at most
$|\mS|<\infty$ before computing a probability for some $y$. Thus the
procedure terminates in a finite number of steps with all proba\-bilities.
\end{pf*}
\begin{pf*}{Proof of Theorem \ref{Thm:Rearrangement}}
The proof follows from the inequality of Lo\-rentz~\cite{Lorentz1953}
applied to each component individually. This follows the strategy of
Chernozhukov, Fern\'{a}ndez-Val and Galichon \cite{CFG2009} that
previously used this inequality to prove a similar result.
\end{pf*}
\begin{pf*}{Proof of Corollary \ref{Corollary:Rearrangement}}
If $x_\tau$ is partial-monotone, by Theorem \ref{Thm:Rearrangement}
we have
\begin{eqnarray*}
\biggl|\int_0^1\|\hat x^r_{u}-\hat x_{u}\|^{\kappa
}\,du\biggr|^{1/\kappa} & \leq&  \biggl| \int_0^1\|\hat
x_{u}^r- x_{u}\|^{\kappa}\,du\biggr|^{1/\kappa} + \biggl|\int_0^1\|
\hat x_{u}- x_{u}\|^{\kappa}\,du \biggr|^{1/\kappa}\\
& \leq&   2 \biggl|\int_0^1\|\hat x_{u}- x_{u}\|^{\kappa
}\,du\biggr|^{1/\kappa}.
\end{eqnarray*}
The second follows by a triangular inequality.
\end{pf*}
\begin{pf*}{Proof of Theorem \ref{Thm:Ind}}
%
Note that by independence and no point mass, we have
$P(X\succcurlyeq x)=\prod_{j=1}^{d}(1-F_{j}(x_{j}))$, $P(X\preccurlyeq
x)=\prod_{j=1}^{d}F_j(x_{j})$ and $p_x = P(X\succcurlyeq
x)+P(X\preccurlyeq x)$. Thus,
$x_\tau\in\arg\max\{\prod_{j=1}^{d}(1-F_{j}(x_{j})) + \prod
_{j=1}^{d}F_{j}(x_{j}) \dvtx\break \tau\prod_{j=1}^{d}(1-F_{j}(x_{j})) =
(1-\tau)\prod_{j=1}^{d}F_j(x_{j})\}$. By the independence, we can
write $a_j = F_j(x_j)$ and recast the problem as
$ \max_a\{ \prod_{j=1}^d a_j + \prod_{j=1}^d(1-a_j) \dvtx (1-\tau)\prod
_{j=1}^d a_j = (1-\tau)\prod_{j=1}^d (1-a_j), 0\leq a_j \leq1\}$. By
inspection, we have that $0<a_j<1$, $j=1,\ldots,d$, at the optimal. By
the optimality conditions, there is a $\lambda$ such that we have for
every $k=1,\ldots,d$
\[
0 = \prod_{j\neq k} a_j - \prod_{j\neq k} (1-a_j) - \lambda(1-\tau)
\prod_{j\neq k} a_j + \lambda\tau\prod_{j\neq k} (1-a_j).
\]
This implies that for every $j=1,\ldots,d$, we have $
\frac{1-\lambda\tau}{1-\lambda(1-\tau)} =\frac{\tau}{1-\tau
}\cdot\frac{1-a_k}{a_k}$.

Therefore, $a_k^*=a(\tau)$ for every $k=1,\ldots,d$. On the other
hand, by feasibility we must\vspace*{1pt} have $\prod_{j=1}^d [a_j^*/(1-a_j^*)] =
a(\tau)^d / (1-a(\tau))^d = \tau/(1-\tau)$. Therefore, $a(\tau
)/(1-a(\tau)) = (\tau/(1-\tau))^{1/d}$, which yields the result.
\end{pf*}
\begin{pf*}{Proof of Theorem \ref{Thm:IndInd}}
By Proposition \ref{Prop:Invariance} with $h(x) =
(F_1(x_1),F_2(x_2),\ldots,\break F_d(x_d))$, we have $\tau_{X} = \tau
_{h(X)}$ so that we can assume that $X$ is a uniform $(0,1)$ random
variable. Therefore,
\[
P( \tau_X \leq\tau) = P\Biggl( \prod_{j=1}^d x_j \leq \tau\Biggl[\prod
_{j=1}^d x_j + \prod_{j=1}^d (1-x_j)\Biggr] \Biggr) = P\Biggl( \prod_{j=1}^d \frac
{x_j}{1-x_j} \leq \frac{\tau}{1-\tau} \Biggr).
\]
The first result follows by taking logs and noting that $Z_j:=\log
(x_j/(1-x_j))$ is distributed as a logistic random variable with zero
mean and variance $\pi^2/3$ when $x_j$ is a uniform $(0,1)$ random variable.

Next, since $Z_j$ is symmetric around zero, $P(\tau_X \geq1/2) =
P(\sum_{j=1}^d Z_j \geq0) = 1/2$. Finally, let $Z^{(d)}:=d^{-1/2}\sum
_{j=1}^d Z_j$ and denote its probability density function by $f_d$. It
follows that $\max_z f_d(z) = f_d(0) \leq1/2$. Since $Z^{(d)}$ is
symmetric, we have, for $t \in(0,1/2)$,
\[
P( |\tau_X - 0.5 | \geq t ) = 2P( \tau_X \geq0.5 + t ) = 2P\biggl( Z^{(d)}
\geq d^{-1/2}\log\biggl(\frac{0.5 + t}{0.5 - t}\biggr) \biggr).
\]
Thus, using that $\log(1+x)\leq x$ and $f_d(z)\leq1/2$,
\begin{eqnarray*}
P( |\tau_X - 0.5 | \geq t ) &\geq& 2P\biggl( Z^{(d)} \geq\frac
{2td^{-1/2}}{0.5-t} \biggr) \geq 1 - 2\int_0^{2td^{-1/2}/(0.5-t)}
f_d(z)\,dz\\
&\geq& 1 -  \frac{2td^{-1/2}}{0.5-t}.
\end{eqnarray*}
Using $t:= 0.5-Cd^{-1/2}$ in the expression above,
\[
P( |\tau_X - 0.5 | \geq0.5-Cd^{-1/2} ) \geq1 -  \frac
{2(0.5-Cd^{-1/2})d^{-1/2}}{Cd^{-1/2}} \geq1 - 1/C.\quad
\]
\upqed\end{pf*}
\begin{pf*}{Proof of Lemma \ref{Lemma:Finite}}
We can compute the partial order and the probabilities $\Pp%
(X\succcurlyeq x|\mathcal{C}(x))$ and $\Pp(X\preccurlyeq x|\mathcal
{C}%
(x))$, which are bounded by $O(|\mathcal{S}|)$ for every fixed $x\in
\mathcal{S}$. Varying over all choices of $|\mathcal{S}|$, we obtain
$O(|%
\mathcal{S}|^{2})$ operations.
\end{pf*}
\begin{pf*}{Proof of Lemma \ref{Lemma:LogConcaveConvexCone}}
First, note that under C.2, we have that $K$ is a convex cone with
nonempty interior. Therefore, $K$ has a strict recession direction,
that is, $\exists w\neq0$ such that $K+w\subset
\operatorname{int} K$. Moreover, if $K\cap-K$ is full dimensional, $K=\RR
^d$ and we have $x \succcurlyeq y$ for every $x, y \in\RR^d$ and the
result holds trivially. Therefore, we can assume that $K\cap-K$ is not
full dimensional.

Since $K\cap-K$ is not full dimensional and $X$ has no point mass, we
have $\Pp(X\succcurlyeq x\succcurlyeq X)=0$ for every $x \in\mathcal
{S}$. Therefore $p_x=\Pp(X\succcurlyeq x)+\Pp%
(X\preccurlyeq x)$.
Moreover, $p_x$, $\Pp(X\preccurlyeq x)$
and $\Pp(X\succcurlyeq x)$ are continuous in $x$.

Note that any pair $(p,x)$ such that $x\in\mathcal{Q}(\tau)$ and
$p=p_x$ is feasible for problem (\ref{Prob:Aux}). By the log-concavity
of the probability density function, $\Pp(X\preccurlyeq x)=\Pp(x-K)$
and $\Pp(X\succcurlyeq x)=\Pp(x+K)$ are log-concave functions of $x$
by the
Pr\'{e}kopa--Leindler inequality (e.g., see \cite{Gadner2002}). This shows
that (\ref{Prob:Aux}) can be recast as a convex programming problem.

Next, we will show that the solution to (\ref{Prob:Aux}) also solves
(\ref%
{Def:Beta-Prob}). If $p^{\ast}=p_{x^{\ast}}$, then
both constraints are active at the optimal point, and the result follows.
Note that at least one constraint must be active at $(p^{\ast},x^{\ast})$.

Suppose $p^{\ast}<p_{x^{\ast}}$, in which case $%
x^{\ast}\notin\mathcal{Q}(\tau)$. Without loss of generality,
assume that $\Pp%
(X\succcurlyeq x^*)>(1-\tau)p^{\ast}$. Define the continuous
functions $u(t)=%
\Pp(X\succcurlyeq x^{\ast}+td)$ and $\ell(t)=\Pp%
(X\preccurlyeq x^{\ast}+td)$, which are, respectively, decreasing and
increasing in $t$. For some $t>0$, we have $u(t)>(1-\tau)p^{\ast}$
and $%
\ell(t)>\tau p^{\ast}$, which contradicts the optimality of $(p^{\ast
},x^{\ast})$.
\end{pf*}
\begin{pf*}{Proof of Theorem \ref{Thm:Comp}}
From Lemma \ref{Lemma:LogConcaveConvexCone}, it follows that we can recast
the problem as the convex programming problem defined in (\ref{Prob:Aux}).
For $\bar{p} < p_\tau$, define the convex set
\begin{eqnarray*}
H(\bar{p})&:=&\{(v,x) \in\RR\times\mS:
\log\Pp(X\succcurlyeq x)\geq\log(1-\tau) + v, \\
&&\hspace*{39.5pt}\log\Pp(X\preccurlyeq x)\geq\log\tau+ v,
\log\bar{p} \leq v \leq0\},
\end{eqnarray*}
where $v = \log p$ for $p$ in (\ref{Prob:Aux}). For an arbitrary
$\varepsilon> 0$, note that for every $x$ we can approximate $\Pp
(X\succcurlyeq x)$ and $\Pp(X\preccurlyeq x)$ up to a multiplicative
factor of $1+\varepsilon$ using the integration procedure for
log-concave distributions based on random walks proposed by Lov\'asz
and Vempala \cite{LV2006}.
Relying on these results, we can construct and $\varepsilon
_0$-approximate a membership oracle whose complexity is given by
\[
O\biggl( \frac{d^{4}\log^{3}d \log(1/\delta)}{\varepsilon_0
^{2}}\biggr),
\]
where $\varepsilon_0 = p_\tau\varepsilon$.
Note that by controlling the error in the computation of $\Pp
(X\preccurlyeq x)$ and $\Pp(X\succcurlyeq x)$ by a factor of
$1+\varepsilon$, we control the error in the computation of $\tau_x$
by an additive error of $\varepsilon$.

Based on this membership oracle, we can apply the results in \cite
{LV2006} for optimization, which requires $O^{\ast}(d^{4.5})$ calls of
the constructed membership oracle.
\end{pf*}
\end{appendix}

\section*{Acknowledgments}
We would like to thank Victor Chernozhukov, Gustavo Didier and Ilia
Tsetlin for helpful comments that improved the paper. We also thank
Saed Alizamir for thorough proofreading of several versions of this
paper. We are grateful to two anonymous referees and the editors for
comments and suggestions that increased the quality of the paper.



%
\printaddresses

\end{document}